\documentclass[12pt]{amsart}

\usepackage{amsmath}
\usepackage{amssymb}

\newtheorem{theorem}{Theorem}[section]

\newtheorem{lemma}[theorem]{Lemma}
\newtheorem{proposition}[theorem]{Proposition}
\newtheorem{corollary}[theorem]{Corollary}

\theoremstyle{definition}

\theoremstyle{remark}

\newcommand{\Proof}{\begin{proof}} 

\def\dots{\mathinner{\ldotp\ldotp\ldotp}}

\def\R{\mathbb R}

\def\C{\mathbb C}
\def\D{\mathcal D}

\def\P{\mathcal P}

\numberwithin{equation}{section}
\allowdisplaybreaks
\begin{document}

\title[Infinitesimal form boundedness 
and Trudinger's subordination]{ 
Infinitesimal form boundedness \\
and Trudinger's subordination\\ for the Schr\"odinger operator}

\author[V.~G. Maz'ya]
{V.~G. Maz'ya}
\address{Department of Mathematics,
Ohio State University, 231 West 18-th
Ave.,  Columbus, OH 43210, USA}
\address{Department of Mathematics,
Link\"oping University, SE-581 83, Link\"oping, 
Sweden}
\email{vlmaz@mai.liu.se}

\author[I.~E. Verbitsky]
{I.~E. Verbitsky$^*$}
\address{Department of Mathematics,
University of Missouri,
Columbia, MO 65211, USA}
\email{igor@math.missouri.edu}

\thanks{$^*$Supported in part 
by NSF Grant DMS-0070623.}

\begin{abstract} 
 We give explicit analytic criteria for two problems 
 associated with the  Schr\"odinger 
operator 
$H = -\Delta + Q$  on $L^2(\R^n)$ where $Q\in D'(\R^n)$ is   
an arbitrary real- or complex-valued potential. 

First, we  obtain necessary and 
sufficient conditions on   $Q$ so that the quadratic form 
$\langle Q \cdot, \, \cdot\rangle$ 
 has zero relative bound with respect to the Laplacian. 
 For $Q\in L^1_{\rm loc}(\R^n)$, this property can be expressed 
in the form of the integral inequality:
$$
 \left \vert \int_{\R^n}
|u(x)|^2 \, Q(x) \, dx\right  \vert \leq \epsilon \, ||\nabla
u||^2_{L^2(\R^n)} + C(\epsilon) \, ||u||^2_{L^2(\R^n)}, \quad \forall u  \in
C^\infty_0(\R^n), 
$$
for an arbitrarily small $\epsilon >0$ and some $C(\epsilon)> 0$. 
One of the major steps here is the reduction to a similar inequality 
with nonnegative  function 
$|\nabla (1-\Delta)^{-1} \, Q|^2 + |(1-\Delta)^{-1} \, Q|$ 
in place of $Q$. This  
 provides a complete solution to the infinitesimal form boundedness 
problem for the Schr\"odinger operator, and leads to new broad classes 
of admissible distributional potentials $Q$, which extend the usual 
$L^p$ and Kato classes, as well as those based on the well-known conditions 
of Fefferman--Phong and Chang--Wilson--Wolff.

Secondly, we characterize Trudinger's subordination property  
where $C(\epsilon)$ in the above inequality 
is subject to the condition 
 $C(\epsilon) \le  c \, {\epsilon^{-\beta}}$ ($\beta>0$) 
as $\epsilon\to +0$.  Such quadratic form inequalities 
can be understood entirely 
 in the framework of Morrey--Campanato spaces, 
using  mean oscillations of 
$\nabla (1-\Delta)^{-1} \, Q$ and $(1-\Delta)^{-1} \, Q$ on balls or cubes. 
A version of this  condition where 
$\epsilon\in(0,  +\infty)$ is equivalent to 
the multiplicative inequality:
$$
 \left \vert \int_{\R^n}
|u(x)|^2 \, Q(x) \, dx\right  \vert \leq C \, 
||\nabla u||^{2p}_{L^2(\R^n)}  \,||u||^{2(1-p)}_{L^2(\R^n)}, 
\quad \forall u 
\in C^\infty_0(\R^n), 
$$
with $p=\frac \beta {1 + \beta} \in (0, \, 1)$. We show that this 
inequality 
 holds if and only if 
$\nabla \Delta^{-1} Q \in \text{BMO}(\R^n)$ if
$p=\frac 1 2$.  For $0<p<\frac 1 2$, it is 
valid whenever $\nabla \Delta^{-1} Q$ is  H\"older-continuous 
of order $1-2p$, or respectively lies in the  Morrey space 
$\mathcal {L}^{2, \lambda}$ with $\lambda = n+2 - 4p$ if $\frac 1 2 <p<1$. 
As a consequence, we characterize completely the class of those $Q$ which 
satisfy an analogous multiplicative inequality of 
Nash's type, with $||u||_{L^1(\R^n)}$ in place of $||u||_{L^2(\R^n)}$. 

These results are intimately connected with spectral theory 
and dynamics of the Schr\"odinger operator, and elliptic PDE theory. 

\end{abstract}

\vfill
\eject

\maketitle

\vfill
\eject

\section{Introduction}\label{Introduction}

In the present paper we obtain necessary and sufficient conditions 
for the  {\em infinitesimal form
boundedness\/}  of the potential energy operator $Q$ with respect to the 
 kinetic energy operator 
$H_0 = -\Delta$ on $L^2(\R^n)$.  
Here $Q$ is an  arbitrary 
real- or
complex-valued  potential (possibly a distribution). 
This notion appeared in relation to the so-called KLMN theorem 
(\cite{RSi}, Theorem X.17), and has become  an indispensable tool in mathematical 
quantum mechanics  and PDE  theory.  
Furthermore, we characterize explicitly a 
related  {\em form subordination property} of Trudinger type. (See \cite{Ka2}, 
\cite{RSi},  \cite{RSS}, 
\cite{Sch}, \cite{Tru}, and the literature cited there.) 

Both of these notions play a major role in 
studies of asymptotics of the spectral counting function
for self-adjoint and non-self-adjoint  Schr\"odinger operators 
 \cite{BiSo}, \cite{Gr1}, \cite{Gr2}, \cite{MM1}, \cite{MM2}, 
generalized eigenvector expansions \cite{Agr}, \cite{RSS}, 
Schr\"odinger semigroups \cite{Dav1}, \cite{Dav2}, \cite{LPS},  \cite{Sim}, 
stochastic processes and  elliptic PDE \cite{AiSi}, \cite{CZh}, \cite{Tru}.  
Nevertheless, a complete analytic characterization of the 
corresponding classes of admissible potentials has been an open problem until now.

More precisely, we characterize the class of potentials $Q \in \D'(\R^n)$ 
which are $-\Delta$-form bounded with relative bound zero, i.e.,  
for every $\epsilon>0$,
there
exists  $C(\epsilon)>0$ such that 
\begin{equation}\label{E:1.1}
\left \vert \langle Q u, \, u\rangle \right \vert  \leq \epsilon \, 
||\nabla u||^2_{L^2(\R^n)} + C(\epsilon) \, ||u||^2_{L^2(\R^n)}, 
\quad \forall u  \in
C^\infty_0(\R^n). 
\end{equation}

The preceding inequality ensures that, 
 in case $Q$ is real-valued,
 a semi-bounded self-adjoint Schr\"odinger operator $H = H_0 + Q$ 
can be defined on $L^2(\R^n)$ so that 
the quadratic form domain $Q(H)$ coincides 
with $Q(H_0)  = 
W^{1,2} (\R^n)$. For complex-valued $Q$, it follows that 
$H$ is  an $m$-sectorial operator on $L^2(\R^n)$ with 
$\mathcal{D}(H)\subset W^{1, \, 2}(\R^n)$ (\cite{RSi}, Sec. X.2; 
\cite{EE}, Sec. IV.4). 

Our characterization of (\ref{E:1.1}) 
uses only the functions $|\nabla (1-\Delta)^{-1} \, Q|$ and 
$|(1-\Delta)^{-1} \, Q|$ (see Sec.~\ref{Section 2}, Theorem I), and is based on the 
representation: 
\begin{equation}\label{E:1.rep}
Q = {\rm div} \, \vec \Gamma + \gamma, \qquad 
\vec \Gamma (x) = -\nabla (1-\Delta)^{-1} \, Q, \quad \gamma = (1-\Delta)^{-1} \, Q.
\end{equation}
In particular, we will demonstrate that, 
necessarily,  $\vec \Gamma \in  L^2_{\rm loc}(\R^n)^n, 
\quad \gamma \in L^1_{\rm loc}(\R^n),$ 
and, when $n\ge 3$, 
\begin{equation}\label{E:1.1b}
\lim_{\delta\to +0} \, 
\sup_{x_0\in \R^n} \, \delta^{2-n} \int_{B_\delta(x_0)} \left ( 
|\vec \Gamma(x)|^2 + |\gamma(x)| \right) \, dx =0,
\end{equation} 
once (\ref{E:1.1}) holds. Here  $B_\delta (x_0)$ is a Euclidean ball 
of radius $\delta$   centered at $x_0$. 

In the opposite direction, it follows from our results that (\ref{E:1.1}) 
holds whenever
\begin{equation}\label{E:1.1c}
\lim_{\delta\to +0} \, 
\sup_{x_0\in \R^n} \, \delta^{2r-n} \int_{B_\delta(x_0)} \left ( 
|\vec \Gamma(x)|^{2} + |\gamma(x)| \right)^r \, dx =0,
\end{equation}
where $r>1$. Such admissible potentials form  a natural   
analogue of the Fefferman$-$Phong class \cite{Fef}  for the infinitesimal 
form boundedness problem, where cancellations between the positive and negative 
parts of $Q$  come into play. It 
includes functions with 
highly oscillatory behavior as well as singular measures, and  contains 
properly the class of potentials  based on the original Fefferman$-$Phong condition 
where  $|Q|$ is used  in (\ref{E:1.1c}) in place of 
$|\vec \Gamma|^{2} + |\gamma|$. Moreover, one can expand  
this class further using a sharp condition due to Chang, Wilson, and Wolff 
\cite{ChWW}  applied to $|\vec \Gamma|^{2} + |\gamma|$. 

A complete  characterization of (\ref{E:1.1}) given below involves discrete Carleson measure  
sums over dyadic cubes, together with equivalent local energy and pointwise potential conditions 
(Sec.~\ref{Section 2}, Theorem II). In the proofs  
given in Sec.~\ref{Section 5} we overcome considerable 
technical difficulties using a series of new crucial estimates for powers of 
equilibrium potentials, factorization of functions in Sobolev spaces, 
 and theory of $A_p$-weights, along with appropriate localization arguments, and good 
understanding of trace 
inequalities for nonnegative potentials $Q$ studied respectively in 
Sec.~\ref{Section 3} and 
Sec.~\ref{Section 4}. 

Among well-known sufficient conditions for (\ref{E:1.1}), 
which ignore possible cancellations, we  mention: 
$Q\in L^{\frac n 2} (\R^n) + L^\infty(\R^n)$ 
($n\ge 3$) and  $Q\in L^{r} (\R^2) + L^\infty(\R^2)$, $r>1$ ($n=2$) (see \cite{BrK}), 
as well as Kato's condition $K_n$ introduced in \cite{Ka1}:
\begin{align}
\lim_{\delta \to +0} \, \sup_{x_0\in\R^n} \, & \int_{B_\delta(x_0)} 
 \frac {|Q(x)|} {|x-x_0|^{n-2}}
 \, dx
 = 0, \qquad n \ge 3,  \label{E:K_n}\\
\lim_{\delta \to +0} \, \sup_{x_0\in\R^n} \, & \int_{B_\delta(x_0)} 
\log  \frac 1 {|x-x_0|}  \, |Q(x)| 
 \, dx 
 = 0, \qquad n =2. \label{E:K_2}
\end{align}

Kato's class proved to be  especially important in studies of Schr\"odinger semigroups, 
Dirichlet forms,   
and Harnack  inequalities \cite{Agm}, \cite{AiSi}, \cite{Sim}. Our results yield  
that actually (\ref{E:1.1}) holds for a substantially 
broader class of potentials for which $|\vec \Gamma|^2 + |\gamma| 
\in K_n$. 

We emphasize that  
no  a priori assumptions 
were imposed above on $C(\epsilon)$. An observation due to Aizenman and Simon 
 states that, under the hypothesis 
$C(\epsilon) \le a \, e^{b \, \epsilon^{-p}}$ for some $a, b>0$ and $0<p<1$, 
 all   potentials $Q$ 
which obey (\ref{E:1.1}) with $|Q|$ in place of $Q$,  are contained in Kato's class. 
This was first  
proved in \cite{AiSi}  using the Feynman-Kac formalism. In Sec.~\ref{Section 4}, 
we give a sharp 
result of this kind with a simple analytic proof.  
 We show that 
if (\ref{E:1.1}) holds with $|Q|$ in place of $Q \in L^1_{{\rm loc}}(\R^n)$, 
then for any $C(\epsilon)>0$,  
\begin{align}
& \sup_{x_0\in\R^n} \, \int_{B_\delta(x_0)} 
 \frac {|Q(x)|} {|x-x_0|^{n-2}}
 \, dx \le c \, 
\int_{\delta^{-2}}^{+\infty} \frac {\hat C(s)} {s^2} \, ds, 
\qquad n \ge 3, \label{E:k1} \\
& \sup_{x_0\in\R^2} \, \int_{B_\delta(x_0)} \log \frac 1 {|x-x_0|} \, 
|Q(x)|  \, dx \le c \, 
\int_{\delta^{-2}}^{+\infty} \frac {\hat C(s)} {s^2 \log s} \, ds, 
\qquad n =2,
\label{E:k2} 
\end{align}
where $c$ is a constant which depends only on $n$, and $\delta$ is sufficiently small. Here 
$\hat C(s) = \inf_{\epsilon>0} \, \{C(\epsilon) + s \, \epsilon\}$ is the Legendre transform of $-C(\epsilon)$.  
In particular, it follows that the condition $C(\epsilon) \le a \, e^{b \, \epsilon^{-p}}$ for  {\it any} 
 $p>0$ is enough  to ensure 
that $Q \in K_2$ in the more subtle two-dimensional case.

In the second part of the paper, we study quadratic form inequalities of Trudinger type 
  where $C(\epsilon)$ in (\ref{E:1.1}) has   power growth, 
i.e., there exists $\epsilon_0>0$ such that 
  \begin{equation}\label{E:1.2}
\left \vert \langle Q u, \, u\rangle \right \vert  \leq \epsilon \, 
||\nabla u||^2_{L^2(\R^n)} + c \,   \epsilon^{-\beta}  \,
||u||^2_{L^2(\R^n)},  \quad \forall u 
\in C^\infty_0(\R^n), 
\end{equation}
 for every
$\epsilon\in(0, \epsilon_0)$, where $\beta>0$. Such inequalities 
appear in  studies of elliptic PDE with measurable 
coefficients \cite{Tru}, and have been  used extensively in spectral theory 
 of the  Schr\"odinger operator (see, e.g., \cite{RSS}, Sec. 20). 

As it turns out, it is still possible to 
characterize (\ref{E:1.2}) using only $|\vec \Gamma|$ and 
$|\gamma|$ defined by (\ref{E:1.rep}), provided $\beta>1$. We will show below (Theorem IV) that 
in this case (\ref{E:1.2})
 holds if and only if 
both of the following conditions hold:
\begin{align}
 \sup_{\substack{x_0\in \R^n \\ 0<\delta< \delta_0}} \, 
\delta^{2 \frac {\beta-1}{\beta+1}-n} 
& \int_{B_\delta(x_0)} 
|\vec \Gamma(x)|^{2} \, dx < +\infty, \label{E:1.2a}\\
\sup_{\substack{x_0\in \R^n \\ 0<\delta< \delta_0}} \, 
\delta^{\frac {2\beta}{\beta+1}-n} 
& \int_{B_\delta(x_0)}  |\gamma(x)| \, dx < +\infty, 
\label{E:1.2b}
\end{align}
for some $\delta_0>0$. 
However, in the case $\beta\le 1$  
this is no longer true. For $\beta=1$, (\ref{E:1.2a}) has to be replaced with  
the condition that $\vec \Gamma$ is in the local ${\rm BMO}$ space, 
or respectively is H\"older-continuous of order $\frac {1-\beta}{1+\beta}$ if 
$0<\beta<1$. 

In the homogeneous case $\epsilon_0=+\infty$, (\ref{E:1.2})   
is equivalent to the  {\em multiplicative inequality\/}:
\begin{equation}\label{E:1.3}
\left \vert \langle Q u, \, u\rangle \right \vert  \leq C \, 
||\nabla u||^{2p}_{L^2(\R^n)}  \,||u||^{2(1-p)}_{L^2(\R^n)}, 
\quad \forall u 
\in C^\infty_0(\R^n), 
\end{equation}
where $p = \frac {\beta}{1 + \beta} \in (0, \, 1)$. In 
spectral theory,  (\ref{E:1.3}) 
 is  referred to as the form $p$-subordination property 
\cite{Agr}, \cite{Gr1}, \cite{Gr2}, \cite{MM1}, \cite{RSS}.

For {\em nonnegative\/}  potentials $Q$,  where $Q$ coincides with a 
locally finite measure $\mu$  on $\R^n$,  inequality
(\ref{E:1.3}) is known to hold if and only 
if 
\begin{equation}\label{E:1.4}
\mu \left ( B_\delta (x_0)\right)  \le c \, \delta^{n - 2p}, 
\end{equation}
where the constant $c$ does not depend on $\delta>0$ and
$x_0 \in \R^n$ (\cite{M2}, Sec. 1.4.7).

 For general $Q$, we obtain the following result (Theorem V): If $p> \frac 1 2$,
then  {\rm (\ref{E:1.3})} 
holds if and only if $\nabla \Delta^{-1} Q$ lies in  the 
 Morrey  space 
$\mathcal{L}^{2, \, \lambda}(\R^n)$, where $\lambda= n + 2 - 4p$. 
For $p=\frac 1 2$,  it holds if and only if   
$\nabla \Delta^{-1} Q \in {\rm BMO}(\R^n)$, and for
$0<p<\frac 1 2$, whenever  $\nabla \Delta^{-1} Q \in {\rm Lip}_{1-2p}(\R^n)$.   
These different 
characterizations  are equivalent 
to  {\rm (\ref{E:1.4})} 
 if $Q$ is a nonnegative measure.

As a consequence, we are able to characterize 
those  $Q$ which obey an analogous inequality of Nash's type: 
\begin{equation}\label{E:1.3a}
\left \vert \langle Q u, \, u\rangle \right \vert  \leq C \, 
||\nabla u||^{2p}_{L^2(\R^n)}  \,||u||^{2(1-p)}_{L^1(\R^n)}, 
\quad \forall u 
\in C^\infty_0(\R^n). 
\end{equation}
where $p  \in (0, \, 1)$. In fact,  the preceding inequality has two critical 
exponents,   $p_* = \frac n {n+2}$ and $p^* = \frac {n+1} {n+2}$. 
We will show in Sec.~6 (Corollary~\ref{Corollary 3.9}) that, for 
$0<p<p_*$,  (\ref{E:1.3a}) holds only if $Q=0$, whereas for 
$p=p_*$ it follows that $Q \in L^\infty(\R^n)$, i.e., it is 
equivalent to Nash's inequality (\cite{LL}, Theorem 8.13). For $p> p_*$, 
the validity of (\ref{E:1.3a}) 
is equivalent respectively to: $\, \nabla \Delta^{-1} Q \in {\rm Lip}_{n+1-p(n+2)}(\R^n)$ 
if $p_*<p<p^*$;  $\nabla \Delta^{-1} Q \in {\rm BMO}(\R^n)$ 
if $p=p^*$;  
and   $\nabla \Delta^{-1} Q \in \mathcal{L}^{2, \, \lambda}(\R^n)$, where 
$\lambda= 3n + 2 - 2p(n+2)$, if $p^*<p<1$. 

There is an interesting 
connection of our form subordination theorems, 
in the sufficiency part, with  estimates of the type:
$$
||(\vec u  \cdot \nabla) \, \vec  u||_{\mathcal{H}^1(\R^n)} \le c \, 
||\vec u ||_{L^2(\R^n)} \, ||\nabla \vec u||_{L^2(\R^n)},  
 \quad {\rm div} \, 
\vec u = \vec 0, \quad \forall \vec u \in C^\infty_0(\R^n)^n,
$$ 
where ${\mathcal H}^1(\R^n)$ is a real Hardy space (\cite{St2}). 
Such inequalities and their generalizations 
considered in  Sec.~\ref{Section 6} (see Lemma~\ref{Lemma 3.7}) are 
useful in hydrodynamics 
 (\cite{CLMS}, \cite{Co}).

Our proofs of 
the necessity statements 
are generally more delicate than those of the sufficiency. 
They make use of factorization of functions in Sobolev spaces, 
 and sharp mean oscillation estimates for 
``anti-derivatives'' of $Q$. Our methods based on 
a combination of ideas from  potential theory,  Sobolev spaces,  nonlinear 
and harmonic analyses  have been developed over the past ten years  
 \cite{HMV}, \cite{MV1}-\cite{MV4}.  We believe that this approach has a broader scope, 
and might be useful 
with regards to  higher 
order elliptic operators, pseudodifferential operators, harmonic analysis and 
nonlinear PDE problems (see, e.g.,  \cite{BoBr}, \cite{Dav1}, \cite{Dav2}, \cite{Fef}).

\medskip

We conclude the introduction with a brief outline of the contents of the 
paper. In Sec.~\ref{Section 2}, we introduce some notation and state 
our main results, Theorems I$-$V.  In Sec.~\ref{Section 3}, we present a  useful 
localization principle for quadratic form 
inequalities using a Legendre transform associated with the function $C(\epsilon)$. 
This localization will be employed 
throughout the paper. Sec.~\ref{Section 4} is devoted to a  study  of 
the corresponding integral inequalities for 
nonnegative potentials $Q$ (locally finite measures). In particular, 
we discuss connections with Kato's  
condition which is important to  
Schr\"odinger operators and elliptic PDE.

The proofs of the main results are given in 
Sections~\ref{Section 5} 
and \ref{Section 6}   
where we 
treat general distributional potentials $Q$. In Sec.~\ref{Section 5}, we  
establish necessary and sufficient conditions for  
$-\Delta$-form boundedness  of the potential energy operator $Q$ with 
relative bound zero. In Sec.~\ref{Section 6}, we obtain 
subordination criteria for quadratic forms associated with 
 the Schr\"odinger operator, and discuss  connections with the compensated 
compactness phenomenon.

\vfill
\eject
 
\section{Main results}\label{Section 2}

Throughout the paper we will be using the following notation and conventions. We denote by $W^{1,2}(\R^n)$
 the 
 Sobolev space of 
weakly differentiable functions on $\R^n$ $(n\ge 1)$ such that 
$$
||u||_{W^{1,2}(\R^n)} = ||u||_{L^2(\R^n)} + ||\nabla u||_{L^2(\R^n)}< +\infty, 
$$
and  by $W^{-1,2} (\R^n) = W^{1,2} (\R^n)^*$ the dual Sobolev space. 
For a compact set $e\subset \R^n$, the capacity associated with $W^{1,2}(\R^n)$
 is defined by 
$$
{\rm cap} \, (e) = \inf \, \left\{ \, ||u||^2_{W^{1,2} (\R^n)} \, : \quad 
u \in C^\infty_0(\R^n), \quad u(x)>1 \, \,  {\rm on} \, \,  e \right\}.
$$

For $0< r < \infty$, we denote  by $L^r_{\rm unif}(\R^n)$ 
all $f\in L^r_{\rm loc}(\R^n)$ such that 
$$
||f||_{L^r_{\rm unif}} = \sup_{x_0 \in \R^n} \, 
||\chi_{B_1 (x_0)} \,  f ||_{L^r(\R^n)} < \infty.
$$
By  $L^r (\R^n)^n$
we denote  the class of vector fields 
$\vec \Gamma = \{\Gamma_j\}_{j=1}^n : \, \R^n \to \C^n$, such that 
$\Gamma_j \in L^r (\R^n)$, $j =1, 2, \dots, n$, and use  similar notation for other vector-valued 
function spaces. 

By $M^+(\R^n)$ we denote the class of nonnegative locally finite Borel measures on $\R^n$. 
If $Q\in \D'(\R^n)$ is nonnegative, i.e.,  coincides with  $\mu\in M^+(\R^n)$, we write 
$\int_{\R^n} |u(x)|^2 \, d \mu$ in place of $\langle Q, \, |u|^2\rangle = 
\langle Q u, \, u\rangle$ for the quadratic
form associated with the distribution $Q$, if $u \in C^\infty_0(\R^n)$. 
Sometimes we will use  $\int_{\R^n}
|u(x)|^2 \, Q(x) \, dx$ in place of
$\langle Q u, \, u\rangle$  even if $Q$ is not in
$L^1_{\rm loc} (\R^n)$.

 We set 
$$
m_B (f) = \frac 1 {|B|} \int_B f(x) \, dx
$$
for a ball $B\subset \R^n$, and denote by  
 $\text{BMO}(\R^n)$ the class of 
$f \in L^r_{\rm loc} (\R^n)$ for which 
$$
\sup_{x_0\in \R^n, \, \delta>0} \, \,  \frac {1}{|B_\delta(x_0)|} 
\int_{B_\delta (x_0)} |f(x)-m_{B_\delta(x_0)}(f)|^r \, dx < + \infty,
$$
for any (or, equivalently, all) $1\le r < +\infty$. 

We will also need  
an inhomogeneous version of $\text{BMO}(\R^n)$ (the so-called 
local $\text{BMO}$) which we denote by 
$\text{bmo}  (\R^n)$. It can be  defined in a similar way as 
the set of  $f \in L^r_{\rm unif}(\R^n)$ such that the preceding condition holds 
for $0<\delta\le 1$ (see  \cite{St2}, p. 264). 

The  Morrey space $\mathcal{L}^{r, \, \lambda}$ ($r>0, \, \lambda>0$) 
consists of $f \in  L^r_{\rm loc} (\R^n)$ such that 
$$
\sup_{x_0\in \R^n, \, \delta>0} \,\,  \frac {1}{|B_\delta(x_0)|^{1-\frac \lambda n}} 
\int_{B_\delta(x_0)} |f|^r \, dx < + \infty.
$$
In the corresponding inhomogeneous analogue, we set $0<\delta\le 1$ in the preceding 
inequality. It will be clear from the context which version of the Morrey space is used.

We now state our main results. 
\vskip12pt

\noindent {\bf Theorem I.} 
{\it 
Let $Q \in \mathcal \D'(\R^n)$, $n \ge 2$.
Then the  following statements hold. 

{\rm (i)} Suppose that $Q$ is represented in the form:
\begin{equation}\label{E:1.5}
Q = {\rm div} \,\,  \vec \Gamma + \gamma,
\end{equation} 
 where
 $\vec \Gamma  \in L^2_{\rm loc} (\R^n)^n$ and 
$\gamma \in  L^1_{\rm loc} (\R^n)$
 satisfy respectively the conditions: 
\begin{equation}\label{E:1.6}
  \lim_{\delta \to +0}  \sup_{x_0 \in \R^n} \, \sup_{u}
 \, \frac{
\int_{B_\delta(x_0)} |\vec \Gamma(x)|^2  \, |u(x)|^2 \, dx}{||\nabla u||^2_{L^2(B_\delta(x_0))}}  =0,
\end{equation}
\begin{equation}\label{E:1.7}
  \lim_{\delta \to +0} \sup_{x_0 \in \R^n}  \,  \sup_{u} \,   
\frac { \int_{B_\delta(x_0)} | \gamma(x)| \, |u(x)|^2 \, dx} { ||\nabla u||^2_{L^2(B_\delta(x_0))} }   =0,
\end{equation} 
where $u \in C^\infty_0(B_\delta(x_0))$, $u \not\equiv 0$. 
Then $Q$ is infinitesimally form bounded with respect to
$-\Delta$, i.e., for every $\epsilon >0$ there exists 
$C(\epsilon)>0$ such that  {\rm (\ref{E:1.1})} holds.  

{\rm (ii)}  Conversely, suppose $Q$ is infinitesimally form bounded 
with respect to
$-\Delta$.  Then $Q$ can be represented in the
form {\rm (\ref{E:1.5})} so that both {\rm (\ref{E:1.6})} and {\rm (\ref{E:1.7})}
hold. Moreover, one can set $\vec \Gamma = - \nabla (1-\Delta)^{-1} Q$ and 
$\gamma = (1-\Delta)^{-1} Q$ in the representation {\rm (\ref{E:1.5})}.
}
\vskip12pt

\noindent{\bf Remark 2.1.} In the statement of Theorem I one can 
replace conditions {\rm (\ref{E:1.6})} and {\rm (\ref{E:1.7})} with  
the equivalent condition where $|(1-\Delta)^{-\frac 1 2} Q|^2$ 
is used in place of $|\vec \Gamma|^2$ in {\rm (\ref{E:1.6})}. 
\vskip12pt

The importance of Theorem I is in the means it provides for deducing 
explicit criteria of  
 the infinitesimal form 
boundedness  in terms of the 
{\it nonnegative} locally integrable functions $|\vec \Gamma|^2$ and 
$|\gamma|$. \vskip12pt

\noindent {\bf Theorem II.} 
{\it Let $Q \in \mathcal \D'(\R^n)$, $n \ge 2$. The 
following statements are equivalent:  

{\rm (i)} $Q$ is infinitesimally form bounded with respect to $-\Delta$. 

{\rm (ii)} $Q$ has the form {\rm (\ref{E:1.5})} where 
$\vec \Gamma = - \nabla (1-\Delta)^{-1} \, Q$,
$\gamma = (1-\Delta)^{-1} \, Q$, and the measure $\mu\in M^+(\R^n)$ defined by 
\begin{equation}\label{E:1.8}
d \mu = \left ( |\vec \Gamma (x)|^2 + |\gamma(x)| \right) \, dx
\end{equation} 
has the property that, for every $\epsilon>0$, there exists $C(\epsilon)>0$ 
such that  
\begin{equation}\label{E:1.9}
 \int_{\R^n} |u(x)|^2 \, d \mu \le \epsilon \, ||\nabla u||^2_{L^2(\R^n)} 
+ C(\epsilon) \, ||\nabla u||^2_{L^2(\R^n)}, \quad \forall u \in C^\infty_0(\R^n).
\end{equation} 

{\rm (iii)} For $\mu$ defined by {\rm (\ref{E:1.8})}, 
\begin{equation}\label{E:1.11}
 \lim_{\delta \to +0} \, \sup_{P_0: \, {\rm diam} \,  P_0 \le \delta} \, 
 \frac {1} {\mu(P_0)} \, \sum_{P \subseteq P_0}  \,  \frac 
{\mu(P)^2}  {|P|^{1- \frac 2 n}} 
 = 0, 
\end{equation}
where $P$, $P_0$ are dyadic cubes in $\R^n$.

{\rm (iv)} For $\mu$ defined by {\rm (\ref{E:1.8})},
\begin{equation}\label{E:1.10} 
\lim_{\delta \to +0} \sup_{e: \, \text{{\rm diam}} \,  e \le \delta}
\frac {\mu(e)}
{ \text{\rm{cap}} \, (e)} = 0,
\end{equation}
where $e$ denotes a  compact set of positive capacity in $\R^n$.

 {\rm (v)} For $\mu$ defined by {\rm (\ref{E:1.8})}, 
\begin{equation}\label{E:1.12a}
 \lim_{\delta \to +0} \,  \sup_{x_0\in \R^n} \,  \frac {\left \Vert 
 \mu_{B_\delta(x_0)} \, \right \Vert^2_{W^{-1,2} (\R^n)}}
{ \mu (B_\delta(x_0))}  = 0,
\end{equation}
where  $\mu_{B_\delta(x_0)}$ is the restriction of $\mu$ to the ball 
$B_\delta(x_0)$.

{\rm (vi)} For $\mu$ defined by {\rm (\ref{E:1.8})},
\begin{equation}\label{E:1.13a}
 \lim_{\delta \to +0} \, \sup_{x, \, x_0\in \R^n} \, 
 \frac { G_1 \star \left ( G_1 \star  
\mu_{B_\delta(x_0)} \right )^2 (x)}{G_1 \star \mu_{B_\delta(x_0)}(x) } = 0,
\end{equation}
where
 $G_1 \star \mu = (1-\Delta)^{-\frac 1 2} \mu$ is the Bessel potential of order $1$.
}
\vskip12pt

It is worth noting that although  Theorem II holds in the two-dimensional 
case, its  proof requires certain  
modifications in comparison to $n\ge 3$. In the one-dimensional case,   
 the infinitesimal form boundedness 
of the Sturm-Liouville operator $H  = -\frac{d^2 \,}{d \, x^2} +Q$ on 
$L^2(\R^1)$ 
is actually a  consequence of the form boundedness. \vskip12pt

\noindent {\bf Theorem III.} 
{\it Let $Q \in \mathcal \D'(\R^1)$. Then the 
following statements 
are equivalent.

{\rm (i)} $Q$ is infinitesimally form bounded with respect to 
$-\frac{d^2}{dx^2}$. 

{\rm (ii)} $Q$ is  form bounded with respect to 
$-\frac{d^2}{dx^2}$, i.e.,
$$
|\langle Q \, u, \, u\rangle| \le {\rm const} \, ||u||^2_{W^{1,2}(\R^1)}, \quad 
\forall u \in C^\infty_0(\R^1).
$$

{\rm (iii)} $Q$ can be represented in the form $Q =  \frac {d \Gamma} {dx} + \gamma$, where 
\begin{equation}\label{E:1.12}
\sup_{x \in \R^1} \int_x^{x+1} \left ( |\Gamma(x)|^2 + |\gamma(x)| \right) \, 
dx < +\infty.
\end{equation}

{\rm (iv)} Condition {\rm (\ref{E:1.12})}
  holds where  
$$
\Gamma(x) = \int_{\R^1} {\rm sign} \, (x-t) \, e^{-|x-t|} \, Q(t) \, dt, \qquad  
\gamma (x) =  \int_{\R^1} e^{-|x-t|} \, Q(t) \, dt
$$
are understood in the distributional sense. 
}
\vskip12pt

The statement   (iii)$\Rightarrow$(i) in Theorem III is known 
(\cite{Sch}, Theorem 11.2.1), 
whereas (ii)$\Rightarrow$(iv) 
follows from our earlier results \cite{MV2}, \cite{MV3}.

We now state a characterization of the  form subordination property {\rm (\ref{E:1.2})}. 
It was formulated originally in \cite{Tru},  in the  form of the inequality: 
\begin{equation}\label{E:1.2'}
\left \vert \langle Q u, \, u\rangle \right \vert  \leq \epsilon \, 
||\nabla u||^2_{L^2(\R^n)} + c \, \epsilon^{-\nu} \,
||u||^2_{L^1(\R^n)},  \quad \forall u 
\in C^\infty_0(\R^n), 
\end{equation}
for  $Q \ge 0$. Such $Q$ are called $\epsilon^\nu$-{\it compactly bounded\/} in \cite{Tru}.   
It follows from Nash's inequality 
that {\rm (\ref{E:1.2})} yields {\rm (\ref{E:1.2'})} with $\nu = 
\frac {n+2} 2 \beta + \frac n 2$; the converse is also true, provided $\nu>\frac n 2$, 
and is deduced below using  a localization argument (Sec.~\ref{Section 3}, 
Corollary~\ref{Corollary 2.3}). In the critical case $\nu = \frac n 2$, 
{\rm (\ref{E:1.2'})} holds if and only if $Q \in L^\infty(\R^n)$, while 
for $0<\nu < \frac n 2$, it holds only if $Q=0$.

Necessary and sufficient conditions for {\rm (\ref{E:1.2})}, 
or equivalently  {\rm (\ref{E:1.2'})} with $\nu = 
\frac {n+2} 2 \beta + \frac n 2$, 
can be formulated in terms of  
Morrey-Campanato spaces using mean oscillations of the functions 
 $\vec \Gamma$ and $\gamma$ which have appeared in  Theorems I$-$III.  \vskip12pt

\noindent {\bf Theorem IV.} 
{\it 
Let $Q \in \mathcal \D'(\R^n)$, $n \ge 2$, and let $0<\beta<+\infty$.

{\rm (i)} Suppose there exists $\epsilon_0>0$ such that 
 {\rm (\ref{E:1.2})} holds 
for every $\epsilon\in (0, \epsilon_0)$. Then $Q$ can be represented in the
form  
\begin{equation}\label{E:1.13}
Q = \text{{\rm div}} \,\,  \vec \Gamma + \gamma,
\end{equation} 
 where $\vec \Gamma =  - \nabla (1-\Delta)^{-1} Q\in L^2_{\rm loc} 
(\R^n)^n$ 
and $\gamma = 
(1-\Delta)^{-1} Q\in L^1_{\rm loc} (\R^n)$.
  
Moreover,  there exists $\delta_0>0$ such that 
\begin{equation}\label{E:1.14}
\int_{B_\delta(x_0)} |\vec \Gamma (x)- m_{B_\delta(x_0)} (\vec \Gamma)|^2 
\, dx \le  c \, \delta^{n-2 \frac{\beta-1} {\beta+1}}, \quad 
0< \delta < \delta_0, 
\end{equation}
 \begin{equation}\label{E:1.15}
 \int_{B_\delta (x_0)} |\gamma (x)| \,  dx \le c \, 
\delta^{n - \frac{2 \beta} {\beta+1}},
\quad   0< \delta < \delta_0,
\end{equation} 
where  $c$ does not depend on $x_0\in R^n$ and $\delta_0$. Furthermore, 
$\vec \Gamma \in L^2_{\rm unif}  (\R^n)^n$ if $\beta \ge 1$ and 
$\vec \Gamma \in L^{\infty}  (\R^n)^n$ if $0<\beta<1$.

{\rm (ii)} Conversely, if $Q$ is given by {\rm (\ref{E:1.13})} where
 $\vec \Gamma  \in \mathbf{L}^2_{\rm loc} (\R^n)$, 
$\gamma \in L^1_{\rm loc} (\R^n)$ 
 satisfy   
{\rm (\ref{E:1.14})},  {\rm (\ref{E:1.15})} 
for all $0<\delta<\delta_0$ then  there exists $\epsilon_0>0$ 
such that {\rm (\ref{E:1.2})} holds for all $0<\epsilon<\epsilon_0$.
}
\vskip12pt

\noindent{\bf Remark 2.2.} (a) In the case $\beta=1$, it follows 
 that  {\rm (\ref{E:1.14})} holds if and only if
$\vec \Gamma  \in {\text{bmo}} 
(\R^n)$. In other words, $ Q \in  {\text{bmo}}_{-1}  
(\R^n)= F^{2, \infty}_{-1} (\R^n)$, where 
$F^{p, q}_{\alpha}$ stands for the scale of Triebel-Lizorkin  
spaces (see, e.g.,  \cite{AH}, \cite{Tri}).

(b) Similarly, in the case $0<\beta<1$, 
 {\rm (\ref{E:1.14})} holds if and only if $\vec \Gamma$ is H\"older-continuous:
$$
|\vec \Gamma(x) - \vec \Gamma(x')|\le c \, 
|x-x'|^{\frac{1- \beta} {\beta+1}}, \quad |x-x'|<\delta_0.
$$

(c) For $\beta>1$,  {\rm (\ref{E:1.14})} holds if and only if 
$\vec \Gamma$ is in the 
corresponding Morrey space $\mathcal{L}^{2, \, n-2\frac{\beta-1} {\beta+1}}$:
$$
\int_{B_\delta(x_0)} |\vec \Gamma (x)|^2 
\, dx \le  c \, \delta^{n-2 \frac{\beta-1} {\beta+1}}, \quad 
0< \delta < \delta_0.
$$
Note that, according to {\rm (\ref{E:1.15}),
  $\gamma \in \mathcal{L}^{1, \, n-\frac{2\beta} {\beta+1}}$. 
Equivalently, $Q$ is in the corresponding Campanato-Morrey space  
$\mathcal{L}^{2, \, n-2\frac{\beta-1} {\beta+1}}_{-1}$ 
of negative order.

\noindent{\bf Remark 2.3.} (a) An immediate consequence of Theorem IV is 
that, for all $\beta>0$, {\rm (\ref{E:1.2})} is equivalent to
the following localized condition: 
$$
||(1-\Delta)^{-\frac 1 2} (\eta_{\delta, \, x_0} \, Q)
||^2_{L^2(B_\delta(x_0))} \le c 
\, \delta^{n-2 \frac{\beta-1} {\beta+1}}, \quad 
0< \delta < \delta_0, \, \, x_0 \in \R^n,
$$
where $\eta_{\delta, \, x_0}(x) = \eta \left (\frac {|x-x_0|}{\delta}\right)$, 
and $\eta$ is a smooth cut-off function with compact support.

(b) A similar energy condition, 
$$ 
||(1-\Delta)^{-\frac 1 2} (\eta_{\delta, \, x_0} \, Q)
||^2_{L^2(\R^n)} \le c 
\, \delta^{n-2 \frac{\beta-1} {\beta+1}}, \quad 
0< \delta < \delta_0, \, \, x_0 \in \R^n,
$$
is sufficient but generally  {\it not necessary\/}
in the case $n=2$.

We next state a criterion for the multiplicative inequality 
{\rm (\ref{E:1.3})} to hold, which is equivalent to a 
homogeneous version of {\rm (\ref{E:1.2})} with $\epsilon_0=+\infty$ 
and $p = \frac \beta {\beta+1}$. \vskip12pt

\noindent {\bf Theorem V.} 
{\it 
Let $Q \in \mathcal \D'(\R^n)$, $n \ge 2$, and let $0<p<1$.

{\rm (i)} Suppose that {\rm (\ref{E:1.3})} holds.  
Then $Q$ can be represented in the
form  
\begin{equation}\label{E:1.16}
Q = \text{{\rm div}} \,\,  \vec \Gamma,
\end{equation} 
 where $\vec \Gamma =  \nabla \Delta^{-1} Q$, and one of the 
 following conditions hold:
\begin{equation}\label{E:1.17}
\vec \Gamma \in {\rm BMO}(\R^n) \quad \text{if} \quad  p= \frac 1 2; 
\quad \vec \Gamma \in
{\rm Lip}_{1-2p} (\R^n)\quad \text{if} \quad  0<p< \frac 1 2;  
\end{equation}
 \begin{equation}\label{E:1.18}
 \int_{B_\delta (x_0)} |\vec \Gamma (x)|^2 \,  dx \le c \, \delta^{n+2-4p},
\qquad   \text{if} \quad  \frac 1 2<p< 1; 
\end{equation} 
where  $c$ does not depend on $x_0\in \R^n$, $\delta>0$.

{\rm (ii)} Conversely, if $Q$ is given by {\rm (\ref{E:1.16})} where
 $\vec \Gamma  \in L^2_{\rm loc} (\R^n)^n$ and 
 satisfies  
{\rm (\ref{E:1.17})},  {\rm (\ref{E:1.18})} then  {\rm (\ref{E:1.3})} holds.
}
\vskip12pt

\noindent{\bf Remark 2.4.} In  Theorem V, the 
``antiderivative''  
 $\vec \Gamma =  \nabla \Delta^{-1} Q$ can be replaced 
by $(-\Delta)^{-\frac 1 2} Q$. Furthermore, as a corollary we deduce that 
  (\ref{E:1.3})  holds if and only if  
$Q \in \text{BMO}_{-1} (\R^n)= \dot F^{2, \infty}_{-1}(\R^n)$ 
for $p= \frac 1 2$, and    
$Q \in \dot B^{\infty,  \infty}_{-2p}(\R^n)$ for $0<p<\frac 1 2$. 
Here  $\dot F^{r, q}_{\alpha}$ 
and $\dot B^{r, q}_{\alpha}$ are homogeneous Triebel-Lizorkin and 
Besov spaces respectively (see \cite{Tri}).

In the case  $p=\frac 1 2$, statement (ii) 
of Theorem V (sufficiency of the condition 
$\vec \Gamma \in \text{BMO}$) 
is equivalent via the $\mathcal{H}^1-\text{BMO}$ duality  to the inequality
 \begin{equation}\label{E:1.19}
||u \, \nabla u||_{\mathcal{H}^1(\R^n)}  \le c \, ||u||_{L^2(\R^n)} \, 
||\nabla u||_{L^2(\R^n)}, \qquad 
 \forall u \in C^\infty_0(\R^n).
\end{equation}
Here $\mathcal{H}^1(\R^n)$ is the real Hardy space on $\R^n$ (\cite{St2}).

The preceding estimate yields the 
following vector-valued  
inequality   
which is used in studies of the Navier-Stokes equation, and   
is  related to the  {\it compensated compactness}  
phenomenon   (see \cite{Co}, \cite{CLMS}): 
\begin{equation}\label{E:1.20}
||(\vec u  \cdot \nabla) \, \vec  u||_{\mathcal{H}^1(\R^n)} \le c \, 
||\vec u ||_{L^2(\R^n)} \, ||\nabla \vec u||_{L^2(\R^n)}, \quad 
 \quad {\rm div} \, 
\vec u = \vec 0,
\end{equation} 
for all $\vec u \in C^\infty_0(\R^n)^n$. 

In Sec.~\ref{Section 6} we deduce a nonhomogeneous 
version of the well-known {\it div-curl} lemma  \cite{CLMS} which 
yields both (\ref{E:1.19}) and (\ref{E:1.20}).

\section{Localization of quadratic form estimates}\label{Section 3}

In this section we restate the infinitesimal form boundedness property 
(\ref{E:1.1}) in an 
equivalent localized form.

\begin{lemma}\label{Lemma 2.1} Suppose $Q \in \mathcal \D'(\R^n)$, $n \ge1$. 
Then the following statements are equivalent.

{\rm (i)} There exists a  positive constant $\epsilon_0$ and a function 
$C(\epsilon) : \, (0, \epsilon_0) \to \R_+$ such that 
the inequality 
\begin{equation}\label{E:2.1}
\left \vert \langle Q u, \, u\rangle \right \vert  \leq \epsilon \, ||\nabla
u||^2_{L^2} + C(\epsilon) \, ||u||^2_{L^2}, \quad \forall u  \in
C^\infty_0(\R^n), 
\end{equation}
holds for every $\epsilon\in (0, \, \epsilon_0)$.

{\rm (ii)} 
\begin{equation}\label{E:2.2}
\lim_{\delta \to +0} \, \sup_{x_0 \in \R^n} \, \sup \, \{ \, 
\left | \langle Q u, \, u \rangle \right | \, : \,  u
\in C^\infty_0(B_\delta(x_0)), \, \,   ||\nabla u||_{L^2} \le 1 
\} = 0. 
\end{equation}
Moreover, {\rm (\ref{E:2.1})} implies that 
$$
 \sup_{x_0 \in \R^n} \, \sup \, \{ \, 
\left | \langle Q u, \, u \rangle \right | \, : \,  u
\in C^\infty_0(B_\delta(x_0)), \, \,   ||\nabla u||_{L^2} \le 1 
\} \le c(n) \, \epsilon,
$$
for $\delta = \sqrt{\frac {\epsilon} {C(\epsilon)}}$, where 
the constant $c(n)$ depends only on $n$.

Conversely, if the preceding inequality holds with $\delta = 
\sqrt{\frac {\epsilon} {C(\epsilon)}}$ and (sufficiently small) constant 
$c(n)$ then ${\rm (\ref{E:2.1})}$ holds.

\end{lemma} 

\begin{proof}  Suppose (\ref{E:2.1}) holds. As was mentioned above, 
we can assume without loss of generality that $\epsilon_0=+\infty$. 
Let $\delta >0$, $x_0 \in \R^n$,
and let $u \in C^\infty_0(B_\delta(x_0))$.  
Then clearly 
\begin{equation}\label{E:2.3}
||u||_{L^2} \le k(n) \, \delta \, ||\nabla u||_{L^2},
\end{equation}
where the constant $k(n)$ depends only on $n$. Hence 
\begin{equation}\label{E:2.4}
\left \vert \langle Q u, \, u\rangle \right \vert  \leq \, 
(\epsilon + k(n)^2 \,\delta^2 \, C(\epsilon) )  \, ||\nabla u||^2_{L^2}, 
\quad \forall u  \in C^\infty_0( B_\delta(x_0)),
\end{equation}
which implies   
\begin{equation}\label{E:2.5}
\begin{split}
\sup_{x_0 \in \R^n} \, & \sup \, \{ \, 
\left | \langle Q u, \, u \rangle \right | \, : \,  u
\in C^\infty_0(B_\delta(x_0)), \, \,   ||\nabla u||_{L^2} \le 1 \} \\ 
& \le \epsilon + k(n)^2 \,\delta^2 \, C(\epsilon). 
\end{split}
\end{equation}
Letting $\delta \to +0$, and then  $\epsilon \to +0$,  in the preceding
inequality, we obtain (\ref{E:2.2}).

Conversely, if (\ref{E:2.2}) holds then, for every $\epsilon>0$ there 
exists $\delta = \delta(\epsilon)$ so that, for every $x_0\in \R^n$,   
\begin{equation}\label{E:2.6}
\left \vert \langle Q v, \, v\rangle \right \vert  \leq \, 
\epsilon  \, ||\nabla v||^2_{L^2(\R^n)}, \quad \forall v \in 
C^\infty_0( B_\delta(x_0)). 
\end{equation}
Now fix $u \in C^\infty_0(\R^n)$ which is supported in $B_R(0)$. 
Let $\eta \in C^\infty_0(B_1(0))$ be a cut-off function 
such  that  $0\le \eta (x)
\le 1$, $\eta(x)=1$ for $|x|\le \frac 1 2$,  $|\nabla \eta(x)| \le c(n)$, 
and let 
$\eta_{\delta, x_0} (x) = \eta \left(\frac {x-x_0}{\delta}\right)$. 

Next, let us pick $x_i \in \R^n$  $(i=1, 2, \dots)$ so that $\{x_i\}$ form a
cubic lattice with grid distance $\frac \delta {2 \, \sqrt{n}}$, and let 
 $\eta_i (x) = \eta \left(\frac
{x-x_i}{\delta}\right)$.  Let  
$$
\phi (x) = \sum_i \, \eta_i^2 (x),
$$
where the sum is taken over  a finite 
number of  indices $i$ such that $B_{2R} (0)\subset \, \cup \, 
B_\delta(x_i)$. Notice that 
$$
1 \le \phi(x) \le \kappa_1(n) \quad \text{on} \quad B_{2R}(0),
$$
and
$$
|\nabla \phi(x)| \le \frac {\kappa_2(n)}{\delta}  \quad \text{on} \quad
B_{2R}(0).
$$

Define 
$$
\zeta_i(x) = \frac {\eta_i (x)} {\sqrt{\phi(x)}},
$$
so that 
$$
\sum_i \zeta_i^2(x) \equiv 1 \quad \text{on} \quad B_R(0).
$$
Then by (\ref{E:2.6}), 
\begin{align}
\left \vert \langle \, Q u, \, u\rangle \right \vert 
& = \left \vert \sum_i \,  \langle \, Q, \, |\zeta_i \,  u|^2 \rangle \right
\vert   \notag \\ &\leq  \sum_i \,  \left \vert \langle \, Q, \, |\zeta_i \, 
u|^2 \rangle \right \vert \notag 
\\ & \leq \, \epsilon \, \sum_i \, \int |\nabla u|^2 \zeta_i^2
\, dx  + \epsilon \sum_i \, \int |u|^2 \, |\nabla \zeta_i|^2 \, dx 
\notag \\ & \leq \, \epsilon \, ||\nabla u||_{L^2}^2 + c(n) \, \epsilon \,
\delta^{-2} \, ||u||_{L^2}^2.\notag
\end{align}
In the last line  we have used the estimate 
$$
\sum_i \, |\nabla \zeta_i (x)|^2 \le \frac {c(n)} {\delta^2} \,  \frac 1
{\phi(x)} + \frac {|\nabla \phi(x)|^2} {\phi(x)^2} \le 
c (n) \,  \delta^{-2}, \quad x \in B_R(0).
$$
Setting $C(\epsilon)  = c (n) \, \epsilon \, \, \delta(\epsilon)^{-2}$,
 we get :
$$
\left \vert \langle Q u, \, u\rangle \right \vert  \leq \epsilon \, 
||\nabla u||^2_{L^2(\R^n)} + C(\epsilon) \, ||u||^2_{L^2(\R^n)}.
$$
This proves (\ref{E:2.1}). \end{proof}

Let $0<\epsilon_0\le +\infty$, and let 
$C(\epsilon): \, (0, \, \epsilon_0) \to \R^+$.
We define a modified Legendre transform  by:  
\begin{equation}\label{E:2.7a}
\hat C(s) = \inf_{\epsilon\in (0, \, \epsilon_0)} 
\, \{ \epsilon \, s + C(\epsilon) \}, \qquad s>0, 
\end{equation}
Note that $\hat C(s)$ coincides with the  Legendre transform 
of $-C(\epsilon)$. Obviously, $\hat C(s)$ is increasing, concave, and 
 $\hat C (2 s) \le 2 \, \hat C (s)$  ($s>0$). The following statement is an 
immediate consequence of (\ref{E:2.5}).
\begin{corollary}\label{Corollary 2.2} 
Suppose that {\rm (\ref{E:2.1})} holds. Then 
\begin{equation}\label{E:2.8a}
\begin{split}
 \sup_{x_0 \in \R^n} \, & \sup \, \{ \,
\left | \langle Q u, \, u \rangle \right | \, : \,  u
\in C^\infty_0(B_\delta(x_0)), \, \,   ||\nabla u||_{L^2} \le 1
\} \\ & \le \, k(n) \, \delta^2 \, \hat C (\delta^{-2}), 
\end{split}
\end{equation}
where 
$k(n)$ is a constant which depends only on $n$.
\end{corollary}

In the special case where $C(\epsilon) =  c \, \epsilon^{-\beta}$ and 
consequently $\hat C(s) = c_1 \, s^{\frac \beta {\beta + 1}}$,  
we have the following characterization of Trudinger's subordination property. 

\begin{corollary}\label{Corollary 2.3} Suppose $Q \in D'(\R^n)$, $n \ge 2$. Suppose 
$\beta>0$ and $\nu = \frac {n+2} 2 \beta + \frac n 2$. 
Then the following statements are equivalent. 

{\rm (i)} There exist $\epsilon_0>0$ and $c>0$ such that, 
for every  $\epsilon \in (0, \, \epsilon_0)$, 
\begin{equation}\label{E:2.9a} 
\left | \langle Q u, \, u \rangle \right | \le \, \epsilon \, 
||\nabla u||^2_{L^2(\R^n)} + c \,  \epsilon^{-\beta} \, 
||u||^2_{L^2(\R^n)}, \qquad \forall   u \in C^\infty_0(\R^n). 
\end{equation}

{\rm (ii)}  There exist $\epsilon_0>0$ and $c>0$ such that, for every  
$\epsilon \in (0, \, \epsilon_0)$, 
the inequality
\begin{equation}\label{E:2.10a} 
\left | \langle Q u, \, u \rangle \right | \le \, \epsilon \, 
||\nabla u||^2_{L^2(\R^n)} + c \, \epsilon^{-\nu} \, 
||u||^2_{L^1(\R^n)}, \qquad \forall   u \in C^\infty_0(\R^n). 
\end{equation}

{\rm (iii)}  There exist $\delta_0>0$ and $C>0$ such that, for every  $\delta \in 
(0, \, \delta_0)$, 
\begin{equation}\label{E:2.11u} 
\left | \langle Q u, \, u \rangle \right | \le C \, \delta^{\frac 2 {\beta + 1}} \,  
||\nabla u||^2_{L^2(\R^n)} , \qquad \forall   u \in C^\infty_0(B_\delta(x_0)), 
\end{equation}
where $C$ does not depend on 
$x_0$ and $\delta$.  

Moreover, in statements {\rm (i)}$-${\rm (iii)} one may set $\epsilon_0=1$ 
if $\delta_0 = 1$ and vice versa. Similarly, if $\epsilon_0= + \infty$ in {\rm (i)} 
or {\rm (ii)} then 
$\delta_0 = +\infty$ in {\rm (iii)}, and the converse is also true.
\end{corollary}

\begin{proof} Let $\alpha = \frac n {n+2}$. 
Applying Nash's inequality (\cite{LL}, Theorem 2.13),
\begin{equation}\label{E:2.12u} 
||\nabla u||_{L^2(\R^n)} \le C(n) \, 
||\nabla u||^{\alpha}_{L^2(\R^n)} \, ||u||^{1-\alpha}_{L^1(\R^n)}, \qquad 
\forall u \in  C^\infty_0(\R^n), 
\end{equation}
we see that (\ref{E:2.9a}) yields 
 \begin{equation}\label{E:2.13u} 
\left | \langle Q u, \, u \rangle \right | \le C \, \epsilon \, 
||\nabla u||^2_{L^2(\R^n)} + c \, \epsilon^{-\beta} \, 
||\nabla u||^{2 \alpha}_{L^2(\R^n)}  \, 
||u||^{2(1-\alpha)}_{L^1(\R^n)}. 
\end{equation} 
From this by Young's inequality it follows 
\begin{equation}\label{E:2.14u} 
\left | \langle Q u, \, u \rangle \right | \le 2 \epsilon \, 
||\nabla u||^2_{L^2(\R^n)} + c_1 \, \epsilon^{-\nu} \, 
||\nabla u||^{2 \alpha}_{L^2(\R^n)}  \, 
||u||^{2}_{L^1(\R^n)}. 
\end{equation} 
This proves (i)$\Rightarrow$(ii). 

To show that (ii)$\Rightarrow$(iii), suppose 
  $u \in C^\infty(B_\delta(x_0))$. Then 
by Schwartz's inequality and (\ref{E:2.3}),  
\begin{equation}
\begin{split}
\left | \langle Q u, \, u \rangle \right | & \le \epsilon \, ||u||^2_{L^1 (B_\delta(x_0))}
+ c \, \epsilon^{-\nu} \, \delta^n  \, ||u||^2_{L^2 (B_\delta(x_0))} 
\notag \\ & \le 
(\epsilon + c_1 \, \epsilon^{-\nu} \, \delta^{n+2}) \, 
||\nabla u||^2_{L^2(B_\delta(x_0))}.
\notag
\end{split}
\end{equation}
Minimizing over $\epsilon \in (0, \, \epsilon_0)$, we get (iii). 
The implication (iii)$\Rightarrow$(i) is a direct consequence of Lemma~\ref{Lemma 2.1}. 
\end{proof}

\noindent{\bf Remark 3.1.}  Corollary~\ref{Corollary 2.3} covers the case 
$\nu >\frac n 2$ in  (\ref{E:2.10a}), since 
$\nu = \frac {n+2} 2 \beta + \frac n 2$, where 
$\beta>0$.  However, it is easy to see 
that the validity of  (\ref{E:2.10a})  implies $Q=0$ if $\nu < \frac n 2$, and 
$Q \in L^\infty(\R^n)$ if $\nu = \frac n 2$. 
In the latter case, $Q \in L^\infty(\R^n)$ is also sufficient for (\ref{E:2.10a}) 
by Nash's inequality. 

Indeed, suppose that (\ref{E:2.10a}) holds for 
$0< \nu \le \frac n 2$. Then the same 
inequality holds for $Q_t =Q\star\phi_t$ where $\phi_t (x)= t^{-n} \phi (t^{-1} x)$ 
($t>0$) is a standard mollifier. Now using $Q_t$ in place of $Q$,  letting  
$u=\eta_{\delta, \, x_0}$ in (\ref{E:2.10a}), 
where $\eta_{\delta, \, x_0}$ 
is a smooth cut-off function as in the proof of Lemma~\ref{Lemma 2.1}, and  
  minimizing 
 over $\epsilon>0$, we deduce that 
$ \delta^{- n +2 -\frac {n+2} {\nu + 1}} \, \left | \langle Q_t, \, \eta_{\delta, \, x_0}^2 \rangle \right 
| \le C.$ 
Letting $\delta \to 0$, we conclude that $Q_t(x_0)=0$ if $\nu < \frac n 2$, and 
$\sup_{t>0} \, |Q_t(x_0)| \le C $  if $\nu = \frac n 2$, for every $x_0\in R^n$. 
Hence  
$Q=0$ if $\nu < \frac n 2$, and $Q \in L^\infty(\R^n)$ if $\nu= \frac n 2$.

\section{Nonnegative potentials}\label{Section 4}

 In this section, we give necessary and sufficient conditions 
for the infinitesimal form boundedness for nonnegative 
potentials $Q\in D'(\R^n)$, i.e.,  nonnegative measures. 
The general case of distributional potentials  will be considered 
in the next section.

 Let $M^+(\R^n)$ denote the class of  nonnegative locally finite 
Borel measures on $\R^n$. 
 For a compact set $e\subset \R^n$, the capacity 
${\rm cap} \, (e)$ associated with the Sobolev space $W^{1,2}(\R^n)$ 
 is defined by: 
 \begin{equation}\label{E:3.1a}
{\rm cap} \, (e) = \inf \, \left \{ \int_{\R^n} (|\nabla u|^2 + |u|^2) \, dx: 
\quad u \in C^\infty_0(\R^n), \, \, u (x) >1 \, \, {\rm on} \, \, e \right \}.
\end{equation}
On compact sets,  ${\rm cap} \, (\cdot)$ coincides with the  
so-called Bessel capacity defined by (see, e.g., \cite{AH}):  
$$
{\rm cap} \, (e) = \inf \, \left \{ \int_{\R^n} f(x)^2 \, dx: 
\quad f \in L^2 (\R^n), \quad f \ge 0, \quad  G_1\star f (x) \ge 1 
 \, \, {\rm on} \, \, 
e \right \},
$$
where $G_1 \star f = (1-\Delta)^{- \frac 1 2} f$ is the 
Bessel potential of order $1$.

\begin{theorem}\label{Theorem 3.1a} Suppose $\mu \in M^+(\R^n)$, $n \ge
2$. Then the following statements are equivalent.

{\rm (i)} There exists $\epsilon_0>0$ and a function 
$C(\epsilon) : \, (0, \epsilon_0) \to \R_+$ such that 
the inequality
\begin{equation}\label{E:3.2a}
\int_{\R^n} |u(x)|^2 \, d \mu \leq \epsilon \, ||\nabla
u||^2_{L^2} + C(\epsilon) \, ||u||^2_{L^2}, \quad \forall u  \in
C^\infty_0(\R^n), 
\end{equation}
holds for every $\epsilon\in(0, \, \epsilon_0)$. 

{\rm (ii)} 
\begin{equation}\label{E:3.3a}
\lim_{\delta \to +0} \, \sup_{x_0 \in \R^n} \, \sup \, \left \{ \, 
\int_{\R^n} |u(x)|^2 \, d \mu \, : \,  u
\in C^\infty_0(B_\delta(x_0)), \, \,   ||\nabla u||_{L^2} \le 1 
\right \} = 0. 
\end{equation}

{\rm (iii)} 
\begin{equation}\label{E:3.4a}
\lim_{\delta \to +0} \, \sup \, \left \{ \, \frac {\mu(e)}{{\rm cap} \, (e)}: 
\quad e\subset \R^n, \, \,  {\rm diam} \, (e) \le \delta \right \} = 0,
\end{equation}
where the supremum above is  over  compact sets $e$ 
of positive capacity. 

{\rm (iv)} 
\begin{equation}\label{E:3.5a}
\lim_{\delta \to +0} \, \sup \, \left \{ \, \frac 1 {\mu(P_0)} \, 
 \displaystyle{\sum_{P\subseteq P_0}} \, 
|P|^{\frac 2 n -1} \, \mu(P)^2  \, : 
\quad P_0 \subset \R^n, \, \, 
{\rm diam} \, (P_0) \le \delta \right \} =0, 
\end{equation}
where the supremum is over dyadic cubes $P_0$ such that 
$\mu(P_0)\not=0$, and the sum is over all 
dyadic cubes $P$ contained in $P_0$.  
\end{theorem}

Obviously, in statement (i) one can  assume that $C(\epsilon)$ is defined on 
$\R^+=(0, +\infty)$ so that  (\ref{E:3.2a}) holds for every $\epsilon>0$.  

\begin{proof} By Lemma \ref{Lemma 2.1}, (i)$\Leftrightarrow$(ii).  
We next show that 
(ii)$\Leftrightarrow$(iii).

Let $\Omega \subset \R^n$ be an open set, and let $e$ be a compact subset 
of $\Omega$. We will need the  the Wiener capacity relative 
to the domain $\Omega$   
defined by: 
\begin{equation}\label{E:3.6a}
{\rm cap} \, (e, \Omega) = \inf \, \left  \{ \int_{\Omega} |\nabla u|^2  \, dx: 
\quad u \in C^\infty_0(\Omega), \, \, u >1 \, \, {\rm on} \, \, e \right \}.
\end{equation}
We remark that, for $\Omega=\R^n$, in the case $n\ge 3$ and ${\rm diam} \, (e) \le 1$,  
we have: 
\begin{equation}\label{E:3.7a}
 {\rm cap} \, (e, \R^n) \le  {\rm cap} \, (e) \le c(n) \, {\rm cap} \, (e, \R^n),
\end{equation}
where $c(n)$ depends only on $n$. The left-hand side estimate is obvious, and 
the right-hand side follows, e.g.,  from Hardy's inequality. 
For $n=1, 2$ it is easy to see that 
${\rm cap} \, (e, \R^n) = 0$ for every $e\subset \R^n$ (see \cite{M2}).

We now set 
\begin{equation}\label{E:3.8a}
c_1(\mu, \Omega) = \sup \, \left \{ \, 
\int_{\R^n} |u(x)|^2 \, d \mu \, : \,  u
\in C^\infty_0(\Omega), \, \,   ||\nabla u||_{L^2(\Omega)} \le 1 
\right \}
\end{equation}
and 
\begin{equation}\label{E:3.9a}
c_2(\mu, \Omega) = \sup \, \left \{ \, \frac {\mu(e)}{{\rm cap} \, (e, \Omega)}: 
\quad e\subset \Omega\right\},
\end{equation}
where the supremum above is over compact sets $e\subset \Omega$ of positive capacity. 
As was shown in \cite{M1} (see also \cite{M2}, Sec. 2.5),  
\begin{equation}\label{E:3.10a}
c_2(\mu, \Omega) \le c_1(\mu, \Omega) \le 4 \,c_2(\mu, \Omega).
 \end{equation}

Applying this inequality with $\Omega = B_\delta(x_0)$, 
we see that (\ref{E:3.3a}) holds if and only if 
\begin{equation}\label{E:3.11a}
\lim_{\delta \to +0} \sup_{x_0\in \R^n} 
\sup \, \left \{ \, \frac {\mu(e)}{{\rm cap} \, (e, B_\delta(x_0))}: 
\quad e\subset B_\delta(x_0) \right \} = 0,
\end{equation}
where $e$ is a compact subset of $B_\delta(x_0)$.

To verify (ii)$\Leftrightarrow$(iii), it remains to prove that 
one can replace 
the capacity ${\rm cap} \, (e, B_\delta(x_0))$ in (\ref{E:3.11a})  
with ${\rm cap} \, (e)$  
 where $e$ is a compact set such that  ${\rm diam} \, (e) \le \delta$.  
 Without loss of generality we 
may assume that $0<\delta\le 1$.

Notice that if $e\subset B_\delta(x_0)$ and $\delta\le 1$, then, for every 
$u \in C^\infty_0(B_\delta(x_0))$,   
$$
\int_{\R^n}  (|\nabla u|^2 +|u|^2) \, dx \le 
\int_{\R^n} ( |\nabla u|^2 + {\delta}^{-2} |u|^2) \, dx
\le c(n) \, \int_{\R^n} |\nabla u|^2 \, dx,
$$
by (\ref{E:2.3}). Hence,  minimizing both sides over  $u$
such that 
$u(x) > 1$ on $e$, we get:
$$
{\rm cap} \, (e) \le c(n) 
 \, {\rm cap} \, (e, B_\delta(x_0)).
$$
This estimate, which  holds for every $n \ge 1$,  yields (iii)$\Rightarrow$(ii). 

In the opposite direction, suppose first that  $n \ge 3$. Then we will show that 
\begin{equation}\label{E:3.12a}
{\rm cap} \, (e, B_{\widetilde \delta}(x_0)) \le c \, {\rm cap} \, (e),
\end{equation}
for every $e\subset  B_\delta(x_0)$, where $\widetilde \delta = 2 \delta$, and 
$c$ depends only on $n$.  

Indeed, suppose $u \in C^\infty_0(\R^n)$, where 
$u(x) > 1$ on $e$, and  
$$
\int_{\R^n} (|\nabla u|^2 + |u|^2) \, dx \le 2 \, {\rm cap} \, (e).
$$
Denote by $\eta_{2\delta, x_0} (x) = \eta \left ( \frac{x-x_0}{2 \delta}\right )$ a 
smooth cut-off function 
supported in $B_{2\delta(x_0)}$ so that $ \eta_{2\delta, x_0} (x) =1$ on $B_\delta(x_0)$. 
Then $v(x) = \eta_{2\delta, x_0}(x) \, u(x) \in C^\infty_0(B_{2\delta}(x_0))$, 
and $v(x)>1$ on $e$. We estimate using Hardy's inequality:
\begin{align}
\int_{B_{2\delta}(x_0)} |\nabla v|^2 \, dx & \le 2 \int_{B_{2\delta}(x_0)} 
( |\eta_{2\delta, x_0}|^2 |\nabla u|^2 + 
|\nabla \eta_{2\delta, x_0}|^2 |u|^2 ) \, dx \notag \\ & \le 
c  \int_{B_{2\delta}(x_0)} (|\nabla u|^2 
+  \, \delta^{-2} |u|^2) \, dx \le c_1 \, \int_{\R^n} |\nabla u|^2 \, dx,
\notag
\end{align}
 where $c$, $c_1$ depend only on $n$. Thus,
$$
{\rm cap} \, (e, B_{2 \delta}(x_0))  \le \int_{B_{2\delta}(x_0)} |\nabla v|^2 \, dx 
 \le c_1    \int_{\R^n} |\nabla u|^2 \, dx   \le  \, 2 c_1 \, \, {\rm cap} \, (e),
$$
which proves (\ref{E:3.12a}) with $\widetilde \delta = 2 \delta$.  
From this it is immediate that (ii)$\Rightarrow$(iii) if $n\ge 3$.

In the case $n=2$ we use a more delicate argument with $\widetilde\delta \in
 (\delta, \, 1)$ to be 
determined later.  
We denote by $\nu = \nu_e \in M^+(\R^2)$  
the equilibrium measure associated with a compact set $e\subset \R^2$  
such that (see \cite{AH}, \cite{M2}):
\begin{equation}\label{E:3.13a}
\nu (e) = {\rm cap} \, (e) = c(n) \int_{\R^2} \left (  |\nabla \P \nu|^2 + 
(\P \nu)^2| \right ) \, dx, \qquad {\rm supp} \, (\nu) \subset e.
\end{equation} 
Here $\P \nu = (1-\Delta)^{-1} \nu =G_2 \star \nu$ is the Bessel potential of 
order $2$: 
\begin{equation}\label{E:3.14a}
\P \nu (x) = \int_{\R^2} G_2(x-y) \, d \nu(y), 
\end{equation} 
where the Bessel kernel $G_\alpha$ of order $\alpha$ on $\R^n$ is defined through the Fourier 
transform: 
\begin{equation}\label{E:3.15a}
\widehat G_\alpha \, (\xi) = (1+|\xi|^2)^{-1}, \quad \xi \in \R^n, \quad \alpha>0.
\end{equation}    
 
Without loss of generality we may assume  $\P \nu \in C^\infty (\R^2)$. 
(This assumption, which is not essential to the proof,    
is easily removed by using an appropriate mollifier $\P \nu \star \zeta_r$, where 
$\zeta_r (x) = r^{-2} \zeta(\frac x r)$, and letting $r\to +\infty$.)

Let $v(x) = \eta_{\widetilde \delta} (x) \, \P \nu (x) 
\in C^\infty_0 (B_{\widetilde\delta} (x_0))$. 
It follows that 
$$
{\rm cap} \, (e, B_{\widetilde \delta}(x_0)) \le 
\int_{B_{\widetilde \delta} (x_0)} |\nabla v|^2 \, dx \le 
\int_{B_{\widetilde \delta} (x_0)} |\nabla \P \nu|^2 \, dx + c \, \widetilde\delta^{-2} 
\int_{B_{\widetilde \delta}(x_0)}   (\P \nu)^2 \, dx.
$$
The first term on the right is bounded by ${\rm cap} \, (e)$. We rewrite the 
integral in the second term 
in the form (note that $\nu$ is supported on $e$): 
$$
\int_{B_{\widetilde \delta}(x_0)}   (\P \nu)^2 \, dx = 
\int_{e} \int_{e} \int_{B_{\widetilde \delta}(x_0)} 
 G_2(x-y) \, G_2(x-y') \, dx \, d \nu(y) 
\, d \nu(y'),
$$
and estimate using the fact that 
 $G_2(x) \le c \,  \log \frac 4 {|x|}$ for $|x|\le 2$ in the case $n=2$ (see, e.g., 
\cite{AH}, Sec. 1.2.5).  Since $e\subset B_{\delta}(x_0)$ and $
\max \, (|x-y|, |x-y'|) < 2 \widetilde \delta$, 
we have: 
\begin{align}
& \int_{B_{\widetilde \delta}(x_0)}  G_2(x-y) \, G_2(x-y') \, dx    \le 
\int_{B_{\widetilde \delta}(x_0)}  (G_2(x-y)^2 + G_2(x-y')^2) \,  dx \notag\\
& \le c \int_{|x-y|< 2 \widetilde \delta} \log^2 \frac 4 {|x-y|} \, dx  + 
c \int_{|x-y'|< 2 \widetilde \delta} \log^2 \frac 4 {|x-y'|} \, dy' 
 \le c \, \widetilde\delta^2 \, \log^2 \, \frac 2 {\widetilde \delta}.
\notag
\end{align}
It follows:
$$
\int_{B_{\widetilde \delta}(x_0)}   (\P \nu)^2 \, dx \le c \,
\nu(e)^2 \, \widetilde\delta^2 \, \log^2 \, \frac 2 {\widetilde \delta} 
= c \, {\rm cap} \, (e)^2 \, \widetilde\delta^2 \, \log^2 \, \frac 2 {\widetilde \delta}.
$$
Hence, 
\begin{equation}\label{E:3.16a}
{\rm cap} \, (e, B_{\widetilde \delta}(x_0)) \le c \, 
\left ( {\rm cap} \, (e) + 
\log^2 \, \frac 2 {\widetilde \delta} \, \, {\rm cap} \, (e)^2 \right ),
\end{equation} 
where $e\subset B_{\delta}(x_0)$ and $c$ depends only on $n$. 
Using the known estimate of the capacity of a ball (\cite{M2}, Sec. 7.2.3):
\begin{equation}\label{E:3.17a}
{\rm cap} \, (B_\delta(x_0)) \asymp \frac 1  {\log \, \frac 2 \delta}, \qquad  
0<\delta<1, \qquad n=2,
\end{equation}  
we have: 
$$
{\rm cap} \, (e) \le {\rm cap} \, (B_\delta(x_0)) \le \,  
\frac c {\log \, \frac 2 \delta}.
$$
Combining the preceding estimate and (\ref{E:3.16a}) gives:
\begin{equation}\label{E:3.18a}
  {\rm cap} \, (e, B_{\widetilde \delta}(x_0)) \le c \, {\rm cap} \, (e) \, 
\left ( 1 + \frac {\log^2 \, \frac 2 {\widetilde \delta}} 
{\log \, \frac 2 \delta}\right), \quad e \subset B_\delta(x_0), \quad n=2.
\end{equation}  

Now choosing $\widetilde \delta$ so that 
$\log^2 \,  \frac 2 {\widetilde \delta}  = 
\log \,  \frac 2 {\delta} $,
we get (\ref{E:3.12a}). Letting $\delta\to+0$, so that 
$\widetilde \delta\to+0$, we obtain that (ii)$\Rightarrow$(iii) 
in the case 
$n=2$. This proves (ii)$\Leftrightarrow$(iii).

 To prove (iii)$\Leftrightarrow$(iv), 
we set 
\begin{equation}\label{E:3.19a}
\kappa_1 (\mu, \delta) = \sup \, \left \{ \frac {\mu(e)}{{\rm cap} \, (e)} : \quad 
{\rm diam} \, (e) \le \delta \right\},
\end{equation}
where $e$ is a compact set of positive capacity in $\R^n$, 
and 
\begin{equation}\label{E:3.20a}
\kappa_2 (\mu, \delta) = \sup \left \{  \frac 1 {\mu(P_0)} 
\displaystyle{\sum_{P\subseteq P_0}} \, 
|P|^{\frac 2 n -1} \, \mu(P)^2 \, : 
\quad 
{\rm diam} \, (P_0) \le \delta  \right\},
\end{equation}
where the supremum is over dyadic cubes $P_0$  in $\R^n$ such that 
$\mu(P_0) \not= 0$.  

For a set 
$e\subset\R^n$,  denote by $\mu_{e}$ the restriction of the measure 
$\mu$ to $e$, i.e., $d \mu_{e} = \chi_{e} \, d \mu$. Denote by 
$$
\mathcal E (\mu) = ||\mu||^2_{W^{-1,2}(\R^n)} = \int_{\R^n} (G_1\star \mu)^2 \, dx
$$
the energy of the measure $\mu$. 

We first show that, for any compact set $e$ such that 
 ${\rm diam} \, (e) \le \delta$,   
\begin{equation}\label{E:3.21a}
\mathcal E (\mu_{e}) \le c \,  \kappa_1 (\mu, \delta) \,  \mu(e).
\end{equation}
Such inequalities, without the restriction on the diameter of $e$, are known 
(see, e.g., \cite{M2}, \cite{V2})   
and we will be brief here. Clearly, there exists  $g\in C^\infty_0(\R^n)$,
$g \ge 0$,  such 
that $||g||_{L^2(\R^n)} = 1$, and 
$$
\mathcal E (\mu_{e})^{\frac 1 2} \le 2 \int_{\R^n} G_1\star \mu_e \, g \, dx = 
2 \int_{e} G_1\star  g \, d \mu.
$$
For any $\lambda>0$, set $e_\lambda = \{ x: \, G_1\star  g(x)>\lambda\}$. Then 
$$
\mu_e (e_\lambda) = \mu(e_\lambda\cap e) \le  \kappa_1 (\mu, \delta) \, {\rm cap}
 (e_\lambda\cap e) \le c \, \kappa_1 (\mu, \delta) \, \left \Vert \frac g \lambda 
\right \Vert^2_{L^2(\R^n)},
$$
where the last estimate follows from the definition of the capacity since 
$$
u(x) = \left( G_1\star  \frac g \lambda\right) (x)>1, \quad x \in e_\lambda,
$$
and $||u||^2_{W^{1,2}(\R^n)} = c \, ||g||^2_{L^2(\R^n)} =c$, where 
$c$ depend only on $n$. 

In other words, $G_1\star  g\in L^{2, \infty} (\mu_e)$, and 
$$
||G_1\star  g||_{L^{2, \infty} (\mu_e)} \le c \,
 \kappa_1 (\mu, \delta)^{\frac 1 2},
$$
 where 
$L^{2, \infty} (\nu)$ denotes the corresponding Lorentz (weak $L^2$) space 
on $\R^n$ equipped with the measure $\nu$. Hence, 
$$
\int_{e} G_1\star  g \, d \mu \le c \, 
||G_1\star  g||_{L^{2, \infty} (\mu_e)} \, 
\mu(e)^{\frac 1 2} \le c \, \kappa_1 (\mu, \delta)^{\frac 1 2} 
\, \mu(e)^{\frac 1 2},
$$
which proves (\ref{E:3.21a}). 

Letting $e=P_0$ in (\ref{E:3.21a}) 
where $P_0$ is a dyadic cube such that ${\rm diam} \, (P_0) 
\le \delta$, we get
$$
\mathcal E (\mu_{P_0}) \le \kappa_1 (\mu, \delta) \, \mu(P_0).
$$

For $x, \, y \in P_0$, we have 
$|x-y|\le \delta \le 1$, and it is easy to see  that 
$$
 \int_{P_0} G_1(x-y) \, d \mu(y) \ge \int_{P_0} 
\frac {d \mu(y)} {|x-y|^{n-1}} \ge c \, \sum_{P\subseteq P_0} 
|P|^{\frac 1 n -1} \mu(P) \, \chi_P (x).
$$
We estimate:
\begin{align}
& \mathcal E (\mu_{P_0})  \ge c \int_{P_0} (G_1\star \mu_{P_0})^2 \, dx  
\ge  c \, \int_{P_0} \left (\sum_{P\subseteq P_0} 
|P|^{\frac 1 n -1} \mu(P) \, \chi_P (x) \right)^2 dx \notag 
\\ & \ge \int_{P_0} \sum_{P\subseteq P_0} |P|^{\frac 2 n -2} \mu(P)^2 
\chi_P (x) \, dx \notag  = 
\sum_{P\subseteq P_0} |P|^{\frac 2 n -1} \mu(P)^2. \notag 
 \end{align}
Thus, 
$$
\sum_{P\subseteq P_0} |P|^{\frac 2 n -1} \mu(P)^2 \le c \, \kappa_1 
(\mu, \delta) \, \mu(P_0),
$$
that is, $\kappa_2 (\mu, \delta) \le c \, \kappa_1 (\mu, \delta)$.  
This estimate yields (iii)$\Rightarrow$(iv). 

The converse can be proved as in \cite{V2} 
in the case of Riesz kernels using Th. Wolff's inequality \cite{HW}. (See also \cite{COV} where 
Wolff's inequality is proved for general 
dyadic and radial kernels.) Here we sketch a direct proof 
based on  the dyadic Carleson measure theorem (see, e.g., \cite{V1})  
 which yields the inequality 
$$
\sum_{P \subseteq P_0} \, |P|^{\frac 2 n -1} \mu(P)^2 
\left ( \frac 1 {\mu(P)} \int_P g \, d \mu \right )^2 
\le c \, \kappa_2 (\mu, \delta)  \, 
\int_{P_0} g^2 \, d \mu,
$$
where $g \in L^2(\mu_{P_0})$, $g\ge 0$, 
for every dyadic cube $P_0$ such that ${\rm diam} \, (P_0) \le \delta$. 
(This inequality follows by interpolation between the trivial $L^\infty(\mu_{P_0})\to l^\infty$ 
 estimate 
and the weak-type $(1,1)$ estimate from $L^1(\mu_{P_0})$ to 
$l^{1, \infty}$. Note that the left-hand side of the preceding 
inequality is bounded from below by 
$$
\int_{P_0} \left ( \sum_{P \subseteq P_0} \, |P|^{\frac 1 n -1} 
\chi_P (x) \int_P g \, d \mu \right )^2 \, dx.
$$
To verify this estimate, which is closely related to Wolff's inequality, 
we use the obvious 
pointwise inequality: 
\begin{align}
& \left ( \sum_{P \subseteq P_0} \, |P|^{\frac 1 n -1} 
\chi_P (x) \int_P g \, d \mu \right )^2 \notag \\ & \le 2 \, 
\sum_{P \subseteq P_0}  |P|^{\frac 1 n -1} \left (\int_P g \, d \mu \right) 
\, \sum_{P' \subseteq P} |P'|^{\frac 1 n -1} \, \chi_{P'}(x) \, 
 \left ( \int_{P'} g \, d \mu
\right ). 
\notag
\end{align}
Integrating both sides of the preceding inequality 
 over $P_0$, and  the estimate 
$$
\sum_{P' \subseteq P} |P'|^{\frac 1 n} \, 
  \left (\int_{P'} g \, d \mu \right) 
 \le c(n) \, |P|^{\frac 1 n} \,  \int_{P} g \, d \mu,
$$
 we obtain 
\begin{align}
\int_{P_0} \left ( \sum_{P \subseteq P_0} \, |P|^{\frac 1 n -1} 
\chi_P (x) \int_P g \, d \mu \right )^2  dx & \le c \, 
\sum_{P \subseteq P_0} \, |P|^{\frac 2 n -1} 
\left (\int_P g \, d \mu \right )^2 \notag
 \\ & \le c \,
 \kappa_2 (\mu, \delta)  \, 
\int_{P_0} g^2 \, d \mu. 
\notag
\end{align}
By duality, this gives 
$$
\int_{P_0} \left ( \sum_{P \subseteq P_0} \, |P|^{\frac 1 n -1} 
\chi_P (x) \int_P f (y) \, d y \right )^2  d \mu(x) \le c \,
 \kappa_2 (\mu, \delta)  \, 
\int_{P_0} f(y)^2 \, d y, 
$$
for every $f \in L^2(dx)$, $f \ge 0$. Note that 
$$
I^{P_0}_1 f(x) = \sum_{P \subseteq P_0} \, |P|^{\frac 1 n -1} 
\chi_P (x) \int_P f(y) \, d y
$$
is the dyadic Riesz potential of order $1$ (see \cite{HW}) scaled to the 
cube $P_0$. It is easy to see that this estimate, which holds for every 
dyadic cube $P_0$ such that ${\rm diam} \, (P_0) \le \delta \le 1$, implies 
the inequality 
$$
\int_{P_0}  (G_1\star f)^2 \, d \mu \le c \, \kappa_2 (\mu, \delta)  \, 
||f||^2_{L^2(dx)}, 
$$
for every $f\in L^2(dx)$. (See details in \cite{V2} or \cite{COV}.)

 From this and the definition of the Bessel capacity it is immediate that 
$$
\mu(e) \le c \, \kappa_2 (\mu, \delta) {\rm cap} \, (e),
 \quad {\rm diam} \, (e) \le \delta. 
$$
Hence, $
\kappa_1 (\mu, \delta)  \, \le c \, \kappa_2 (\mu, \delta),
$ 
where $c$ depends only on $n$.   From this it follows that 
(iv)$\Rightarrow$(iii). This completes the proof of Theorem \ref{Theorem 3.1a}. 
\end{proof}

\begin{corollary} \label{Corollary 3.2a} Let $\mu \in M^+(\R^n)$. 
If  inequality {\rm (\ref{E:3.2a})}  holds then 
\begin{align}\label{E:3.22a}
& \lim_{\delta \to +0} \, \sup \, \left \{ \, \frac {\mu(B_\delta(x_0))}
{\delta^{n-2}}: \quad 
x_0 \in  \R^n  \right \} = 0, \quad n \ge 3, \\ 
& \lim_{\delta \to +0} \, \sup \, \left\{ \,  \log  \tfrac {1} {\delta}  \, \,  \mu(B_\delta(x_0)) \, : \quad 
x_0 \in  \R^n  \right \} = 0, \quad n =2.\label{E:3.22b}
\end{align}
\end{corollary}

Corollary 3.2  follows by letting $e = B_\delta(x_0)$  in Theorem 
\ref{Theorem 3.1a} (iii),  
and using the estimates  ${\rm cap} \, (B_\delta(x_0)) \asymp \delta^{n-2}$ 
if $n\ge 3$ 
and ${\rm cap} \, (B_\delta(x_0)) \asymp \frac 1 {\log \, \frac 2 \delta }$ 
if $n=2$, provided $0<\delta \le 1$ (see \cite{M2}).

Stronger necessary conditions are obtained 
by replacing $\mu (B_\delta(x_0))$ above with 
$\int_{B_\delta(x_0)} (G_1 \star \mu)^2 \, dx$.

\begin{corollary} \label{Corollary 3.3a} Let $\mu \in M^+(\R^n)$. 
If  inequality {\rm (\ref{E:3.2a})}  holds then 
\begin{align}
& \lim_{\delta \to +0} \, \sup \, \left \{ \, \frac { \int_{B_\delta(x_0)} 
(G_1\star \mu)^2 \, dx }
{\delta^{n-2}}: \quad 
x_0 \in  \R^n  \right \} = 0, \quad n \ge 3, \label{E:3.22c}\\ 
& \lim_{\delta \to +0} \, \sup \, 
\left \{ \, \log  \tfrac {1} {\delta}   \,  
 \, \int_{B_\delta(x_0)} 
(G_1\star \mu)^2 \, dx: \quad 
x_0 \in  \R^n  \right \} = 0, \quad n =2.\label{E:3.22d}
\end{align}
\end{corollary} 

This corollary is a consequence of the results of \cite{MV2} (
see also Sec. \ref{Section 5}
 below).  

Various related results can be deduced from Theorem \ref{Theorem 3.1a} using 
known form boundedness criteria (see, e.g., \cite{Fef}, 
\cite{ChWW}, \cite{KS}, 
  \cite{M2}, \cite{MV1}, \cite{MV2}, \cite{V2}). In particular, 
we can give more equivalent 
conditions which are necessary and sufficient for (\ref{E:3.2a}).

\begin{corollary} \label{Corollary 3.4a} Let $\mu \in M^+(\R^n)$. 
Then inequality {\rm (\ref{E:3.2a})}  
is equivalent to any one of the following conditions:

\noindent {\rm (i)} 
$$
 \lim_{\delta \to +0} \, \sup \, \left \{ \, \frac { {\mathcal E}  
\left (\mu_{B_\delta (x_0)}\right )}{ \mu (B_\delta(x_0))} : \quad 
x_0 \in  \R^n  \right \} = 0.
$$
 {\rm (ii)} 
$$
 \lim_{\delta \to +0} \, \sup \, \left \{ \, \frac {\left \Vert 
\chi_{B_\delta(x_0)} \, \mu \right \Vert^2_{W^{-1,2} (\R^n)}}
{ \mu (B_\delta(x_0))} : \quad 
x_0 \in  \R^n  \right \} = 0.
$$
{\rm (iii)} 
$$
 \lim_{\delta \to +0} \, \sup \, 
\left \{ \, \frac { G_1 \star \left ( G_1 \star  
\mu_{B_\delta(x_0)} \right )^2 (x)}{G_1 \star \mu_{B_\delta(x_0)}(x) } : 
\quad x, \, x_0 \in  \R^n  \right \} = 0.
$$
\end{corollary} 

Statements (i) and (ii) of Corollary \ref{Corollary 3.3a} follow directly 
from the proof of Theorem~\ref{Theorem 3.1a}  
while (iii) can be proved in a similar way using results of \cite{MV1}. 

Simpler sufficient conditions can be deduced easily from Theorem~\ref{Theorem 3.1a}.   
For instance, it is well known that if $d \mu = \rho(x) \, dx$ where 
$\rho \in L^{\frac n 2} (\R^n) + L^\infty(\R^n)$, $n \ge 3$, then 
{\rm (\ref{E:3.2a})}  holds. Moreover, it can be replaced with the ``locally 
uniform'' $L^{\frac n 2}$-condition (see \cite{BrK}):
$$\lim_{\delta\to +0} \, \sup_{x_0 \in \R^n} \, ||\rho||_{L^{\frac n 2} 
(B_\delta(x_0))} = 0.
$$
In fact, one can use here a local 
weak-$L^{\frac n 2}$, and even  a local Fefferman-Phong norm.

\begin{corollary} \label{Corollary 3.4b} Let $d \mu = \rho(x) \, dx$, where 
$\rho\in L^r_{{\rm loc}} (\R^n)$ for some $r>1$ and $n \ge 3$.  
Then  {\rm (\ref{E:3.2a})}  holds if
\begin{equation}\label{E:3.22e}
 \lim_{\delta \to +0} \, \sup_{x_0\in \R^n}  \, \frac {\int_{B_\delta(x_0)} 
\rho(x)^r \, dx}
{\delta^{n-2r}} = 0.
\end{equation}
Furthermore,   this condition can be replaced with a weaker one, with
$(G_1 \star \mu)^2$  in place of $\rho$.
\end{corollary}

Another well-known sufficient condition for {\rm (\ref{E:3.2a})} 
which is generally not covered by the preceding corollary is Kato's  
condition $K_n$:
\begin{align}
\lim_{\delta \to +0} \, \sup_{x_0\in\R^n} \, & \int_{B_\delta(x_0)} 
 \frac {\rho(x)} {|x-x_0|^{n-2}}
 \, dx
 = 0, \qquad n \ge 3,  \label{E:3.23b}\\
\lim_{\delta \to +0} \, \sup_{x_0\in\R^n} \, & \int_{B_\delta(x_0)} 
\log  \frac 1 {|x-x_0|}  \, \rho(x)
 \, dx 
 = 0, \qquad n =2, \label{E:3.23c}
\end{align}
where $\rho \ge 0$ and  $\rho \in L^1_{{\rm loc}} (\R^n)$. 

\noindent{\bf Remark 4.1.} This class can be broadened  in the same way 
as above by replacing $\rho$ with $(G_1\star \mu)^2$ and 
using Theorem~\ref{Theorem 3.1a}.

In the opposite direction, the following remarkable fact was 
proved in \cite{AiSi} using a probabilistic argument: 
Kato's condition $K_n$ 
follows from {\rm (\ref{E:3.2a})} provided 
\begin{equation}\label{E:3.23d}
C(\epsilon) \le a \, e^{b \epsilon^{-p}}, 
\qquad \epsilon \in (0, \, \epsilon_0), \qquad 0<p<1.
\end{equation}
Such results  can also be deduced by operator semigroup methods \cite{LPS}. 
In fact,  the same conclusion has been 
obtained  under a more general assumption 
in terms of the Legendre transform $\hat C(s)$ defined by (\ref{E:2.7a}): 
\begin{equation}\label{E:3.23e}
\int_{\delta_0}^{+\infty} \frac {\hat C(s)} {s^2} \, ds < +\infty,
\end{equation} 
for some $\delta_0>0$. (See \cite{Gr1}, \cite{Gr2}; an analogous fact is now known for general  
elliptic operators with good  heat 
kernel bounds \cite{Dav1}.) It is worth noting that, for a decreasing $C(\epsilon)$ such that $
\lim_{\epsilon\to +0} C(\epsilon) = +\infty$, the preceding condition is equivalent 
to (\cite{Ko}, Sec. VII.D.2):
\begin{equation}\label{E:3.23ko}
\int_0^{\epsilon_0} \log {C(\epsilon)}  \, d \epsilon < +\infty,
\end{equation} 
for some $\epsilon_0>0$. 

We should emphasize that such  results can be deduced directly 
from Corollary~\ref{Corollary 2.2} in a sharper form. In particular,  for $n=2$, 
 it suffices to assume that $C(\epsilon) \le a \, e^{b \epsilon^{-p}}$ 
 for {\it any} $p>0$.

\begin{proposition} \label{Proposition 3.5a} 
Let $d \mu = \rho(x) \, dx$, where $\rho \ge 0$, and 
$\rho \in L^1_{{\rm loc}} (\R^n)$. 
Suppose that 
{\rm (\ref{E:3.2a})} holds  with $C(\epsilon)$ 
obeying condition {\rm (\ref{E:3.23e})} if $n\ge 3$, or 
 \begin{equation}\label{E:3.23f}
\int_{\delta_0}^{+\infty} \frac {\hat C(s)} {s^2 \log s} \, ds < +\infty, \qquad n=2,
\end{equation}
for some $\delta_0>1$. Then Kato's condition $K_n$ is valid. Moreover, 
\begin{align}
& \sup_{x_0\in\R^n} \, \int_{B_\delta(x_0)} 
 \frac {\rho(x)} {|x-x_0|^{n-2}}
 \, dx \le c \, 
\int_{\delta^{-2}}^{+\infty} \frac {\hat C(s)} {s^2} \, ds, 
\qquad n \ge 3, \label{E:3.23g} \\
& \sup_{x_0\in\R^n} \, \int_{B_\delta(x_0)} \log \frac 1 {|x-x_0|} \, 
\rho(x) 
 \, dx \le c \, 
\int_{\delta^{-2}}^{+\infty} \frac {\hat C(s)} {s^2 \log s} \, ds, 
\qquad n =2,
\label{E:3.23h}
\end{align}
where $c >0$ is a constant which depends only on $n$, and $\delta$ is sufficiently small.  
\end{proposition}
 
\begin{proof} Suppose that {\rm (\ref{E:3.23f})} holds. By 
Corollary~\ref{Corollary 2.2}, 
$$\int_{B_\delta(x_0)}  |u(x)|^2 \, d \mu  \le 
k \, \delta^2 \, \hat C(\delta^{-2}) \, ||\nabla u||^2, \qquad 
\forall u \in C^\infty_0(B_\delta(x_0)), 
$$
where $k$ depends only on $n$. Let $u>1$ on $B_{\frac \delta 2} 
(x_0)$, where $u \in C^\infty_0(B_\delta(x_0))$. Taking the infimum over all 
 such test functions $u$ in the preceding estimate, 
we obtain:
$$\mu \, ( B_{\frac \delta 2}(x_0) ) \le  \,
c \, \delta^2 \, \hat C(\delta^{-2}) \, \, 
{\rm cap} \,  (B_{\frac \delta 2}(x_0), \, B_{\delta}(x_0)),
$$ 
where ${\rm cap} \, (\cdot, \, B_{\delta}(x_0))$ is the capacity (\ref{E:3.6a}). 
Using the estimates (see \cite{M2}):
\begin{align}
& {\rm cap (B_{\frac \delta 2}(x_0), \, B_{\delta}(x_0))} \asymp 
\delta^{n-2}, \qquad n \ge 3,  \notag \\ 
& {\rm cap (B_{\frac \delta 2}(x_0), \, B_{\delta}(x_0))} \asymp \, (\log 
\tfrac {2} {\delta})^{-1}, \qquad n =2, \notag
\end{align}
where $0<\delta<1$, together with the inequality $\hat C (2s) \le 2 \, \hat C(s)$ ($s>0$), we deduce:
\begin{align}
& \mu (B_{\delta}(x_0) )  \le c \,
 \, \delta^n \, \hat C(\delta^{-2}), \qquad n\ge 3, \notag \\ 
 & \mu( B_{\delta}(x_0) )  \le c \,
 \, \delta^2 \, \hat C(\delta^{-2}) \, \left (\log 
\frac 2 \delta \right)^{-1}, \qquad n=2.
\notag
\end{align}

Suppose $d \mu = \rho(x) \, dx$ and $n \ge 3$. We estimate: 
\begin{align}
\int_{B_{\delta}(x_0)} \frac{\rho (x)} {|x-x_0|^{n-2}}  \, dx
& \le c \, \int_0^\delta \frac{\mu(B_{s}(x_0))} {s^{n-2}} \,  \frac {ds} {s}
 \notag \\ & \le 
c \, \int_0^\delta  s \, \hat C(s^{-2}) \, ds.
\notag
\end{align}
In the case $n=2$ it follows:
\begin{align}
\int_{B_{\delta}(x_0)} \log \, \frac 1 {|x-x_0|}  \rho(x) \, dx
& \le c \int_0^\delta \mu(B_{s}(x_0)) \,  \frac {ds} {s} +  c \, 
\frac {\mu(B_{\delta}(x_0))}{\delta} \notag \\
& \le   c_1 \, \int_0^\delta \frac{s \, \hat C(s^{-2})} {\log s^{-1}} \,  ds.
\notag
\end{align}
Obviously, these estimates yield (\ref{E:3.23g}) and (\ref{E:3.23h}).
\end{proof}

We conclude this section with a characterization of the  subordination 
inequality 
for nonnegative measures $\mu$: 
\begin{equation}\label{E:3.24a}
\int_{\R^n} |u(x)|^2 \, d \mu \le \epsilon \, ||\nabla u||^2_{L^2(\R^n)} + 
\frac C {\epsilon^\beta} \, ||u||^2_{L^2(\R^n)}, 
\quad u \in C^\infty_0(\R^n).
\end{equation}

\begin{theorem} \label{Theorem 3.4a}  Let $\mu \in M^+(\R^n)$, $n \ge 2$, 
and let  $0<\beta<+\infty$. 

{\rm (i)} There exists a constant $C>0$ such that 
  {\rm (\ref{E:3.24a})} holds for every $\epsilon>0$  
if and only if $\mu$ satisfies the Frostman condition: 
\begin{equation}\label{E:3.25a}
\mu \left(B_\delta(x_0)\right) \le c \, \delta^{n-\frac{2\beta}{1+\beta}},
\end{equation}
for every ball $B_\delta(x_0)$ in $\R^n$, 
where $c$ does not depend on $\delta>0$ and $x_0$. 

{\rm (ii)} There exists a constant $C>0$ such that 
   {\rm (\ref{E:3.24a})} holds for every $\epsilon\in (0,1)$ 
if and only if {\rm (\ref{E:3.25a})} is valid for every $\delta\in (0, 1)$ 
and $x_0\in\R^n$. 
\end{theorem}

Statement (i) of Theorem \ref{Theorem 3.4a} follows from \cite{M2}, Lemma 1.4.7,  
where it is  proved that {\rm (\ref{E:3.25a})} holds for 
all $\delta>0$    
if and only if the multiplicative inequality   
\begin{equation}\label{E:3.26a}
||u||_{L^2(\R^n, \, d \mu)} \le C \, ||\nabla u||^p_{L^2(\R^n, \, d x)} \, 
||u||^{1-p}_{L^2(\R^n, \, d x)}, \quad u \in C^\infty_0(\R^n),
\end{equation}
holds, where $p = \frac \beta{1+\beta} \in (0, \, 1)$. Clearly, 
this inequality is equivalent to   {\rm (\ref{E:3.24a})} 
provided it holds for every $\epsilon>0$. 

There is also an inhomogeneous analogue  (\cite{M2}, Corollary 1.4.7/1) 
which states that {\rm (\ref{E:3.25a})} 
holds for all $\delta \in(0, \, 1)$  if and only if 
\begin{equation}\label{E:3.27a}
||u||_{L^2(\R^n, \, d \mu)} \le C \, ||u||^p_{W^{1,2} (\R^n)} \, 
||u||^{1-p}_{L^2(\R^n, d x)}, \quad u \in C^\infty_0(\R^n).
\end{equation}
 It is easy to see using a localization argument, as in the proof 
of Theorem \ref{Theorem 3.1a}, that this inequality holds if and only if 
{\rm (\ref{E:3.24a})} is valid for every $\epsilon\in (0, \, 1)$, which 
yields statement (ii) of Theorem \ref{Theorem 3.4a}.  
 
\section{The infinitesimal form boundedness  criterion}\label{Section 5}

We are now in a position to prove Theorems I and II stated in 
Sec. \ref{Section 2}.

\noindent{\it Proof of Theorem I.} 
Suppose  $Q \in \mathcal \D'(\R^n)$, $n \ge 2$, is represented in the form
\begin{equation}\label{E:2.0}
Q = \text{div} \,  \vec \Gamma + \gamma, \qquad \vec \Gamma \in {\mathbf
L}^2_{\rm loc}(\R^n), \quad \gamma \in L^1_{\rm loc}(\R^n),
\end{equation}  
where $\vec \Gamma$ and $\gamma$ respectively satisfy the conditions: 
\begin{align}\label{E:2.7}
& \lim_{\delta \to +0} \sup_{x_0 \in \R^n} \, \sup \, 
\left \{
\left\Vert \,  |\vec \Gamma|  \, u \right\Vert_{L^2(\R^n)} : \,
 u \in C^\infty_0(B_\delta(x_0)), \, \,   ||\nabla u||_{L^2(\R^n)} \le 1 
\right \} =0; \notag \\
& \lim_{\delta \to +0} \sup_{x_0 \in \R^n} \, \sup \, 
\left \{
\left\Vert \,  \gamma  \, |u|^2 \right\Vert_{L^1(\R^n)} : \,
 u \in C^\infty_0(B_\delta(x_0)), \, \,   ||\nabla u||_{L^2(\R^n)} \le 1 
\right \} =0.
\end{align} 
Then, for every $\epsilon>0$, there exists $\delta >0$ 
so that 
\begin{equation}\label{E:2.8}
\left\Vert \,  ( \, |\vec \Gamma|^2 + |\gamma|\,) \, |u|^2 \right\Vert_{L^1(\R^n)} \le
\epsilon \,  ||\nabla u||^2_{L^2(\R^n)}, 
\qquad \forall u \in C^\infty_0 (B_\delta(x_0)).
\end{equation} 
Hence by the Cauchy-Schwarz inequality, 
\begin{align}
\left \vert \langle Q u, \, u\rangle \right \vert   & \le  2 \, 
\left \vert \langle \vec \Gamma, \, u \,  \nabla u\rangle \right \vert  
+ \left \vert \langle \gamma u, \, u\rangle \right \vert  \notag \\ & \le 
2 \, \left \Vert \, |\vec \Gamma| \, u \right \Vert_{L^2(\R^n)} \, ||\nabla u||_{L^2(\R^n)}
+   \left \Vert \, \gamma \, |u|^2\right \Vert_{L^1(\R^n)}  \notag \\
& \le  (2 \sqrt{\epsilon} +  \epsilon)  \, ||\nabla u||^2_{L^2(\R^n)}.\notag 
\end{align}
Now we take the supremum on the left-hand side 
over all $u \in C^\infty_0 (B_\delta(x_0))$ such that 
$||\nabla u||^2_{L^2(\R^n)}\le 1$. Letting $\delta \to +0$,
 and then   $\epsilon \to +0$, we obtain:
$$
\lim_{\delta \to +0} \, \sup_{x_0 \in \R^n} \, \sup \, \{ \, 
\left | \langle Q u, \, u \rangle \right | \, : \,  u
\in C^\infty_0(B_\delta(x_0)), \, \,   ||\nabla u||_{L^2(\R^n)} \le 1 
\} = 0.
$$
By 
Lemma \ref{Lemma 2.1} this is equivalent to (\ref{E:2.1}). 

Conversely, suppose that,  for every $\epsilon>0$, 
there exists $C(\epsilon)>0$ such that 
  (\ref{E:2.1}) holds. By 
 the polarization identity, (\ref{E:2.1})   is
equivalent to  
\begin{align}\label{E:2.9}
 \left \vert \langle Q, \, u \, v\rangle \right \vert & \leq \, 
\epsilon^2  \, \left ( ||\nabla u||^2_{L^2(\R^n)} \notag  + 
\frac {C(\epsilon)}{\epsilon} ||u||^2_{L^2(\R^n)}\right)^{\frac 1 2} \\ & \times 
\left ( ||\nabla v||^2_{L^2(\R^n)} + 
\frac {C(\epsilon)}{\epsilon} ||v||^2_{L^2(\R^n)}\right)^{\frac 1 2}, 
\quad \forall u, \, v
\in  C^\infty_0(\R^n). 
\end{align}

We will show that (\ref{E:2.9}) implies  
 that, 
for every $\epsilon>0$, there exists $\delta >0$ such that 
\begin{equation}\label{E:2.10}
\int_e \, ( \, |\vec \Gamma (x)|^2 + |\gamma(x)|) \, dx \le \epsilon \,
\text{cap} \, (e), \end{equation}
for every compact set $e$, $\text{diam} \, (e) < \delta$, 
where 
\begin{equation}\label{E:vec}
\vec \Gamma = -\nabla (1-\Delta)^{-1} \, Q, 
\qquad \gamma = (1-\Delta)^{-1} \,
Q,
\end{equation}
so that $Q = \text{div} \, \vec \Gamma + \gamma$. 
 Here  $\text{cap}(\cdot)$ is the Bessel capacity associated with
the Sobolev space $W^{1,2}(\R^n)$. 
By Theorem~\ref{Theorem 3.1a},  inequality (\ref{E:2.10}) yields 
  (\ref{E:2.7}).

To deduce (\ref{E:2.10}) from  (\ref{E:2.9}), we fix an arbitrary 
compact set $e$ of positive  capacity
 such that $\text{diam} \, e \le \delta$. Without loss of generality we may assume that 
 $e \subset B_{\delta/2} (x_0)$, $x_0 \in \R^n$, where  $0<\delta\le 1$. 
The rest of the proof makes use of the factorization method developed  in
 \cite{MV2} combined with some new estimates for Bessel 
potentials of equilibrium measures. 

 Denote by
$\mu$ the equilibrium  measure associated with $e$ (see, e.g., \cite{AH},
\cite{M2}). By $P(x) $  denote the equilibrium 
potential of $\mu$ defined by:  
\begin{equation}\label{E:pot}
P(x) = (1-\Delta)^{-1} \mu (x) =  G_2\star \mu (x).
\end{equation}
Here the  Bessel kernel $G_\alpha(\cdot)$ of order $\alpha>0$
is defined by (\ref{E:3.15a}). 

We will need the following well-known properties of $\mu$ and $P$:
\begin{align}
& {\rm(a)} \quad  {\rm supp} \, (\mu) \subseteq e;\notag\\
& {\rm(b)} \quad  \mu(e) = {\rm cap} \, (e); \notag\\
& {\rm(c)} \quad  ||(1-\Delta)^{-\frac 1 2} P||^2_{L^2(\R^n)} = {\rm cap} \, (e);\notag\\ 
& {\rm(d)} \quad   P(x) = 1 \quad d \mu-{\rm a.e.}\notag\\
& {\rm(e)} \quad   P(x) \le 1 \quad {\rm on} \quad \R^n\notag
\end{align}

We will also need the asymptotics (see \cite{AH}, Sec. 1.2.5):
\begin{align}\label{E:2.11a}
& G_2(x) \asymp |x|^{2-n} \quad {\rm if} \, \, n\ge 3, \quad 
G_2(x) \asymp \log \,  \frac 1 {|x|} \quad {\rm if} \, \, n=2,  
 \quad |x| \to 0;\\   
& G_2(x) \asymp |x|^{\frac{1-n} 2} \, e^{-|x|}, \quad |x| \to +\infty, \quad n \ge 2.
\label{E:2.11b}
\end{align} 
Similar asymptotics hold for the derivatives of $G_2$ which will be used below as well. 

Sometimes, it will be more convenient to use a modified kernel, 
\begin{equation} \label{E:mod}
\widetilde G_2 (x) = \max \,  \left ( G_2(x) , \, 1 \right),
\end{equation} 
which does not have the exponential decay at infinity. 
Obviously, both $G_2$ and $\widetilde G_2$ are positive nonincreasing radial kernels. 
Moreover, $\widetilde G_2$ has the doubling property: 
$$ 
 \widetilde G_2 (2 x) \le \widetilde G_2(x) \le C(n) \, \widetilde G_2 (2 x).
$$
The corresponding modified potential is defined by:
\begin{equation} 
\widetilde P (x) = \widetilde G_2 \star \mu (x).
\end{equation}

Next, we fix a constant $\tau$, which  ultimately will be picked in the range 
\begin{equation} 
1 < 2 \tau < \min \, \left ( \frac n {n-2}, \, 2 
\right), \quad n \ge 2.
\label{E:R1}
\end{equation} 
The proof of statement (ii) of Theorem~I is based on  a series of propositions 
establishing  some  estimates 
for  the powers of the equilibrium potential, $P(x)^{2 \tau}$.

\begin{proposition}\label{Proposition 2.1}   Let $n\ge 2$, 
and let $0< 2\tau < \frac {n}{n-2}$.  
Then $\widetilde P^{2\tau}$ lies in the Muckenhoupt class $A_1$ on 
$\R^n$, i.e., 
\begin{equation}\label{E:A1}
M \, \widetilde P^{2\tau} (x) \le \, C(\tau, n) \, \widetilde P^{2\tau} (x), \quad 
dx-{\rm a.e.}
\end{equation} 
where $M$ denotes the Hardy-Littlewood 
maximal operator on $\R^n$, 
and the corresponding $A_1$-bound $C(\tau, n)$  depends
only on  $n$ and $\tau$. 
\end{proposition} 

\begin{proof} Let $\rho: \, \R^+\to \R^+$ be a nonincreasing function which satisfies 
the doubling condition, $\rho(2 s) \le c \, \rho(s)$, $s>0$. It is easy to see that the 
radial weight $\rho(|x|) \in A_1$   if and only if 
\begin{equation}
\int_0^r s^{n-1} \rho (s)^{2 \tau} \,  d s \le C   \, r^n \, \rho(r), 
\qquad r>0.
\label{radial}
\end{equation}
Moreover, the $A_1$-bound of $\rho$ is bounded by a constant which 
depends only on  $C$ 
in the preceding estimate and the doubling constant $c$ (see \cite{St2}). 

It follows from (\ref{E:2.11a}) that ${\widetilde G}_2(s)\asymp |s|^{2-n}$ if $n\ge 3$, 
${\widetilde G}_2(s)\asymp \log \, \frac 2 {|s|}$ if 
$n=2$, for $0<s<1$, and ${\widetilde G}_2(s)\asymp 1$ for $s\ge 1$. 
Hence, $\rho (|s|)={\widetilde G}_2^\tau (s)$ is a radial  
nonincreasing  kernel with the doubling property.  Clearly, (\ref{radial}) holds 
if and only if  $0< 2\tau < \frac {n}{n-2}$, and 
 the $A_1$-bound of ${\widetilde G}_2^\tau$ depends only on $\tau$ and $n$.  

By Jensen's inequality, 
${\widetilde G}_2^{\tau_1} \in A_1$ implies ${\widetilde G}_2^{\tau_2} \in A_1$ if 
$\tau_1\ge\tau_2$.
Hence, without loss of generality we may assume  $1\le 2 \tau \le \frac {n}{n-2}$. 
Then by Minkowski's integral inequality and the $A_1$-estimate for ${\widetilde G}_2^{2\tau}$ 
established above, 
it follows:
$$
M (\widetilde P^{2\tau}) (x)
 \le \left ( (M \, {\widetilde G}_2^{2\tau})^{\frac 1 {2\tau}} \star \mu(x) \right)^{2\tau}
\le C (\tau, n) \, ({\widetilde G}_2\star \mu)^{2\tau}(x) = C (\tau, n) 
 \widetilde P^{2\tau} (x).
$$
\end{proof}

\begin{proposition}\label{Proposition 2.2} 
Let $n \ge 2$ and $1<2 \tau < \min \, \left (\frac n {n-2}, 2 \right)$. Let  
$P = (1-\Delta)^{-1} \mu$, where $\mu\in M^+(\R^n)$ is  a compactly supported Borel measure 
on $\R^n$.  Then 
\begin{equation}\label{E:2.12}
 ||\nabla P^\tau||^2_{L^2(\R^n)}= \frac { \tau^2} {2 \tau-1} \, \left ( 
\int_{\R^n} P^{2\tau-1} \, d \mu - \int_{\R^n} P^{2\tau} \, d x \right),
\end{equation} 
provided $P^\tau\in W^{1,2}(\R^n)$.
\end{proposition}

\noindent{\bf Remark 5.1.} Suppose $\mu$ is the equilibrium measure of a compact set 
$e\subset \R^n$, so that $\mu(e) = {\rm cap} \, (e)$, $P(x)=1$ $d \mu$-a.e., and 
$0\le P(x) \le 1$ on $\R^n$.   
Then it is easy to see that $P^\tau\in W^{1,2}(\R^n)$ if 
 $1<2 \tau < \min \, \left (\frac n {n-2}, 1\right)$, and 
(\ref{E:2.12}) yields the estimates:
\begin{align}\label{E:2.12a}
 ||\nabla P^\tau||^2_{L^2(\R^n)} & \le  \frac {\tau^2} {2\tau-1} \, 
{\rm cap} \, (e), \\
||P^\tau||^2_{L^2(\R^n)} & \le {\rm cap} \, (e).
\label{E:2.12b}
\end{align}

\begin{proof} The proof of Proposition \ref{Proposition 2.2} is based on 
a multiple  integration by parts argument, 
and is analogous to  the proof of 
Proposition 2.5  in \cite{MV2}. Clearly, 
$$ 
\Delta P = P - (1-\Delta)P = P - \mu 
$$ 
 in the sense of distributions. Furthermore,  for $\tau>0$, 
\begin{align} 
\Delta P^\tau & = \tau (\tau -1) P^{\tau-2} |\nabla P|^2 + \tau \Delta P \, P^{\tau-1} \notag
\\ & = \tau (\tau -1) P^{\tau-2} |\nabla P|^2 + 
\tau P^\tau - \tau P^{\tau-1} \, \mu.
\notag 
\end{align} 
Using integration by parts and the preceding equation, we obtain:
\begin{align} 
 ||\nabla P^\tau||_{L^2(\R^n)}^2 & = \int_{\R^n} \nabla P^\tau \cdot \nabla P^\tau \, 
dx \notag\\
& = - \int_{\R^n} P^\tau \Delta P^\tau \, dx \notag \\ & = 
- \tau(\tau-1) \int_{\R^n} P^{2\tau-2} |\nabla P|^2 \, dx 
- \tau \int_{\R^n} P^{2 \tau} \, dx + \tau \int_{\R^n} P^{2 \tau-1} \, d\mu.
\notag 
\end{align}
Here integration by parts is easily justified, as long as 
$1<2 \tau < \min \, \left (\frac n {n-2}, 1\right)$, by looking at the behavior of 
$P(x)$ and $\nabla P(x)$ at infinity:
 \begin{align}  P(x) & \asymp |x|^{\frac{1-n}2} \, e^{-|x|} \, \mu(e), & \qquad |x|\to +\infty
\notag \\
|\nabla P(x)| & \le c(n) \, |x|^{\frac{1-n}2} \, e^{-|x|} \, \mu(e), & \qquad |x|\to +\infty. 
\notag 
\end{align}
(See details in the proof of Lemma 4.3 in \cite{MV2}.) Next, 
we use integration by parts again to obtain: 
\begin{align} 
\int_{\R^n} P^{2\tau-2} |\nabla P|^2 \, dx &= \int_{\R^n} P^{2\tau-2} \, \nabla P \cdot 
\nabla P \, dx \notag\\& = - \int_{\R^n} \Delta P \, P^{2\tau-1}  \, dx - 
(2 \tau -2) \int_{\R^n} P^{2\tau-2} |\nabla P|^2 \, dx.\notag
\end{align}
From the preceding equation we deduce:
\begin{align} 
(2\tau -1) \int_{\R^n} P^{2\tau-2} |\nabla P|^2 \, dx & = - \int_{\R^n} \Delta P \, P^{2\tau-1} 
 \, dx \notag\\ & =  \int_{\R^n} (1-\Delta) P \, P^{2\tau-1} \, dx -
\int_{\R^n} P^{2\tau} \, dx \notag\\& = \int_{\R^n}  P^{2\tau-1} \, d \mu -
\int_{\R^n} P^{2\tau} \, dx.\notag
\end{align}
Combining the preceding equations, we finally have:
\begin{equation} 
||\nabla P^\tau||^2_{L^2(\R^n)} = \frac {\tau^2} {2\tau-1}   
\left (\int_{\R^n}  P^{2\tau-1} \, d \mu - 
 \int_{\R^n} P^{2\tau} \, dx \right). \notag
\end{equation}
This completes the proof of Proposition~\ref{Proposition 2.2}. \end{proof}

We will need a more precise estimate of $||P^\tau||^2_{L^2(\R^n)}$ than 
 (\ref{E:2.12b}) in  Remark 5.1. Note that estimate (\ref{E:2.12a}) 
is sharp.

\begin{proposition}\label{Proposition 2.3} Let $n\ge 2$, and  $1 \le 2 \tau < \, \frac n
{n-2}$.  Let $\mu\in M^+(\R^n)$ be  a measure supported on 
 a compact set $e\subset B_\delta(x_0)$,
 where $0<\delta\le 1$.  Let $P = G_2\star \mu$. 
Then the following estimates hold:
\begin{align}
 & \int_{\R^n} P^{2\tau} \, dx   \le c(n, \tau) \, \mu (e)^{2 \tau},
 & \qquad n \ge 2;\label{E:2.13}\\
& \int_{B_{2 \delta}(x_0)} P^{2\tau} \, dx   \le c(n, \tau) \, \delta^{n
-(n-2) 2 \tau} \, \mu (e)^{2 \tau}, & \qquad n \ge 3;\label{E:2.13a}\\ 
& \int_{B_{2 \delta}(x_0)} P^{2\tau} \, dx   \le c(2, \tau) \, \delta^2 \, 
\log^{2\tau} \left(\frac 2 \delta\right) \, 
\, \mu (e)^{2 \tau}, & \qquad n =2. \label{E:2.13b}
\end{align}
\end{proposition}

\noindent{\bf Remark 5.2.} Estimates (\ref{E:2.13})$-$(\ref{E:2.13b}) are sharp,  
and can be reversed. 

\begin{proof} Suppose $e \subset B_{\delta} (x_0)$, where $0<\delta \le 1$.
Then 
$$
\int_{\R^n} P^{2\tau} \, dx  =  \int_{B_{2 \delta}(x_0)}
 P(x)^{2 \tau} \, dx +  \int_{B_{2 \delta}(x_0)^c} P (x)^{2 \tau} \, dx =
I+II.
$$ 
Using  Minkowski's integral inequality and the estimate $G_2(x) \le c(n) \, 
|x|^{2-n}$ for $|x|\le 3$ and $n \ge 3$, we get:
\begin{align}
I &=  c(n) \int_{B_{2 \delta}(x_0)} \left ( \int_e G_2(x-t) \, d \mu (t)
\right)^{2 \tau} dx \notag \\
& \le c(n, \tau) \, \int_{B_{2 \delta}(x_0)} \left ( \int_e \frac {d \mu 
(t)} {|x-t|^{n-2}}\right)^{2 \tau} \, dx \notag\\ & \le c(n, \tau) \, \left \{
\int_e  \left ( \int_{B_{2 \delta}(x_0)} \frac {dx} {|x-t|^{(n-2) 2 \tau}}
\right)^{\frac 1 {2 \tau}} d \mu (t)\right\}^{2 \tau}.\notag
\end{align}
Since $2 \tau < \, \frac n
{n-2}$, it follows: 
$$
\int_{B_{2 \delta}(x_0)} \frac {dx} {|x-t|^{(n-2) 2 \tau}} \le 
\int_{|x-t|< 3 \delta} \frac {dx} {|x-t|^{(n-2) 2 \tau}} \le c(n, \tau) \, 
\delta^{n -(n-2) 2 \tau}.
$$
Hence, 
 \begin{equation}\label{E:2.14} 
I \le c(n, \tau) \, \delta^{n
-(n-2) 2 \tau} \, \mu(e)^{2 \tau}.
\end{equation} 
This proves (\ref{E:2.13a}). In the case $n=2$, similar estimates using Minkowski's 
integral inequality and 
 $G_2(x) \asymp \log \frac C {|x|}$ as $|x| \to +0$, give  (\ref{E:2.13b}). 

To  estimate  $II$,  notice that  
$$ 
\int_{\R^n}  G_2 (x)^{2 \tau} \, dx < +\infty, \qquad 0<2 \tau < \frac n {n-2}.
$$
If $|x-x_0|\ge 2 \delta$ and $|x_0-t|< \delta$, it follows that 
$|x-t|\ge \frac 1 2 |x-x_0|$, and  $G_2(x-t)\le G_2 \left (\frac 1 2 |x-x_0|\right )$.  
Hence,  
\begin{align}
II & = \int_{B_{2 \delta}(x_0)} \left ( \int_e G_2(x-t) \, d \mu(t)
\right)^{2 \tau} \, dx \notag \\ & \le  \mu(e)^{2 \tau} \, \int_{\R^n}
    G_2 \left (\frac 1 2 |x-x_0|\right )^{2 \tau} \, dx 
\le  C(\tau, n) \, \mu(e)^{2 \tau}.\notag 
\end{align}
Combining the estimates for $I$ and $II$, we complete the proof of 
(\ref{E:2.13}). \end{proof}

The proof of the next  proposition is contained in the proof of Theorem 4.2 
 in \cite{MV2}.

\begin{proposition} \label{Proposition 2.4}  {\rm (\cite{MV2}).} 
 Let $1 < 2 \tau < \frac n {n-2}$. Then 
\begin{equation}\label{E:2.15} 
||\nabla (w \, P^{-\tau})||_{L^2(\R^n)} \le c(n, \tau) \, 
\left(  \int_{\R^n} ( \, |\nabla
w|^2 + |w|^2) \,  P^{-2\tau} \, dx \right)^{\frac 1 2}. 
\end{equation} 
\end{proposition}

\begin{proposition}\label{Proposition 2.5} Let 
$e\subset B_\delta(x_0), \, 0<\delta\le 1$, and $P(x) = (1-\Delta)^{-1}\mu$, 
where $\mu$ is the equilibrium measure associated with  $e$. 
Let $1 < 2 \tau < \min \, \left ( \frac n {n-2}, \, 2\right)$. Let 
 $w = \nabla (1-\Delta)^{-1} \psi$, where $ \psi \in C^\infty_0  (B_{\delta}
(x_0))$, $0<\delta\le 1$.   Then   the following estimates hold:
\begin{align}
& \int_{\R^n} \, \frac {|\nabla w|^2}  {P^{2\tau}} \, dx  \le C(n, \tau) \, 
\int_{B_{\delta} (x_0)} \frac {|\psi|^2} {P^{2\tau}} \, dx \qquad & n \ge 2;\label{E:2.16}\\ 
& \int_{\R^n} \, \frac {|w|^2}  {P^{2\tau}} \, dx  \le C(n, \tau) \, \delta^{n-(n-2) 2
\tau} \int_{B_{\delta} (x_0)} \frac {|\psi|^2} {P^{2\tau}} \, dx, 
 \qquad & n \ge 3;\label{E:2.17a}\\
& \int_{\R^2} \, \frac {|w|^2}  {P^{2\tau}} \, dx  \le C(2, \tau) \, \delta^{2} 
 \, \log^{2 \tau}
 \left (\frac 2 \delta\right)  \int_{B_{\delta} (x_0)}
 \frac {|\psi|^2} {P^{2\tau}} \, dx, \qquad & n=2.\label{E:2.17b}
\end{align}
\end{proposition}

\begin{proof} Estimate (\ref{E:2.16}) is actually contained in Lemma~4.3 \cite{MV2}  
 where a detailed proof was given only in the case $n\ge 3$.   
 For the sake of completeness, and since a similar argument is 
needed anyway in the proof of estimates 
(\ref{E:2.17a})$-$(\ref{E:2.17b}),   
we give a proof of 
(\ref{E:2.16}) here which is valid for $n\ge 2$.  

We split the integral on the left-hand side of (\ref{E:2.16}) into two parts:
$$
\int_{\R^n} \, \frac {|\nabla w|^2}  {P^{2\tau}} \, dx   = 
\int_{B_{2}(x_0)} \, \frac {|\nabla w|^2}  {P^{2\tau}} \, dx + 
\int_{B_{2}(x_0)^c} \, \frac {|\nabla w|^2}  {P^{2\tau}} \, dx = I + II, 
$$
where $w= \nabla (1-\Delta)^{-1} \psi$, and $\psi \in C^\infty_0 (B_{\delta}(x_0))$,  
$0<\delta\le 1$. 

To estimate $I$,   we begin by replacing $P^{2\tau}$ above with $\widetilde P^{2\tau}$, where 
$\widetilde P = \widetilde G_2\star \mu$ is  the modified equilibrium 
potential, and 
$$
 P(x) \asymp \widetilde P(x) \qquad  {\rm for}  \quad x\in B_{2}(x_0). 
$$
 Notice that $\nabla w = \nabla \nabla (1-\Delta)^{-1} \psi$, where $\nabla \nabla$ 
denotes the Hessian. Clearly, 
$$
\nabla \nabla (1-\Delta)^{-1} = \left \{ i R_j \,  i R_k \, \Delta (1-\Delta)^{-1} \right\},
\qquad j, k =1, \ldots, n,
$$
where $\{R_j\}_{j=1}^n$ are the Riesz transforms \cite{St2}, which are known to be 
bounded 
operators on  
$L^2(\R^n, \, \rho \, dx)$ with weights $\rho$ in the Muckenhoupt class $A_2$ defined by:   
$$
\sup \, \left \{ \, m_{B_\delta(x)} (\rho) \,\, 
m_{B_\delta(x)} (\rho^{-1}) : \quad
\delta>0, \, x\in \R^n \right\} < +\infty.
$$
Moreover, the operator norms of $R_j$ in  $L^2(\R^n, \, \rho \, dx)$ are bounded by 
a constant which 
 depends only on the $A_2$-bound of $\rho$ given
by the preceding formula. 

By  Proposition~\ref{Proposition 2.1}, the weight $\rho=\widetilde P^{2 \tau}$ 
 is in the Muckenhoupt class $A_1$, and hence to $A_2$ (see \cite{St2}), 
and its $A_2$-bound 
depends only on $\tau$, $n$.  Hence, $\rho^{-1}=\widetilde P^{-2 \tau} \in A_2$, and 
$$
||R_j f ||_{L^2 ( \R^n, \, \widetilde P^{-2 \tau} \, dx)} \le C(\tau, n) \, 
||f ||_{L^2 (\R^n, \, \widetilde P^{-2 \tau} \, dx)}, \qquad 
\forall f \in L^2 (\R^n, \, \widetilde P^{-2 \tau} \, dx),
$$
for every  $j=1, \ldots, n$. It follows:
\begin{equation}\label{I_1}
I \le c(\tau, n) \int_{B_{2}(x_0)} \frac {|\nabla w|^2}
 {\widetilde P^{2\tau}} \, dx \le C(\tau, n) \, \int_{\R^n}  \frac  
{| \Delta (1-\Delta)^{-1} \, \psi|^2}{\widetilde P^{2\tau}} \, dx. 
\end{equation}

It remains to notice that $\Delta (1-\Delta)^{-1} = 1- (1-\Delta)^{-1}$ is a Fourier 
multiplier operator on the space $L^2 (\R^n, \, \rho \, dx)$ for any weight 
$\rho\in A_2$, and its 
norm depends only on 
the Muckenhoupt bound of $\rho$. (See, e.g., \cite{MV2} or \cite{MV4}.) Thus, for 
$\rho=\widetilde P^{-2 \tau}$,
we have:
\begin{equation}\label{I_2}
I \le C(n, \tau)  \, \int_{\R^n}  \frac {|\psi|^2}{\widetilde P^{2\tau}} \, dx \le 
 C(\tau, n) \, \int_{B_\delta(x_0)}  \frac {|\psi|^2}{P^{2\tau}} \, dx,
\end{equation}
since $\widetilde P^{2\tau}(x) \asymp P^{2\tau}(x)$ on ${\rm supp} \, \psi \subset 
B_\delta(x_0)$, where $0<\delta\le 1$.

The term $II$ can be estimated directly using known inequalities for derivatives of Bessel 
potentials (see, e.g.,  \cite{AH}, Sec. 1.2.5):
$$
|\nabla^k G_2(x)| \le C(n) \, |x|^{\frac{1-n}{2}} \, e^{-|x|}, \qquad |x|\ge 1,
$$
where $k=0, 1, 2$. Notice that both $\psi$ and $\mu$ are 
supported in $B_\delta(x_0)\subset B_1(x_0)$. Hence, by the preceding estimate, 
$$
|\nabla \nabla (1-\Delta)^{-1} \psi (x)| \le C(n) \,  |x-x_0|^{\frac{1-n}{2}} \, e^{-|x-x_0|} \, 
||\psi||_{L^1(B_\delta(x_0)}, 
\qquad \forall x \in B_{2}(x_0)^c.
$$
From the lower estimate  $G_2(x) \ge c(n) \,  |x|^{\frac{1-n}{2}} \, e^{-|x|}$ for $|x|\ge 1$, 
we have: 
$$
P (x) = \int_{e} G_2(x-y) \, d \mu(y) \ge c(n) \,   |x-x_0|^{\frac{1-n}{2}} \, e^{-|x-x_0|} 
 \, \mu(e), \qquad \forall x \in  {B_{2}(x_0)^c}.
$$
Since $\mu(e) = {\rm cap} \, (e)>0$, it follows:
\begin{align} 
& II =\int_{B_{2}(x_0)^c} \frac{|\nabla \nabla (1-\Delta)^{-1} \psi (x)|^2} 
{P^{2 \tau}(x)} \, dx  \le C(\tau, n) \, ||\psi||_{L^1(B_\delta(x_0)}^2 {\rm cap} 
\, (e)^{-2\tau} \notag\\ & \times 
\int_{B_2(x_0)^c}  |x-x_0|^{(1-n)(1-\tau)} \, e^{-(2-2\tau)|x-x_0|} \, dx \le 
 C(\tau, n) \, ||\psi||_{L^1(B_\delta(x_0))}^2 \, {\rm cap} \, (e)^{-2\tau}, 
\notag
\end{align}
 where the constant $C(\tau, n)$ is finite   as long as  $\tau<1$. 

Now by Schwarz's inequality,
$$
||\psi||_{L^1(B_\delta(x_0))}^2 \le \int_{B_{\delta}(x_0)} 
\frac {|\psi(x)|^2}{P^{2\tau}(x)} \, dx \,  \int_{B_{\delta}(x_0)} P^{2\tau}(x) \, dx.
$$
By Proposition~\ref{Proposition 2.3}, $\int_{\R^n} P^{2\tau} \, dx \le 
c(\tau, n) \, {\rm cap} \, (e)^{2\tau}.$ Combining these estimates, 
we obtain: 
$$
II \le C(\tau, n) \, \int_{B_{\delta}(x_0)} 
\frac {|\psi(x)|^2}{P^{2\tau}(x)} \, dx.
$$ 
This completes the proof of estimate (\ref{E:2.16}).

We now prove estimates (\ref{E:2.17a})$-$(\ref{E:2.17b}). Splitting the integral 
$\int_{\R^n} \, \frac {|w|^2}  {P^{2\tau}} \, dx$ 
into two parts depending on $\delta$, 
we have:
$$
\int_{\R^n} \, \frac {|w|^2}  {P^{2\tau}} \, dx   = 
\int_{B_{2 \delta}(x_0)} \, \frac {|w|^2}  {P^{2\tau}} \, dx + 
\int_{B_{2 \delta}(x_0)^c} \, \frac {|w|^2}  {P^{2\tau}} \, dx = I + II, 
$$
where as before $w= \nabla (1-\Delta)^{-1} \psi$, and $\psi \in C^\infty_0 (B_{\delta}(x_0))$,  
$0<\delta\le 1$.  Using the estimate (\cite{AH}, Sec. 1.2.5):
$$
|\nabla G_2(x)| \le C(n) \, |x|^{1-n}, \qquad |x| \le 3, \qquad n\ge 2,
$$
we have:
$$
|w(x)| = | \nabla (1-\Delta)^{-1} \psi(x)|  \le \int_{B_\delta(x_0)} \frac {|\psi(t)|} 
{|x-t|^{n-1}} \, dt, \qquad \forall x \in B_{2 \delta}(x_0).
$$
We next use a version of Hedberg's inequality (see \cite{AH}): 
$$
\int_{B_\delta(x_0)} \frac {|\psi(t)|} 
{|x-t|^{n-1}} \, dt \le  \sum_{j=0}^{+\infty} \delta 2^{-j} 
\int_{ \delta 2^{-j-1}<|x-t|\le \delta 2^{-j}}
\frac {|\psi(t)|} 
{|x-t|^{n}} \, dt \le C(n) \, \delta \, M \psi (x), 
$$
for every $x \in B_{2 \delta}(x_0)$. Here $M$ is the Hardy-Littlewood maximal function. 

Since $P^{2\tau}\asymp  
\widetilde P^{2\tau}$ on $B_{2 \delta} (x_0)$ we have:
$$
\int_{B_{2 \delta}(x_0)} \, \frac {|w|^2}  {P^{2\tau}} \, dx \le 
C(\tau, n) \, \delta^2 \int_{\R^n} \frac {|M \psi|^2}  {\widetilde P^{2\tau}} \, dx
\le C(\tau, n) \, \delta^2 \int_{B_{\delta}(x_0)} \frac {|\psi|^2}  {P^{2\tau}} \, dx.
$$
Here we have used the fact that $M$ is a bounded operator on $L^2(\R^n, \, \rho)$ 
with the weight $\rho = \widetilde P^{-2\tau}\in A_2$, and its operator norm 
is bounded by a constant which depends only 
on the $A_2$-bound of $P^{-2\tau}$ (see \cite{St2}). Hence, by 
Proposition~\ref{Proposition 2.1} it is  bounded by $C(\tau, n)$. Note that we can 
use $P^{2\tau}$ in place of 
$\widetilde P^{2\tau}$ on ${\rm supp} \,  \psi \subset B_\delta(x_0)$ since 
$0<\delta\le 1$. 

Thus, 
$$
I   \le C(n, \tau) \, \delta^{2} 
\int_{B_{\delta} (x_0)} \frac {|\psi|^2} {P^{2\tau}} \, dx, 
 \qquad  n \ge 2. 
$$

To estimate $II$,  notice that  
$$
|\nabla (1-\Delta)^{-1} \psi (x)| \le \int_{B_\delta(x_0)} |\nabla G_2(x-t)| \, |\psi(t)| 
\, dt, 
\qquad \forall x \in B_{2 \delta}(x_0)^c,
$$
where $|x-t|\ge |x-x_0|-\delta\ge \frac 1 2 |x-x_0|$. 
Using  estimates  for  derivatives of Bessel kernels (\cite{AH}, Sec. 1.2.5), we have:
\begin{align}
|\nabla G_2(x-t)| & \le C(n) \, G_1 (x-x_0), \qquad |x-x_0| \le 1; \notag\\
|\nabla G_2(x-t)| & \le C(n) \, G_2 (x-x_0), \qquad |x-x_0| \ge 1,
\notag
\end{align}
where $x\in B_{2 \delta}(x_0)^c$ and $t\in B_{\delta}(x_0)$, $0<\delta\le 1$. 

Since $\mu(e) = {\rm cap} \, (e)>0$, and $|x-t|\le |x-x_0| + \delta$, we estimate:
$$
P (x) = \int_e G_2(x-t) \, d \mu(t) \ge c(n) \,  G_2(x-x_0) \, {\rm cap} \, (e), 
\qquad \forall x \in  {B_{2 \delta}(x_0)^c}.
$$
It follows:
\begin{align} 
II & =\int_{B_{2 \delta}(x_0)^c} \frac{|\nabla (1-\Delta)^{-1} \psi (x)|^2} 
{P(x)^{2 \tau}} \, dx  \le C(\tau, n) \, ||\psi||_{L^1(B_\delta(x_0))}^2 {\rm cap} 
\, (e)^{-2\tau} \notag\\ & \times 
\left ( \int_{B_1(x_0)} \frac {G_1(x-x_0)^2} {G_2(x-x_0)^{2 \tau}} \, dx + 
\int_{B_1(x_0)^c} G_2(x-x_0)^{2- 2 \tau} \, dx \right). 
\notag
\end{align}
It is easy to see that, for $2\tau\in(1, \, 2)$, 
\begin{equation}
\int_{B_1(x_0)} \frac {G_1(x-x_0)^2} {G_2(x-x_0)^{2 \tau}} \, dx + 
\int_{B_1(x_0)^c} G_2(x-x_0)^{2- 2 \tau} \, dx < +\infty. 
\notag
\end{equation}
Indeed, estimates  for  Bessel kernels yield: 
\begin{align}
\frac {G_1(x)^2} {G_2(x)^{2 \tau}} & \le C(\tau,n)  \, |x|^{2\tau(n-2)-2(n-1)}, 
 \qquad |x|\le 1, 
 \qquad n \ge 3;\notag\\
\frac {G_1(x)^2} {G_2(x)^{2 \tau}} & \le C(\tau,n)  \,
 |x|^{-2} \, \log^{-2\tau} \,  \frac 2 {|x|}, 
\qquad |x| \le 1, \qquad n =2;\notag \\ 
G_2(x)^{2- 2 \tau} & \le C(\tau,n)  \, |x|^{(1-n)(1-\tau)} 
e^{(2\tau-2)|x|} , \qquad |x|\ge 1,  \qquad n \ge 2,
\notag
\end{align}
where the exponents satisfy 
the  inequalities $2\tau(n-2)-2(n-1)>-n$ in the case $n\ge 3$ and $-2\tau< -1$ 
in the case $n=2$ 
for $|x|\le 1$; 
and $2\tau-2<0$ for $|x|\ge 1$. Thus,  
$$
II \le C(\tau, n) \, ||\psi||_{L^1(B_\delta(x_0))}^2 \, {\rm cap} \, (e)^{-2\tau}.
$$
Now by Schwarz's inequality, 
$$
||\psi||_{L^1(B_\delta(x_0))}^2 \le \int_{B_{\delta}(x_0)} 
\frac {|\psi|^2}{P^{2\tau}} \, dx \,  \int_{B_{\delta}(x_0)} P^{2\tau} \, dx,
$$
and by Proposition~\ref{Proposition 2.3}, 
\begin{align}
& \int_{B_{\delta}(x_0)} P^{2\tau} \, dx \le C(\tau, n) \,  
\delta^{n-(n-2) 2\tau} \, {\rm cap} \, (e)^{2\tau}, 
 \qquad & n \ge 3;\notag\\
& \int_{B_{\delta}(x_0)} P^{2\tau} \, dx \le C(\tau, 2) \,  \delta^2 \,  \log^{2 \tau}
 \left (\frac 2 \delta\right) \, {\rm cap} \, (e)^{2\tau}, \qquad & n=2.\notag
\end{align}
Combining these estimates, 
we obtain: 
\begin{align}
& II \le C(\tau, n) \, \delta^{n-(n-2) 2
\tau} \int_{B_{\delta} (x_0)} \frac {|\psi|^2} {P^{2\tau}} \, dx, 
 \qquad & n \ge 3;\notag\\
& II \le C(\tau, 2) \, \delta^{2} \, \log^{2 \tau}
 \left (\frac 2 \delta\right)  \int_{B_{\delta} (x_0)}
 \frac {|\psi|^2} {P^{2\tau}} \, dx, \qquad & n=2.\notag
\end{align}

Since $2 \tau>1$, obviously $\delta^2 \le \delta^{n-(n-2) 2\tau}$ if $n\ge 3$, and  
 $\delta^2 \le \delta^2 \, \log^{2 \tau}
 \left (\frac 2 \delta\right)$ in the case $n=2$. Thus, the above inequalities
 for  $I$ and $II$
yield the desired  estimate for $\int_{\R^n} \frac {|w|^2} {P^{2\tau}} \, dx$.
\end{proof}
 
The proof of the next proposition is very 
similar to that of Proposition~\ref{Proposition 2.5}, with obvious changes in the 
orders of Bessel potentials, and so we omit it.

\begin{proposition}\label{Proposition 2.6} 
Let $1 < 2 \tau < \min \, \left ( \frac n {n-2}, \, 2\right)$. Let 
 $u = (1-\Delta)^{-1} \psi$, where $ \psi \in C^\infty_0  (B_{\delta}
(x_0))$, $0<\delta\le 1$.   Then   
\begin{align} 
 \int_{\R^n}  \, \frac {|u|^2}  {P^{2\tau}} \, dx  & \le C(\tau, n) \, 
\delta^{n-(n-2) 2\tau} \int_{B_{\delta} (x_0)} \frac {|\psi|^2} {P^{2\tau}} \, dx, 
\qquad n\ge 3;
\label{E:2.19a}\\ 
 \int_{\R^2}  \, \frac {|u|^2}  {P^{2\tau}} \, dx  & \le C(\tau, 2) \,  \delta^2 \, 
\log^{2\tau} \left(\frac 2 \delta\right) \int_{B_{\delta} 
(x_0)} \frac {|\psi|^2} {P^{2\tau}} \, dx, \qquad n=2;
\label{E:2.19b}
\end{align}
\end{proposition}

We can now  complete the proof of the implication (\ref{E:2.9})$\Rightarrow$(\ref{E:2.10}).
Let $e\subset B_\delta(x_0)$ 
be a set of positive capacity,  and  let $\vec \psi = \{\psi_k\}_{k=1}^n$, where 
$\psi_k \in C^\infty_0  (B_{\delta}
(x_0))$, $k=1, \ldots, n$, and $0<\delta\le 1$. Denote by $P = {G_2}\star\mu$ the 
potential of the equilibrium measure 
$\mu$ associated with $e$ which satisfies properties (a)$-$(e) 
listed above. Pick $\tau$ so that (\ref{E:R1}) holds. 

We next apply (\ref{E:2.9}) with $u= P^\tau$ and $v = \frac w {P^\tau}$, 
where 
$w = {\rm div} \,  (1-\Delta)^{-1} \vec \psi$. (Since  $u, v \in W^{1,2}(\R^n)$, we 
can actually 
use that inequality for  mollifiers of $u$ and $v$ and then pass to the limit as 
in \cite{MV2}.) We obtain:
\begin{align}
\left \vert \langle Q, \,  w \rangle \right \vert 
& \leq \epsilon \left ( ||\nabla P^\tau||^2_{L^2(\R^n)} + \frac
{C(\epsilon)} {\epsilon} \, ||P^\tau||^2_{L^2(\R^n)}\right)^{ \frac 1 2} 
\notag \\ &\times \left ( \left \Vert \nabla \left (\frac {w} 
{P^\tau} \right)\right \Vert^2_{L^2(\R^n)} + \frac
{C(\epsilon)} {\epsilon}   \left \Vert  \frac {w} 
{P^\tau}\right \Vert^2_{L^2(\R^n)}\right)^{ \frac 1 2}.
\label{E:P1}
\end{align} 

Let us assume  that $n\ge 3$. From the preceding inequality where 
$w = {\rm div} \,  (1-\Delta)^{-1} \vec \psi$, we deduce
using Propositions \ref{Proposition 2.2}$-$\ref{Proposition 2.6}:
 \begin{align}\label{E:2.20}
\left \vert \langle Q, \,   {\rm div} \, (1-\Delta)^{-1} \vec \psi \rangle \right \vert 
& \leq \, c(n, \tau) \,  \epsilon \left ( \text{cap} \, (e) + 
\text{cap} \, (e)^{2 \tau} \, \frac
{C(\epsilon)} {\epsilon} \right)^{ \frac 1 2} \notag \\ &\times \left ( 1 +  
\delta^{n - (n-2) 2 \tau} \frac
{C(\epsilon)} {\epsilon} \right)^{ \frac 1 2} \, 
\left(  \int \frac{|\vec \psi|^2}
 {P^{2\tau}} \, 
{dx} 
\right)^{\frac 1 2},
\end{align} 
 for all $\vec \psi
\in  C^\infty_0(B_\delta(x_0))^n$. 

Letting $\vec \Gamma = - \nabla (1- \Delta)^{-1} Q$, we have:
$$
 \langle Q, \,   {\rm div} \, (1-\Delta)^{-1} \vec \psi \rangle  
=  \langle \vec \Gamma, \,  \vec \psi \rangle, 
$$
where $\vec \Gamma \in L^2_{\rm loc}(\R^n)^n$. (This formal equation was 
justified  in the sense of distributions  in \cite{MV2} 
under a weaker assumption that $Q$ is form 
bounded with respect to $-\Delta$.) 

Now taking the 
supremum over all $\vec \psi \in C^\infty_0(B_\delta(x_0))^n$,
we deduce, as in  \cite{MV2}:
\begin{align}
 \int_{B_\delta(x_0)} |\vec \Gamma|^2 \,  P^{2\tau}  \, dx \notag & \le c(n, \tau) \, 
\epsilon^2  \left ( \text{cap} \, (e) + \text{cap} \, (e)^{2 \tau} \, \frac
{C(\epsilon)} {\epsilon} \right) \notag \\ & \times \left ( 1 +  
\delta^{n - (n-2) 2 \tau} \frac
{C(\epsilon)} {\epsilon} \right) .
\notag
\end{align}
 Since $P (x) \ge 1$ quasieverywhere, and hence 
$dx$-a.e. on $e$ \cite{AH}, it follows that, 
for every compact set  $e\subset B_\delta(x_0)$ of positive capacity,
$$
\int_e \, |\vec \Gamma|^2 \, dx \le \, c(n, \tau) \, 
\epsilon^2  \left ( \text{cap} \, (e) + \text{cap} \, (e)^{2 \tau} \, \frac
{C(\epsilon)} {\epsilon} \right) \left ( 1 +  
\delta^{n - (n-2) 2 \tau} \frac
{C(\epsilon)} {\epsilon} \right).
$$
Thus,
$$
\frac{\int_e \, |\vec \Gamma|^2 \, dx }{\text{cap} \, (e)} \le 
c(n, \tau) \, 
\epsilon^2 \left (1 + \text{cap} \, (e)^{2 \tau-1} \, \frac
{C(\epsilon)} {\epsilon} \right) \left ( 1 +  
\delta^{n - (n-2) 2 \tau} \frac
{C(\epsilon)} {\epsilon} \right).
$$
Since $\text{cap} \, (e)^{2 \tau-1} \le c(n, \tau)\,  \delta^{(n-2) (2
\tau-1)}$, it follows that 
$$
\frac{\int_e \, |\vec \Gamma|^2 \, dx }{\text{cap} \, (e)} \le 
c(n, \tau) \, 
\epsilon^2 \, \left (1 + \delta^{(n-2)(2 \tau-1)} \, \frac
{C(\epsilon)} {\epsilon} \right) \left ( 1 +  
\delta^{n - (n-2) 2 \tau} \frac
{C(\epsilon)} {\epsilon} \right).
$$
Taking the supremum over all $e$ and letting $\delta \to +0$, we obtain:
$$
\lim_{\delta \to +0} \sup_{e: \, \text{diam} \, e \le \delta} 
\frac{\int_e \, |\vec \Gamma|^2 \, dx }{\text{cap} \, (e)} \le \epsilon.
$$

The case $n=2$ requires  usual 
modifications where a logarithmic term appears, as in Proposition~\ref{Proposition 2.3}.

A similar argument using Proposition~\ref{Proposition 2.6} 
shows that the same
inequality  holds with $|\gamma|$ in place of $|\vec \Gamma|^2$. 
This proves (\ref{E:2.10}), which  
yields  (\ref{E:2.7}). The proof of Theorem~I is complete.
\qed

Theorem~II stated in Sec.~\ref{Section 2} is an immediate consequence of  
Theorem~I  combined with Theorem~\ref{Theorem 3.1a} and Corollary~\ref{Corollary 3.3a}.  


\section{Form subordination criteria}\label{Section 6}

Suppose $Q \in \mathcal \D'(\R^n)$, $n \ge
2$, and $\beta >0$. In this section we study quadratic form  inequalities 
of Trudinger type; i.e., we characterize $Q$ for which  there exist 
positive constants $c$ and $\epsilon_0$ such that, 
for every $\epsilon\in (0, \, \epsilon_0)$, 
\begin{equation}\label{E:3.1}
\left \vert \langle Q u, \, u\rangle \right \vert  \leq \epsilon \, ||\nabla
u||^2_{L^2(\R^n)} + c \, \epsilon^{-\beta} \, ||u||^2_{L^2(\R^n)}, 
\qquad  \forall u \in
C^\infty_0(\R^n).  
\end{equation}

By Corollary~\ref{Corollary 2.3}, this is equivalent to the existence of 
$C>0$ and  $\delta_0>0$ such that 
\begin{equation}\label{E:3.2}
\left\vert \langle Q u, \, u\rangle \right\vert  \leq \, 
 C \, \delta^{\frac 2 {1 + \beta}} \,  ||\nabla u||^2_{L^2}, 
\qquad  \forall u  \in C^\infty_0( B_\delta(x_0)), 
\end{equation}
for all $\delta \in (0, \, \delta_0)$ and $x_0 \in \R^n$.
Moreover, if $\epsilon_0 =1$ then one may set $\delta_0=1$;  
 if $\epsilon_0 = + \infty$ then  $\delta_0 = + \infty$, 
 and the converse is also true.

We are now in a position to prove the following theorem 
which combines Theorems IV and V stated in Sec.~\ref{Section 2}.

\begin{theorem} \label{Theorem 3.2}
Let $Q \in \mathcal \D'(\R^n)$, $n \ge 2$, and $\beta > 1$. 
Then the  following statements hold. 

{\rm (i)} Suppose there exist $\delta_0>0$, and 
$\vec \Gamma \in 
L^2_{\rm loc}(\R^n)^n$, 
$\gamma \in L^2_{\rm loc}$ such that 
\begin{equation}\label{E:3.3}
Q = \text{{\rm div}} \,\, \vec \Gamma + \gamma,
\end{equation}
where $\vec \Gamma$ and $\gamma$ satisfy respectively the 
inequalities: 
\begin{align}
& \int_{B_\delta(x_0)} \left\vert \vec \Gamma (x)\right\vert^2 
\, dx \le 
C_1 \, \delta^{n-2 +\frac 4 {1 +\beta}}, \qquad 0<\delta <\delta_0,
 \label{E:3.4} \\
& \int_{B_\delta(x_0)} \left\vert \gamma (x) \right\vert 
\, dx \le 
C_2 \, \delta^{n-2 + \frac 2 {1 +\beta}}, \qquad 0<\delta <\delta_0, 
\label{E:3.5}
\end{align}
with the constants $C_1, \, C_2$ which do not depend on $x_0$ and $\delta$.
Then there exist $c>0$ and $\epsilon_0>0$ such 
that {\rm (\ref{E:3.1})} holds  for every $\epsilon \in (0, \,
\epsilon_0)$.  

{\rm (ii)} Conversely, 
suppose {\rm (\ref{E:3.1})} holds  for every $\epsilon \in (0, \,
\epsilon_0)$.
Then $Q$ can be represented in the form {\rm (\ref{E:3.3})} 
 where $\vec \Gamma  \in \mathbf{L}^2_{\rm loc} (\R^n)$,  
$\gamma \in  L^2_{\rm loc} (\R^n)$ satisfy respectively 
conditions {\rm (\ref{E:3.4})} and {\rm (\ref{E:3.5})}.

{\rm (iii)}  If $\epsilon_0=1$, then one can set 
\begin{equation}\label{E:3.6}
\vec \Gamma = 
- \nabla (1 - \Delta)^{-1} \, Q, \qquad \gamma = (1 - \Delta)^{-1} \, Q,
\end{equation}
in {\rm (\ref{E:3.3})}, where 
$\delta_0=1$ 
in {\rm (\ref{E:3.4})}, {\rm (\ref{E:3.5})}.

If $\epsilon_0= +\infty$, then one can set  
\begin{equation}\label{E:3.7}
\vec \Gamma = 
\nabla \Delta^{-1} \, Q, \qquad \gamma = 0, 
\end{equation}
in {\rm (\ref{E:3.3})}, where 
$\delta_0=+\infty$ in 
{\rm (\ref{E:3.4})}, {\rm (\ref{E:3.5})}.
\end{theorem}

\begin{proof} We first prove statement (i). Suppose $\beta>1$ and 
$u \in C^\infty_0 (B_\delta(x_0))$. Suppose that $Q$ 
is represented in the form (\ref{E:3.3}), and estimates 
(\ref{E:3.4}), (\ref{E:3.5}) hold for every $\delta>0$, 
i.e., $\delta_0=+\infty$. Applying the multiplicative 
inequality for nonnegative measures 
(\cite{M2}, Theorem 1.4.7) to 
$|\vec \Gamma|^2 \, dx$ 
and  $|\gamma| \, dx$ respectively, we get:
$$
 \int_{\R^n} |\vec  \Gamma(x)|^2 \, |u(x)|^2 \, dx  \le C_1 \,  
||\nabla u||^{\frac{2 (\beta-1)}{1 + \beta}}_{L^2} \, 
||u||^{\frac{4}{1 + \beta}}_{L^2},
$$
and
$$
\int_{\R^n} |\gamma (x)| \, |u(x)|^2 \, dx  \le C_2 \, 
||\nabla u||^{\frac{2 \beta}{1 + \beta}}_{L^2} \, 
||u||^{\frac{2}{1 + \beta}}_{L^2}.
$$
Hence, 
\begin{align}
\left\vert \langle Q u, \, u\rangle \right\vert & \leq \, 
\left\vert \langle \vec \Gamma u, \, \nabla u\rangle 
\right\vert + \left\vert \langle \gamma u, \,  u\rangle 
\right\vert    \notag \\ & \le 
||\vec  \Gamma \, u||_{L^2} \, 
||\nabla u||_{L^2}
+ || \gamma \, |u|^2||_{L^1} 
 \notag \\ 
& \le C_1^{\frac 1 2} \, 
||\nabla u||^{1 + \frac{\beta-1}{1 + \beta}}_{L^2} \,
 ||u||^{\frac 2 {1+\beta}}_{L^2} + C_2 \, 
||\nabla u||^{\frac{2 \beta}{1 + \beta}}_{L^2} \, 
||u||^{\frac{2}{1 + \beta}}_{L^2}.
\notag
\end{align}
Combining the preceding estimates with  the inequality (\ref{E:2.3}),
$$
||u||_{L^2} \le c(n) \,  \delta \, 
||\nabla u||_{L^2}, \qquad u \in C^\infty_0(B_\delta(x_0)),
$$
we obtain:
$$
\left\vert \langle Q u, \, u\rangle \right\vert \le 
C \, \delta^{\frac {2}{1 + \beta}} \,  
||\nabla u||_{L^2}, \qquad u \in C^\infty_0(B_\delta(x_0)).
$$
 By Lemma \ref{Lemma 2.1} this implies (\ref{E:3.1}). 

If $\delta_0=1$, then  in the proof above one has to replace 
 $||\nabla u||_{L^2}$ with the inhomogeneous Sobolev norm
$||u||_{W^{1,2}}$, and apply Corollary 1, Sec. 1.4.7 \cite{M2} 
together with  (\ref{E:2.3}). For an 
arbitrary $\delta_0$, one can use $||\nabla u||_{L^2} + \delta_0^{-1} \, 
||u||_{L^2}$ in place of $||u||_{W^{1,2}}$ and  the corresponding analogue 
of the multiplicative inequality for positive measures 
in \cite{M2}. This completes 
the proof of statement (i).

We now prove statement (ii). We first establish  a localized 
version of (ii) which holds for every $\beta>0$.
 Let 
$\eta_{\delta, x_0} (x) =\eta \left(\frac {x-x_0}{\delta}\right)$
 where $\eta\in C^\infty_0(B_1(0))$ is a 
smooth cut-off function such that $|\eta (x)| \le 1$ and 
$|\nabla \eta (x)| \le 1$ for $x \in B_1(0)$. It is worthwhile 
to observe 
that in all localized estimates below the constants will not depend on 
a particular choice of $\eta$.

\begin{proposition} \label{Proposition 3.3}
Suppose $Q \in \mathcal \D'(\R^n)$, $n \ge 2$, 
and $\beta>0$. Suppose that {\rm (\ref{E:3.2})} holds for every 
$\delta\in(0, \, \delta_0)$, where either $\delta_0=1$ or 
 $\delta_0=+\infty$. Then the 
following inequalities hold:
\begin{equation}\label{E:3.8i}
\vert \langle Q \, u, \, v \rangle \vert \le c \, 
\delta^{\frac {2}{1 + \beta}} \, ||\nabla u||_{L^2} \, ||\nabla v||_{L^2}, 
\qquad u, \, v \in C^\infty_0(B_\delta(x_0)), 
\end{equation}
and, in particular, 
\begin{equation}\label{E:3.8}
\vert \langle Q \, \eta_{\delta, x_0}, \, v
\rangle \vert \le c \, 
\delta^{\frac {2}{1 + \beta} + \frac {n-2} 2} \,||\nabla v||_{L^2}, 
\qquad v \in C^\infty_0(B_\delta(x_0)). 
\end{equation}
Moreover, if $\delta>\delta_0 =1$, then 
\begin{equation}\label{E:3.9}
\vert \langle Q \, \eta_{\delta, x_0}, \, v
\rangle \vert \le c \, 
\delta^{2+ \frac {n-2} 2} \,||\nabla v||_{L^2}, \qquad 
v \in C^\infty_0(B_\delta(x_0)). 
\end{equation}

\end{proposition}

\begin{proof} Clearly, (\ref{E:3.8i}) follows from  (\ref{E:3.2}) by the polarization 
identity. Letting $u = \eta_{\delta, x_0}$ in (\ref{E:3.8i}) one deduces (\ref{E:3.8}). 
Finally, in the case $\delta_0=1$ and $\delta>1$, (\ref{E:3.9}) follows by using 
a polarized form of (\ref{E:3.1}) where one sets  $\epsilon =1$, and $u = \eta_{\delta, x_0}, \, 
v \in  C^\infty_0(B_\delta(x_0))$, together with estimate (\ref{E:2.3}): 
\begin{align}
\vert \langle Q \, \eta_{\delta, x_0}, \, v
\rangle \vert & \le 
\sqrt{||\nabla \eta_{\delta, x_0}||^2_{L^2} + c \, ||\eta_{\delta, x_0}||^2_{L^2}} \, \, 
\sqrt{ ||\nabla v||^2_{L^2} + c \, ||v||^2_{L^2}} \notag \\ 
& \le C \, \sqrt{ \delta^{n-2} + c \, \delta^n} \,  \, 
\sqrt{ ||\nabla v||^2_{L^2} + c \, \delta^2 \, ||\nabla v||^2_{L^2}} \notag \\
& \le C \, \delta^{2+ \frac {n-2} 2} \,||\nabla v||_{L^2}.\notag
\end{align}
\end{proof} 

We now establish certain localized versions 
of the necessary condition for {\rm (\ref{E:3.1})} and {\rm (\ref{E:3.2})}.

\begin{proposition} \label{Proposition 3.4}
Suppose $Q \in \mathcal \D'(\R^n)$, $n \ge
2$, and $\beta>0$. Suppose that {\rm (\ref{E:3.2})} holds for every 
$\delta\in(0, \, \delta_0)$, where either $\delta_0 =1$ or 
$\delta_0=+\infty$.  Let $\vec \Gamma$ and $\gamma$
 be defined by {\rm (\ref{E:3.6})} if $\delta_0 =1$, or by 
{\rm (\ref{E:3.7})} if $\delta_0 =+\infty$.

{\rm(i)} For  $n\ge 3$, 
$$
\int_{\R^n} | \nabla (1-\Delta)^{-1} (\eta_{\delta, x_0} \, Q)|^2 \, dx 
\le c \, \delta^{\frac 4 {1+\beta} + n-2}, \qquad 
0<\delta\le 1,
$$
and 
$$
\int_{\R^n} |(1-\Delta)^{-1} (\eta_{\delta, x_0} \, Q)| \, dx 
\le c \, \delta^{\frac 2 {1+\beta} + n-2}, \qquad 
0<\delta\le 1.
$$

{\rm(ii)} For  $n\ge 3$ and $\delta_0=+\infty$, 
$$
\int_{\R^n} | \nabla \Delta^{-1} (\eta_{\delta, x_0} \, Q)|^2 \, dx 
\le c \, \delta^{\frac 4 {1+\beta} + n-2}, \qquad 
0<\delta<+\infty.
$$

{\rm(iii)} For  $n\ge 2$, 
$$
\int_{B_\delta(x_0)} | \nabla (1-\Delta)^{-1} (\eta_{\delta, x_0} \, Q)|^2 \, dx 
\le c \, \delta^{\frac 4 {1+\beta} + n-2}, \qquad 
0<\delta\le 1,
$$
and 
$$
\int_{\R^n} |(1-\Delta)^{-1} (\eta_{\delta, x_0} \, Q)| \, dx 
\le c \, \delta^{\frac 2 {1+\beta} + n-2}, \qquad 
0<\delta\le 1.
$$

 {\rm(iv)} For  $n\ge 2$ and $\delta_0=+\infty$, 
$$\int_{B_\delta(x_0)} | \nabla \Delta^{-1} 
(\eta_{\delta, x_0} \, Q)|^2 \, dx 
\le c \, \delta^{\frac 4 {1+\beta} + n-2}, \qquad 
0<\delta<+\infty.
$$
\end{proposition}

\begin{proof} We first consider the case $\delta_0=1$. The
homogeneous case $\delta_0 =+\infty$ will be treated in a similar way.

We  pick another cut-off function $\zeta \in C^\infty_0(B_1(0))$ 
such that $\zeta (x) \, \eta (x) = \eta (x)$, and 
let $\zeta_{\delta, \, x_0} (x) = \zeta  \left(\frac {x-x_0}{\delta}\right)$. 
Let $\vec \psi \in C^\infty_0(\R^n)^n$. Setting 
$v =\zeta_{\delta, \, x_0} \,  \text{div} \, (1 - \Delta)^{-1} \vec \psi$ 
in (\ref{E:3.8}) we obtain:
\begin{align}\label{E:3.10}
 & \vert \langle Q \, \eta_{\delta, x_0}, \, 
v  \rangle \vert  \le c \, 
\delta^{\frac {2}{1 + \beta}+ \frac {n-2} 2} \notag \\ & \times 
\left ( \left \Vert  (\nabla \zeta_{\delta, x_0}) 
 \text{div} \, (1 - \Delta)^{-1} \vec \psi\right \Vert_{L^2}  
  +  ||\nabla \text{div} \, (1 - \Delta)^{-1} \vec \psi||_{L^2} 
\right).
\end{align}

It follows, e.g. from the Plancherel theorem,
$$||\nabla \text{div} \, (1 - \Delta)^{-1} \vec \psi||_{L^2} \le c \, 
||\vec \psi||_{L^2}.
$$
Obviously,  
$$|\nabla  \zeta_{\delta, \, x_0} (x)| \le c \, \delta^{-1} \le  c \,
|x-x_0|^{-1}, \qquad x \in B_\delta(x_0).
$$ 
For $n \ge 3$, applying the preceding estimates together 
with  Hardy's inequality, we obtain:
$$\int_{B_\delta (x_0)}  |\nabla \zeta_{\delta, x_0}|^2 \,  
 | \text{div} \, (1 - \Delta)^{-1} \vec \psi|^2 \, dx  \le c \, 
\int |\nabla \text{div} \, (1 - \Delta)^{-1} \vec \psi|^2 \, dx  \le c \, 
||\vec \psi||^2_{L^2}.
$$
Thus, in case $n \ge 3$,
\begin{equation}\label{E:3.11}
\vert \langle Q \, \eta_{\delta, x_0}, \, 
\text{div} \, (1 - \Delta)^{-1} \vec \psi  \rangle \vert \le c \, 
\delta^{\frac {2}{1 + \beta} + \frac{n-2} 2} \, ||\vec \psi||_{L^2}, 
\end{equation}
or equivalently,
\begin{equation}\label{E:3.12}
\left \vert \int_{\R^n} \nabla (1 - \Delta)^{-1}  
( \eta_{\delta, x_0} \, Q) \cdot \vec \psi (x) \, dx \right \vert \le c \, 
\delta^{\frac {2}{1 + \beta} + \frac {n-2}2} \, ||\vec \psi||_{L^2}. 
\end{equation}
Minimizing over all $\vec \psi$, we  obtain:
\begin{equation}\label{E:3.13}
\left \Vert \nabla (1 - \Delta)^{-1}  ( \eta_{\delta, x_0} \, Q)
\right \Vert_{L^2(\R^n)}^2 \le c \, \delta^{\frac {4}{1 + \beta} + n-2}.
\end{equation}

\noindent{\bf Remark 6.1.}\label{Remark 6.1} In the case $n=2$, 
the preceding estimate fails 
 even for nonnegative $Q$.
It can be replaced with a weaker estimate:
\begin{equation}\label{E:3.13'}
\left \Vert \nabla (1 - \Delta)^{-1}  ( \eta_{\delta, x_0} \, Q)
\right \Vert_{L^2(\R^2)}^2 \le c \, \delta^{\frac {4}{1 + \beta}} 
\log \frac 2 \delta, \qquad 0<\delta\le 1.
\tag{6.13$'$}
\end{equation}

Nevertheless, the following localized version still holds for 
every $n \ge 2$. 
We set 
$v =\zeta_{\delta, \, x_0} \,  \text{div} \, (1 - \Delta)^{-1} \vec \psi$ 
in (\ref{E:3.8}) where in this case we assume that $ \vec \psi \in 
C^\infty_0(B_\delta(x_0))^n$.  Since 
$\text{supp} \, \vec \psi \subset B_\delta(x_0)$, 
we deduce from (\ref{E:3.10}) that, for $x \in B_\delta(x_0)$,
\begin{align}
& | \text{div} \, (1 - \Delta)^{-1} \vec \psi (x)|  \le 
c \, \int_{B_\delta(x_0)} |\nabla G_2(x-y)| \, |\vec \psi (y)| \, dy \notag
\\ & \le 
c \, \int_{|x-y|< 2 \delta} |x-y|^{1-n} 
\, |\vec \psi (y)| \, dy, \notag 
\end{align}
where $G_\nu$ is a Bessel kernel of order $\nu>0$. 
Here we have used the known 
estimate (see, e.g., \cite{AH}, Sec. 1.2.5):
$$| \nabla G_2(x)| \le c \, G_1(x) 
\le c \, |x|^{1-n}, \qquad |x|\le 1, \quad n \ge 2.
$$

We now apply the inequality (see \cite{AH}):
$$
\int_{|x-y|< 2 \delta} |x-y|^{1-n} 
\, |\vec \psi (y)| \, dy \le c \, \delta \,  M \, |\vec \psi (x)|,
$$
where $M$ is the 
 Hardy-Littlewood maximal function operator. 
Hence, by the maximal function inequality,
\begin{equation}\label{E:3.14}
\int_{B_\delta (x_0)}  |\nabla \zeta_{\delta, x_0} (x)|^2 \,  
 | \text{div} \, (1 - \Delta)^{-1} \vec \psi (x)|^2 \, dx  \le c \, 
||M |\vec \psi| \, ||^2_{L^2} \le c ||\vec \psi||^2_{L^2} .
\end{equation}
This gives:
\begin{equation}\label{E:3.15}
\left \vert \int_{B_\delta(x_0)} \nabla (1 - \Delta)^{-1}  
( \eta_{\delta, x_0} \, Q) \cdot \vec \psi \, (x) \, dx \right \vert \le c \, 
\delta^{\frac {2}{1 + \beta} + \frac{n-2}2} \, ||\vec \psi||_{L^2}. 
\end{equation}
Minimizing over all such $\vec \psi$, we arrive at the estimate:
\begin{equation}\label{E:3.16}
 \int_{B_\delta(x_0)} |\nabla (1 - \Delta)^{-1}  
( \eta_{\delta, x_0} \, Q) (x)|^2 \, dx  \le c \, 
\delta^{\frac {4}{1 + \beta} + n-2}. 
\end{equation}

We now estimate 
$\int |(1-\Delta)^{-1} ( \eta_{\delta, x_0} \, Q) (x)| \, dx.$ 
Let  $\phi$ be a compactly supported function, 
$\phi \in L^\infty(\R^n)$, 
and set 
$v = \zeta_{\delta, x_0} \, (1-\Delta)^{-1} \phi$ in (\ref{E:3.8}),
where $\zeta_{\delta, x_0}$ was defined above. Then
\begin{align}
| \langle \eta _{\delta, x_0} \, Q, & \,  (1-\Delta)^{-1} \phi \rangle | 
\le  c \, \delta^{\frac 2 {1+\beta} + \frac{n-2} 2}\notag  \\ 
& \times ( \delta^{-2} \, 
||(1-\Delta)^{-1} \phi||^2_{L^2(B_\delta(x_0))} + 
||\nabla (1-\Delta)^{-1} \phi||^2_{L^2(B_\delta(x_0))})^{\frac 1 2}.
\notag
\end{align}

Notice that $G_2 \in L^1(\R^n)$, $n \ge 2$, and hence   
$$ ||(1-\Delta)^{-1} \phi||_{L^\infty} \le ||\phi||_{L^\infty} \, 
||G_2||_{L^1}.
$$
This gives:
$$ \delta^{-2} \, ||(1-\Delta)^{-1} \phi||^2_{L^2(B_\delta(x_0))} \le c \, 
\delta^{n-2} \, ||\phi||^2_{L^\infty}.
$$
Similarly, $\nabla G_2 \in L^1(\R^n)$, $n \ge 2$, and so
$$ ||\nabla (1-\Delta)^{-1} \phi||_{L^\infty} \le ||\phi||_{L^\infty} \, 
||\nabla G_2||_{L^1},
$$
which implies:
$$ ||\nabla (1-\Delta)^{-1} \phi||^2_{L^2(B_\delta(x_0))} \le c \, 
\delta^n \, ||\phi||^2_{L^\infty}.
$$

Since $\delta\le 1$, we can combine the preceding estimates to get:
$$| \langle \eta _{\delta, x_0} \, Q, \,  (1-\Delta)^{-1} \phi \rangle | 
= | \langle (1-\Delta)^{-1} (\eta _{\delta, x_0} \, Q), \,  \phi \rangle | 
\le c \, \delta^{\frac 2 {1+\beta} + n-2} ||\phi||_{L^\infty}.
$$
This yields
\begin{equation}\label{E:3.17}
 \int_{\R^n} |(1 - \Delta)^{-1}  
( \eta_{\delta, x_0} \, Q) (x)| \, dx  \le c \, 
\delta^{\frac {2}{1 + \beta} + n-2}, \qquad 0<\delta\le 1.
\end{equation}

Analogous inequalities hold in the homogeneous case $\delta_0=+\infty$. 
If $n \ge 3$, we set  in (\ref{E:3.8})
$v= \zeta_{\delta, x_0} \, \text{div} \,  \Delta^{-1} \, \vec \psi$ in 
(\ref{E:3.8}), 
where $\vec \psi\in C^\infty_0(\R^n)^n$. Estimating exactly as above,
using Plancherel's theorem and Hardy's inequality, we get: 
\begin{equation}\label{E:3.18}
 \int_{\R^n} |\nabla \Delta^{-1}  
( \eta_{\delta, x_0}, \, Q) (x)|^2 \, dx  \le c \, 
\delta^{\frac {4}{1 + \beta} + n-2}, \qquad 0<\delta<+\infty. 
\end{equation}
For $n \ge 2$, we assume that $\psi\in C^\infty_0(B_\delta(x_0))^n$,
and notice that, for $x \in B_\delta(x_0)$, 
$$|\text{div} \,  \Delta^{-1} \, \vec \psi (x)| \le 
c \, \int_{|y-x|< 2 \delta} \frac{|\vec \psi (y)|}{|x-y|^{n-1}} dy \le 
c \, \delta \, M (|\vec \psi|)(x).
$$
Applying the maximal function inequality, we deduce 
as in the inhomogeneous 
case:
$$\int_{B_\delta(x_0)} |\nabla \zeta_{\delta, x_0} (x)|^2  \, 
|\text{div} \,  \Delta^{-1} \, \vec \psi (x)|^2 \, dx + 
\int_{B_\delta(x_0)} |\nabla \text{div} \,  \Delta^{-1}
 \, \vec \psi (x)|^2 \, dx \le c \, ||\vec \psi||^2_{L^2}.
$$
This yields:
\begin{equation}\label{E:3.19}
 \int_{B_\delta(x_0)} |\nabla \Delta^{-1}  
( \eta_{\delta, x_0}, \, Q) (x)|^2 \, dx  \le c \, 
\delta^{\frac {4}{1 + \beta} + n-2}, \qquad 0<\delta<+\infty. 
\end{equation}
\end{proof} 

We now state the key lemma which gives estimates of the global 
``antiderivative'' of $Q$.

\begin{lemma} \label{Lemma 3.5}
Suppose $Q \in \mathcal \D'(\R^n)$, $n \ge
2$, and $\beta>0$. 
Suppose that {\rm (\ref{E:3.2})} holds for every 
$\delta\in(0, \, \delta_0)$, where either $\delta_0=1$ or 
$\delta_0= +\infty$.  
Then the following statements hold:

{\rm(i)} For   $\delta_0=1$ and $\beta>1$, 
$$\int_{B_\delta(x_0)} | \nabla (1-\Delta)^{-1} \, Q|^2 \, dx 
\le c \, \delta^{\frac 4 {1+\beta} + n-2}, \quad 
0<\delta\le 1,
$$
and 
$$\int_{B_\delta(x_0)} |(1-\Delta)^{-1} \, Q| \, dx 
\le c \, \delta^{\frac 2 {1+\beta} + n-2}, \quad 
0<\delta\le 1.
$$

{\rm(ii)} For  $\delta_0=1$ and $0<\beta \le 1$, 
$$\int_{B_\delta(x_0)} |\nabla (1-\Delta)^{-1} \, Q - m_{B_\delta(x_0)} 
(\nabla (1-\Delta)^{-1} \, Q)|^2 \, dx 
\le c \, \delta^{\frac 4 {1+\beta} + n-2}, \quad 
0<\delta\le 1,
$$
and
$$\int_{B_\delta(x_0)} |(1-\Delta)^{-1} \, Q| \, dx 
\le c \, \delta^{\frac 2 {1+\beta} + n-2}, \quad 
0<\delta\le 1.
$$

{\rm(iii)} For  $\delta_0=+\infty$ and $\beta>1$, 
$$\int_{B_\delta(x_0)} |\nabla \Delta^{-1} \, Q|^2 \, dx 
\le c \, \delta^{\frac 4 {1+\beta} + n-2}, \quad 
0<\delta<+\infty.
$$

{\rm(iv)} For  $\delta_0=+\infty$ and $0<\beta \le 1$, 
$$\int_{B_\delta(x_0)} |\nabla \Delta^{-1} \, Q - m_{B_\delta(x_0)} 
(\nabla \Delta^{-1} \, Q)|^2 \, dx 
\le c \, \delta^{\frac 4 {1+\beta} + n-2}, \quad 
0<\delta<+\infty.
$$
\end{lemma}

\begin{proof} We fix $x_0\in\R^n$ and $\delta>0$, and define a special  
partition of unity $\{\psi_j\}_{j=0}^{+\infty}$
associated with $\delta$ and $x_0$. (For the sake of 
convenience,  we will suppress the dependence 
of $\psi_j$ on 
$\delta, x_0$ in our notation used below.)  
 
Let $\eta^0_{\delta, x_0} (x) = \eta
 \left (\frac {2 |x-x_0|} {\delta}\right)$, where 
$\eta\in C^\infty_0(\R_+)$ is such that 
$\eta(x) = 1$ if $0\le x\le \frac 1 2$ and $\eta(x)=0$ if $x\ge 1$. 
Next, we pick  
$\zeta \in C^\infty_0(\R_+)$ so that $\zeta(x) = 0$ if $0\le x \le \frac 1 4$ 
or $x\ge 2$, and 
$\zeta(x) = 1$ if $\frac 1 2 \le x \le 1$, and set 
$$\eta^j_{\delta, x_0} (x) = 
\zeta \left (\frac { |x-x_0|} {2^{j-2}\delta}\right), 
\qquad 
j = 1, 2, \ldots.
$$ 

Letting $\phi_{\delta, x_0}(x) = \sum_{j=0}^{+\infty} \, 
\eta^j_{\delta, x_0} (x),$ we see that 
clearly $\phi_{\delta, x_0} \in C^\infty(\R^n)$, and 
$1 \le \phi_{\delta, x_0}(x)\le 3$. We now set 
\begin{equation} \label{E:3.20}
\psi_j (x) = \frac {\eta^j_{\delta, x_0}(x)} {\phi_{\delta, x_0}(x)}, 
\qquad  j =0,1, \ldots.
\end{equation} 
We observe that \begin{equation} \label{E:3.21}
0 \le \psi_j (x) \le  1, \qquad 
|\nabla \psi_j (x)| \le  \frac {c }{ 2^j \delta}, \qquad 
j =0,1, \ldots,
\end{equation} 
and 
\begin{equation} \label{E:3.22}
\sum_{j=0}^{+\infty} \, \psi_j  (x) \equiv 1, \qquad 
\psi_j  \in C^\infty_0(B_{2^j \delta} (x_0) 
\setminus B_{2^{j-2} \delta} (x_0)), 
\qquad j = 1, 2,  \ldots.
\end{equation} 

We also denote by 
$\widetilde \psi_j$, $j= 1, 2,  \ldots$, a function in 
$C^\infty_0 (B_{2^j \delta}(x_0)\setminus B_{2^{j-2} \delta} (x_0))$
 with a slightly larger 
support so that 
$$\widetilde \psi_j (x) \, \psi_j (x) =\psi_j (x), \qquad 
|\nabla  \widetilde \psi_j (x) | \le \frac c {2^j \delta}.
$$

We now prove statement (i). Let $0<\delta\le 1$ and $\beta>1$. 
Using (\ref{E:3.22}),
 we obtain:
$$||\nabla (1-\Delta)^{-1}  Q||_{L^2(B_\delta(x_0))} \le  \, 
\left \Vert \sum_{j=0}^{+\infty} \, 
|\nabla (1-\Delta)^{-1} (\psi_j \, Q)| \, \right \Vert_{L^2(B_\delta(x_0))}.
$$

Suppose that $j=0, 1, 2$. Since 
$\psi_j \in C^\infty(B_{4 \delta}(x_0))$, 
by Proposition \ref{Proposition 3.4},
$$
||\nabla (1-\Delta)^{-1} (\psi_j \, Q)||_{L^2(B_\delta(x_0))}
\le c \, \delta^{\frac 2 {1+\beta} + \frac{n-2} 2}.
$$

For $j\ge 3$, notice that 
$$| \nabla (1-\Delta)^{-1} (\psi_j \, Q) (x)| = 
| \nabla (1-\Delta)^{-1} (\psi_j \widetilde \psi_j \, Q) (x)|
= |\langle  \psi_j \, Q, \, \nabla G_2 (x-\cdot) \, \widetilde 
\psi_j  \rangle|.
$$

Fix $x \in B_\delta(x_0)$. Applying  (\ref{E:3.8})
with $\psi_j$ in place of $\eta_{\delta, x_0}$, 
$\nabla G_2 (x-\cdot) \widetilde \psi_j$ in place of $v$, and $2^j \delta$ 
in place of $\delta$ respectively, we obtain:
\begin{align}
|\nabla (1-\Delta)^{-1} (\psi_j \, Q) (x)| & \le  
\left \vert  \langle \psi_j \, Q, \, 
\nabla (1-\Delta)^{-1} G_2 (x- \cdot) \, \widetilde \psi_j \rangle 
\right \vert \notag \\ & \le c \,
 (2^j \delta)^{\frac 2 {1+\beta} + \frac {n-2} 2} \, 
|| \nabla (\nabla G_2 (x-\cdot)  
\widetilde \psi_j)||_{L^2(B_{2^j\delta}(x_0))},
\notag
\end{align}
for $2^j \delta \le 1$. For $2^j \delta > 1$, a similar estimate 
follows from (\ref{E:3.9}):
$$|\nabla (1-\Delta)^{-1} (\psi_j \, Q) (x)| \le  
c \, (2^j \delta)^{2  + \frac {n-2} 2} \, 
|| \nabla (\nabla G_2 (x-\cdot)  
\widetilde \psi_j)||_{L^2(B_{2^j\delta}(x_0))}.
$$
 
Let $y \in B_{2^j\delta }(x_0)\setminus
B_{2^{j-2}\delta }(x_0)$. Then
$2^{j-3} \delta \le |y-x| \le 2^{j+1} \delta$, $j \ge 3$. 
Using  the preceding inequalities, together with (\ref{E:3.21}) 
and   
the  estimates for the Bessel kernel
(see, e.g.,  \cite{AH}, Sec 1.2.5): 
\begin{align}
|\nabla G_2 (x-y)| & \le \frac c {|x-y|^{n-1}}, \qquad 
|\nabla \nabla G_2 (x-y)| \le \frac c {|x-y|^{n}}, \qquad |x-y| \le 1, 
\notag \\ |\nabla G_2 (x-y)| & \le \, c \, e^{-|x-y|}, 
\qquad |\nabla \nabla G_2 (x-y)| \le \, c \, e^{-|x-y|}, \qquad |x-y| > 1,
\notag
\end{align}
we get:
\begin{align}
& |\nabla_y ( (\nabla G_2) (x-y) \, \widetilde \psi_j (y))| 
 \le |\nabla G_2 (x-y)| \, |\nabla \widetilde \psi_j (y)| \notag\\
& + |\nabla \nabla G_2 (x-y)| \, 
 \widetilde \psi_j (y)    \le \frac {c} {(2^j \delta)^n}, \quad \text{for} 
\quad 2^j \delta \le 1,\notag 
\end{align}
and 
$$|\nabla_y ( (\nabla G_2) (x-y) \, \widetilde \psi_j (y))| 
\le c \, e^{-2^j \delta}, \quad \text{for}\quad  2^j \delta >1.
$$
Hence, 
\begin{align}
& || \nabla (\nabla G_2 (x-\cdot) 
 \widetilde \psi_j)||_{L^2(B_{2^j\delta}(x_0))} \le 
c \, (2^j \delta)^{-\frac n 2} \quad & \text{if} \quad 2^j
 \delta \le 1,\notag\\
& || \nabla (\nabla G_2 (x-\cdot) 
 \widetilde \psi_j)||_{L^2(B_{2^j\delta}(x_0))} \le 
c \, (2^j \delta)^{\frac n 2} e^{-2^j \delta} \quad & \text{if} \quad
 2^j \delta > 1.\notag
\end{align}

Since $\beta>1$, it follows that $\frac 2 {1+\beta}-1<0$. 
Thus, for $x \in B_\delta(x_0)$, we have a uniform estimate:
\begin{align}
& \sum_{j=3}^{+\infty} \, 
|\nabla (1-\Delta)^{-1} (\psi_j \, Q) (x)|   \le 
c \sum_{j=3}^{[\log_2 \frac 2 \delta]} \, 
(2^j \delta)^{\frac 2 {1+\beta} + \frac {n-2} 2}(2^j \delta)^{-\frac n 2}
\notag\\ & +c 
 \sum_{j=[\log_2 \frac 2 \delta]}^{+\infty} \, (2^j \delta)^{n+1} \, 
e^{-2^j \delta}
\le c \, \delta^{\frac 2 {1+\beta} -1}.\notag
\end{align}
Hence, 
$$\left \Vert \sum_{j=3}^{+\infty} \, 
|\nabla (1-\Delta)^{-1} (\psi_j \, Q) (x)| \, 
\right \Vert_{L^2(B_\delta(x_0))} \le c \, \delta^{\frac 2 {1+\beta} 
+ \frac 
{n-2} 2}.
$$
From  the preceding estimates it follows:
$$
||\nabla (1-\Delta)^{-1}  Q||_{L^2(B_\delta(x_0))}
\le c \, \delta^{\frac 2 {1+\beta} + \frac{n-2} 2}.
$$

A similar argument yields  the estimate:
$$\int_{B_\delta(x_0)} |(1-\Delta)^{-1} \, Q| \, dx 
\le c \, \delta^{\frac 2 {1+\beta} + n-2}, \qquad 
0<\delta\le 1.
$$
We have:
$$ 
\int_{B_\delta(x_0)} |(1-\Delta)^{-1} \, Q| \, dx \le 
\sum_{j=0}^{+\infty} \int_{B_\delta(x_0)} 
|(1-\Delta)^{-1} \, (\psi_j Q)| \, dx.
$$
By Proposition \ref{Proposition 3.4}, 
$$
\int_{B_\delta(x_0)} 
|(1-\Delta)^{-1} \, (\psi_j Q)| \, dx \le c \, \delta^{\frac 2 {1+\beta} + 
n-2}, \qquad 0<\delta \le 1,
$$
for $j=0, 1, 2$. If $j\ge 3$ and $x \in B_\delta(x_0)$, we estimate 
using Proposition 3.3 as above: 
\begin{align}
|(1-\Delta)^{-1} \, (\psi_j Q) (x)| & \le | \langle \psi_j Q, \, 
G_2(x-\cdot) \, \widetilde \psi_j \rangle | \notag \\ & \le 
c \, (2^j \delta)^{\frac 2 {1+\beta} + \frac {n-2}2} \, 
||\nabla (G_2(x-\cdot) \, \widetilde \psi_j) ||_{L^2}, \notag
\end{align}
if $2^j \delta\le 1$. For $2^j \delta> 1$,
$$
|(1-\Delta)^{-1} \, (\psi_j Q) (x)| \le c \, (2^j \delta)^{2 
 + \frac {n-2}2} \, 
||\nabla (G_2(x-\cdot) \, \widetilde \psi_j) ||_{L^2}.
$$

Suppose $n \ge 3$.
For $y \in B_{2^j\delta }(x_0)\setminus
B_{2^{j-2}\delta }(x_0)$,  
we have  $2^{j-3} \delta \le |y-x| \le 2^{j+1} \delta$, $j \ge 3$. 
Hence, for $2^j \delta\le 1$ and $|x-y|\le 1$, 
$$|\nabla (G_2(x-y) \, \widetilde \psi_j (y)) | \le 
|\nabla G_2(x-y)| \, \widetilde \psi_j (y) + 
G_2(x-y) \, |\nabla \widetilde \psi_j (y)| \le 
c \, (2^j \delta)^{1-n}.
$$
 Similarly, for $2^j \delta> 1$ and $|x-y|> 1$,
$$|\nabla (G_2(x-y) \, \widetilde \psi_j (y)) | \le c \,  
e^{-2^j \delta}.
$$

Thus,
\begin{align} 
& ||\nabla (G_2(x-\cdot) \, \widetilde \psi_j) ||_{L^2} \le c \, 
(2^j \delta)^{1-\frac n 2}, \quad & \text{if} \quad 
 2^j \delta \le 1,\notag\\
 & ||\nabla (G_2(x-\cdot) \, \widetilde \psi_j) ||_{L^2} \le c \, 
(2^j \delta)^{\frac n 2} e^{-2^j \delta}, \quad 
& \text{if} \quad  2^j \delta >1.\notag
\end{align}

It follows:
$$\sum_{j=3}^{+\infty} \, |(1-\Delta)^{-1} \, (\psi_j Q) (x)| \le 
c \, \sum_{j=3}^{\log \frac 2 \delta} (2^j \delta)^{\frac 2 {1+\beta}} 
+ \sum_{j > \log \frac 2 \delta} (2^j \delta)^{n +1} 
e^{-2^j \delta} \le C,
$$
where  $C$ does not depend on $\delta$.

A similar  uniform estimate holds in the case $n=2$ as well, since  
$G_2(x-y)\le c \, \log \frac 2 {2^j \delta}$ if 
$2^j \delta \le |x-y| \le 2^{j+2} \delta$ and $2^j \delta\le 1$, and hence 
$$\sum_{j=3}^{+\infty} \, |(1-\Delta)^{-1} \, (\psi_j Q) (x)| \le 
c \, \sum_{j=3}^{\log \frac 2 \delta} 
\log \frac 2 {2^j \delta} \, (2^j \delta)^{\frac 2 {1+\beta}} 
+ \sum_{j > \log \frac 2 \delta} (2^j \delta)^{3} 
e^{-2^j \delta} \le C,
$$
where  $C$ does not depend on $\delta$.

Combining the preceding estimates we obtain:
$$\int_{B_\delta(x_0)} 
|(1-\Delta)^{-1} \,  Q| \, dx \le c \, \delta^{\frac 2 {1+\beta} +n-2},
 \qquad 0<\delta\le 1,
$$
which proves statement (i) of Lemma \ref{Lemma 3.5}.

To prove statement (ii), we  modify the above argument as follows. 
For $j=0,1,2$, it follows from Proposition 3.4 that 
$$\int_{B_\delta(x_0)} |\nabla (1-\Delta)^{-1} \,  (\psi_j Q)|^2 \, dx 
\le c \, \delta^{\frac 4 {1+\beta} + n-2}, \quad 0<\delta\le 1.
$$
Obviously, this gives as well the estimate:
$$\int_{B_\delta(x_0)} |\nabla (1-\Delta)^{-1} \,  (\psi_j Q) - 
m_{B_\delta(x_0)} (\nabla (1-\Delta)^{-1} \,  (\psi_j Q))|^2 \, dx 
\le c \, \delta^{\frac 4 {1+\beta} + n-2}, 
$$
where $0<\delta\le 1$.

For $j=3, 4, \ldots$ and $x, x' \in B_\delta(x_0)$, we deduce 
the following uniform estimate:
\begin{align}
 & |\nabla (1-\Delta)^{-1} \,  (\psi_j Q) (x) - 
\nabla (1-\Delta)^{-1} \,  (\psi_j Q) (x')| \notag \\ & \le 
|\langle \psi_j Q, \, (\nabla G_2(x-\cdot) - \nabla G_2(x'-\cdot))
 \vec \psi_j \rangle| \notag \\
& \le c \, (2^j \delta)^{\frac 2 {1+\beta} +\frac {n-2} 2} 
||\nabla ( (\nabla G_2(x-\cdot) - \nabla G_2(x'-\cdot))\vec \psi_j)||_{L^2},
\notag
\end{align} 
if $2^j \delta \le 1$; for $2^j \delta > 1$, there is a similar estimate 
with $\beta=0$ on the right-hand side. 

The rest of the proof 
is similar to the case $\beta>1$. 
For $x, x' \in B_\delta(x_0)$ and $y \in B_{2^j\delta }(x_0)\setminus
B_{2^{j-2}\delta }(x_0)$, in the case $2^j \delta \le 1$, 
we use the estimates
\begin{align} 
& |\nabla G_2(x-y) - \nabla G_2(x'-y)|  \le c \,  \delta \, 
(2^j \delta)^{-n},
\notag \\
& |\nabla \nabla G_2(x-y) - \nabla \nabla G_2(x'-y)| \le c \, \delta \, 
(2^j \delta)^{-n-1},\notag
\end{align} 
together with similar estimates:
\begin{align} 
& |\nabla G_2(x-y) - \nabla G_2(x'-y)|  \le c \,  \delta \, 
e^{-2^j \delta},
\notag \\
& |\nabla \nabla G_2(x-y) - \nabla \nabla G_2(x'-y)| \le c \, \delta \, 
e^{-2^j \delta},\notag
\end{align} 
in the 
case $2^j \delta >1$.
From this it follows exactly as in the proof of statement (i):
\begin{align} 
& |\nabla (1-\Delta)^{-1} \,  (\psi_j Q) (x) - 
\nabla (1-\Delta)^{-1} \,  (\psi_j Q) (x')| \le c \, \delta 
(2^j \delta)^{-2 + \frac 2 {1+\beta}}, & \quad 2^j \delta \le 1,\notag\\
& |\nabla (1-\Delta)^{-1} \,  (\psi_j Q) (x) - 
\nabla (1-\Delta)^{-1} \,  (\psi_j Q) (x')| \le c \, \delta 
e^{-2^j \delta}, & \quad 2^j \delta > 1.\notag
\end{align} 

Hence, for every $\beta>0$, and $x, \, x' \in B_\delta(x_0)$: 
$$\sum_{j\ge 3} \,  |\nabla (1-\Delta)^{-1} \,  (\psi_j Q) (x) - 
\nabla (1-\Delta)^{-1} \,  (\psi_j Q) (x')| \le c \, 
\delta^{- 1 + \frac {2}{1 + \beta}}.
$$
Consequently,
$$\int_{B_\delta(x_0)} | \nabla (1-\Delta)^{-1} \,  Q - 
m_{B_\delta(x_0)} (\nabla (1-\Delta)^{-1} \,  Q)|^2 \, dx \le 
c \, \delta^{n- 2 + \frac {4}{1 + \beta}}, \quad 0<\delta\le 1.
$$

Notice that a lower term estimate 
$$\int_{B_\delta(x_0)} |(1-\Delta)^{-1} \,  Q| \, dx \le 
c \, \delta^{n- 2 + \frac {2}{1 + \beta}}, \quad 0<\delta\le 1,
$$
was proved above for every $\beta>0$. 
This completes the proof of statement (ii)
 of Lemma \ref{Lemma 3.5}.
 
The proofs of statements (iii) and (iv) in the homogeneous case  are 
similar but simpler, and 
follow 
by using analogous estimates for the derivatives of 
Riesz kernels in place of Bessel kernels. 
\end{proof}

Statement (ii) 
of Theorem \ref{Theorem 3.2} now follows from Propositions 
\ref{Proposition 3.3} 
 and \ref{Proposition 3.4}. The proof of 
Theorem \ref{Theorem 3.2} is complete. \end{proof}

\noindent{\bf Remark 6.2.}\label{Remark 6.2} Estimates similar to those 
used in the 
proof of Lemma \ref{Lemma 3.5} in the inhomogeneous 
case $\delta_0=1$,  give additionally that ({\ref{E:3.2}) implies:
$$\int_{B_1(x_0)} |\nabla (1-\Delta)^{-1} Q|^2 \, dx \le \, A
$$
for all $\beta>0$, 
where the constant $A$ does not depend on $x_0\in \R^n$. In particular, 
for $\beta=1$, this inequality together with the mean oscillation 
condition proved in Lemma \ref{Lemma 3.5} (ii), gives  that 
$ \nabla (1-\Delta)^{-1} 
Q \in \text{bmo}(\R^n)$. For $0<\beta<1$, similar estimates yield:
$$||\nabla (1-\Delta)^{-1} Q||_{L^\infty(R^n)} \le B.
$$
Together with the  mean oscillation condition of  Lemma \ref{Lemma 3.5} (ii) 
this gives: 
 $\nabla (1-\Delta)^{-1} Q \in {\rm Lip}_{\frac{1-\beta}{1+\beta}}(\R^n)$.

The following theorem, which also uses  Lemma \ref{Lemma 3.5} in the 
necessity part, treats  
the case $0<\beta\le 1$.

\begin{theorem} \label{Theorem 3.6}
Let $Q \in \mathcal \D'(\R^n)$, $n \ge 2$, and $0<\beta\le 1$. 
Then the  following statements hold. 

{\rm (i)} Suppose that there exist 
$\vec \Gamma \in 
L^2_{\rm loc}(\R^n)^n$ and 
$\gamma \in L^2_{\rm loc}$ such that 
\begin{equation}\label{E:3.23}
Q = \text{{\rm div}} \,\, \vec \Gamma + \gamma,
\end{equation}
where $\vec \Gamma$ and $\gamma$ satisfy respectively the 
inequalities: 
\begin{align}
& \int_{B_\delta(x_0)} \left\vert \vec \Gamma (x) - 
m_{B_\delta(x_0)} ({\vec \Gamma}) \right\vert^2 
\, dx \le 
C_1 \, \delta^{n-2 +\frac 4 {1 +\beta}}, \qquad 0<\delta <\delta_0,
 \label{E:3.24} \\
& \int_{B_\delta(x_0)} \left\vert \gamma (x) \right\vert 
\, dx \le 
C_2 \, \delta^{n-2 + \frac 2 {1 +\beta}}, \qquad 0<\delta <\delta_0, 
\label{E:3.25}
\end{align}
with the constants $C_1, \, C_2$ which do not depend on $x_0$ and $\delta$.
Then there exist $c>0$ and $\epsilon_0>0$ such 
that {\rm (\ref{E:3.1})} holds  for every $\epsilon \in (0, \,
\epsilon_0)$.  

{\rm (ii)} Conversely, 
suppose {\rm (\ref{E:3.1})} holds  for every $\epsilon \in (0, \,
\epsilon_0)$.
Then $Q$ can be represented in the form {\rm (\ref{E:3.23})} 
 where $\vec \Gamma  \in L^2_{\rm loc} (\R^n)^n$,  
$\gamma \in  L^2_{\rm loc} (\R^n)$ satisfy respectively 
conditions {\rm (\ref{E:3.24})} and {\rm (\ref{E:3.25})}.

{\rm (iii)}  If $\epsilon_0=1$, then one can set 
\begin{equation}\label{E:3.26}
\vec \Gamma = 
- \nabla (1 - \Delta)^{-1} \, Q, \qquad \gamma = (1 - \Delta)^{-1} \, Q,
\end{equation}
in {\rm (\ref{E:3.23})}, where 
$\delta_0=1$ 
in {\rm (\ref{E:3.24})}, {\rm (\ref{E:3.25})}.

If $\epsilon_0= +\infty$, then one can set  
\begin{equation}\label{E:3.27}
\vec \Gamma = 
\nabla \Delta^{-1} \, Q, \qquad \gamma = 0, 
\end{equation}
in {\rm (\ref{E:3.23})}, where 
$\delta_0=+\infty$ in 
{\rm (\ref{E:3.24})}, {\rm (\ref{E:3.25})}.
\end{theorem}

\begin{proof} We need only to prove statement (i) since 
(ii) and (iii) follow from Lemma \ref{Lemma 3.5}. Suppose $0< \beta \le 1$ and 
$Q$ is represented in the form {\rm (\ref{E:3.23})} so that 
{\rm (\ref{E:3.24})}, {\rm (\ref{E:3.25})} hold. The proof is 
very similar to the case $\beta>1$ except for the estimate 
 $$|\langle Q u, \, u\rangle |=|\langle \vec \Gamma \, u, \, \nabla u\rangle | \le 
c \,||\nabla u||_{L^2}^{\frac {2 \beta} {1+\beta}} \, 
||u||_{L^2}^{\frac {2} {1+\beta}}, \quad u \in C^\infty_0(\R^n), 
$$
or, equivalently by Corollary~\ref{Corollary 2.3}, 
\begin{equation}\label{E:3.28} 
|\langle Q u, \, u\rangle |=|\langle \vec \Gamma \, u, \, \nabla u\rangle | \le 
c \, \delta^{\frac {2} {1+\beta}} \, 
||\nabla u||_{L^2}^{2} \, 
\quad u \in C^\infty_0(B_\delta(x_0)).  
\end{equation}

We use known characterizations of the Morrey-Campanato spaces. 
In particular, for $0<\beta<1$, condition {\rm (\ref{E:3.24})} 
is equivalent to the condition $\vec \Gamma \in {\rm Lip}_\alpha(\R^n)$, 
where $\alpha = \frac {1-\beta} 
{1+\beta}$.  In other words,
$$|\vec \Gamma (x) - \vec \Gamma (x')| \le C |x-x'|^\alpha, 
\quad |x-x'|\le \delta_0.
$$
Then, for $u \in C^\infty_0(B_\delta(x_0))$,  
$$|\langle \vec \Gamma \, u, \, \nabla u\rangle | =  \left \vert 
\int_{B_\delta(x_0)} (\vec \Gamma - m_{B_\delta(x_0)} (\vec \Gamma)) \cdot 
\nabla u \, u \, dx \right \vert \le c \,   \delta^{\frac {1-\beta} 
{1+\beta}} 
\int_{B_\delta(x_0)} |\nabla u| \, |u| \, dx.
$$
Using Schwarz's  inequality and the estimate 
\begin{equation}\label{E:3.29} 
\int_{B_\delta(x_0)} |u(x)|^p \, dx \le c \, \delta^p \, 
\int_{B_\delta(x_0)} |\nabla u(x)|^p \, dx, \qquad 
 u \in C^\infty_0(B_\delta(x_0)), 
\end{equation}
with $p=2$, we obtain {\rm (\ref{E:3.28})}. 

In the case $\beta=1$, we have 
$\vec \Gamma  \in \text{BMO}(\R^n)$ 
if $\delta_0=+\infty$, or $\vec \Gamma  \in \text{bmo}(\R^n)$ 
 if $\delta_0=1$ respectively. To prove {\rm (\ref{E:3.28})} 
we apply H\"older's inequality with exponents $\frac{2p}{p-2}$, $2$, and $p$
so that $\frac{2p}{p-2} + \frac 1 2 + \frac 1 p =1$, 
where $2< p < \frac {n} {n-2}$, and $n \ge 2$:
\begin{align}
& |\langle \vec \Gamma \, u, \, \nabla u\rangle | =  \left \vert 
\int_{B_\delta(x_0)} (\vec \Gamma - m_{B_\delta(x_0)} (\vec \Gamma)) \cdot 
\nabla u \, u \, dx \right \vert \notag \\ & \le c  \, 
||\vec \Gamma - m_{B_\delta(x_0)} (\vec \Gamma)||_{L^{\frac {2p}{p-2}} 
(B_\delta(x_0))}
||\nabla u||_{L^2(B_\delta(x_0))} \, ||u||_{L^p(B_\delta(x_0))}.
\notag
\end{align}

It is well-known (and is a consequence of 
the John-Nirenberg inequality) that  $\vec \Gamma \in \text{BMO}(\R^n)$ 
if $\delta_0=+\infty$, and  respectively $\vec \Gamma \in \text{bmo}(\R^n)$ 
if $\delta_0=1$, yield 
$$ \int_{B_\delta(x_0)} |\vec \Gamma - m_{B_\delta(x_0)} (\vec \Gamma)|^p 
\, dx  \le c \, \delta^n, \quad 0<\delta \le \delta_0,
$$
for any $1\le p< \infty$. 
Applying this inequality together with {\rm (\ref{E:3.29})}, we obtain 
 {\rm (\ref{E:3.28})}. This completes the proof of Theorem \ref{Theorem 3.6}. 
\end{proof}

It is easy to see that in the case $\beta =1 $ and $\epsilon_0=+\infty$, 
the sufficiency part of Theorem \ref{Theorem 3.6} is equivalent to the inequality
$$
|\langle \vec \Gamma \, u, \, \nabla u\rangle | \le c \, 
||\vec \Gamma||_{\text{BMO}(\R^n)} \, ||u||_{L^2(\R^n)} \, 
||\nabla u||_{L^2(\R^n)}, \quad \forall u \in C^\infty_0(\R^n). 
$$
By duality, the preceding inequality yields:
$$
|| u  \, \nabla  u||_{\mathcal{H}^1(\R^n)} \le c \, ||u||_{L^2(\R^n)} \, 
||\nabla u||_{L^2(\R^n)}, \quad \forall u \in C^\infty_0(\R^n). 
$$

As an immediate consequence, we obtain the vector-valued quadratic form 
inequality mentioned in Sec. \ref{Section 2}:
$$
||(\vec u  \cdot \nabla) \, \vec  u||_{\mathcal{H}^1(\R^n)} \le c \, 
||\vec u ||_{L^2(\R^n)} \, ||\nabla \vec u||_{L^2(\R^n)}, \quad 
 \quad {\rm div} \, 
\vec u = \vec 0,
$$ 
for all $\vec u \in C^\infty_0(\R^n)^n$.

Both of the preceding inequalities  
are  corollaries of the following nonhomogeneous version of the {\it div-curl} lemma  
in the case $p=2$.     
For the sake of completeness, we give a proof 
which is  similar 
to that in 
\cite{CLMS} where it is assumed additionally that 
${\rm div} \, \vec u = \vec 0$.

\begin{lemma} \label{Lemma 3.7} Let $1<p< +\infty$ and $\frac 1 p + 
\frac 1 {p'} = 1$. Let 
$\vec u  \in W^{1,p}(\R^n)^n$ and $v \in W^{1,p'} (\R^n)$. 
Then 
\begin{align}\label{E:3.30}
||{\rm div} \, (\vec u \, v)||_{\mathcal{H}^1(\R^n)}  & \le c \, 
 ( \, ||\vec u ||_{L^p(\R^n)}  ||\nabla v||_{L^{p'}(\R^n)} \notag \\
& +   ||{\rm div} \,  \vec  u||_{L^p (\R^n)} ||v||_{L^{p'}(\R^n)} \, ), 
\end{align}
where $c$ does not depend on  $\vec u$ and $v$.
\end{lemma}

\begin{proof} Without loss of generality we may assume that 
$\vec u  \in C^\infty_0(\R^n)^n$ and $v \in C^\infty_0(\R^n)$. 
 Fix  $\phi \in C^\infty_0(B_1(0))$ and set 
$\phi_{\delta} = \delta^{-n} \phi(\frac{x}{\delta})$, $\delta>0$.  
Let $B=B_{\delta}(x)$, $x\in \R^n$.  
 Then 
$$|| {\rm div} \, (\vec u \, v) ||_{\mathcal{H}^1(\R^n)} \asymp  
\left \Vert\sup_{\delta>0} |{\rm div} (\vec u \, v)  \star \phi_{\delta}| 
\, \right \Vert_{L^1(\R^n)}.
$$
Note that 
\begin{align}
{\rm div} \,  (\vec u \, v) \star \phi_{\delta}(x) 
& = - \delta^{-1-n} \int_B \nabla \phi \,  
 \left (\frac {x-y} \delta\right) \cdot \vec u(y)  \, 
(v(y) - m_B (v)) \,  dy \notag \\
& + \delta^{-n}  m_{B} (v) \int_B  \phi_{\delta} (x-y) \,  
 {\rm div} \, \vec u (y) \, dy.\notag
\end{align}
As in the proof of Lemma II.1 in 
 \cite{CLMS}, pick $\alpha$, $\beta$ so that 
$\frac 1 \alpha - \frac 1 n = 1 - \frac 1 \beta$, and $1\le \alpha < p$, 
$1< \beta < p'$. By H\"older's 
inequality, 
\begin{align}
|{\rm div} (\vec u \, v) \star \phi_{\delta}(x)| & \le c \, 
\left (\delta^{-n} \int_B |\vec u(y)|^\beta  dy \right)^{\frac 1 {\beta}}
\left(  \delta^{-n-\beta'}  \int_B |v(y) - m_B (v)|^{\beta'}  
dy\right)^{\frac 1 {\beta'}}\notag\\
&  + c \, |m_B(v)| \, \delta^{-n} \int_B |{\rm div} \, \vec u(y)| 
\, dy.\notag
\end{align}
By Poincar\'e's inequality,
$$
\left ( \delta^{-n-\beta'}  \int_B |v(y) - m_B (v)|^{\beta'} 
dy\right)^{\frac 1 {\beta'}} \le c \, 
\left ( \delta^{-n}  \int_B |\nabla v|^\alpha\right)^{\frac 1 
{\alpha}}\le c (M |\nabla v|^\alpha (x))^{\frac 1 \alpha},  
$$
where $M$ is the Hardy-Littlewood maximal operator. 
Also, clearly,   
$$
|m_B(v)| \, \delta^{-n} \int_B |{\rm div} \, \vec u(y)| \, dy \le 
M v(x) \, M ({\rm div}\,  \vec u) (x).
$$
Combining these estimates, we get:
$$
|{\rm div} (\vec u \, v) \star \phi_{\delta}(x)| \le c \, 
(M |\vec u|^\beta(x))^{\frac 1 \beta} \, 
(M |\nabla v|^\alpha (x))^{\frac 1 \alpha} 
+ c \, M v(x) \, M ({\rm div}\,  \vec u) (x).
$$
Applying H\"older's inequality together with 
 the  maximal inequality, we obtain (\ref{E:3.30}).
\end{proof}

The following corollary is an immediate consequence of Theorem \ref{Theorem 3.6}, 
Remark 6.2,  
and the  characterizations of the Morrey-Campanato spaces mentioned 
above.

\begin{corollary} \label{Corollary 3.8} {\rm (i)} Under the assumptions of 
Theorem \ref{Theorem 3.6}, in the case  $\beta=1$, condition {\rm (\ref{E:3.24})} 
holds if and only if 
$ \vec \Gamma \in {\rm BMO} (\R^n)$ if $\delta_0=+\infty$, and 
respectively  
$\vec \Gamma\in {\rm bmo} (\R^n)$ if $\delta_0=1$. 

{\rm (ii)} Similarly, in the case $0<\beta<1$,  condition {\rm (\ref{E:3.24})} 
holds if and only if 
 $ \vec \Gamma$ is  in the 
homogeneous Lipschitz space 
${\rm Lip}_{\frac {1-\beta}{1+\beta}} (\R^n)$ if $\delta_0=+\infty$, 
 and respectively 
in the corresponding 
inhomogeneous space ${\rm Lip}_{\frac {1-\beta}{1+\beta}} (\R^n)$ 
 if $\delta_0=1$.
\end{corollary}

\noindent{\bf Remark 6.3.}\label{Remark 6.3} 
In the case $0< \beta\le 1$, 
Theorem \ref{Theorem 3.6} and Corollary \ref{Corollary 3.8} 
can  be restated  directly in terms  
of conditions imposed on $Q$ using the scale of Triebel-Lizorkin spaces  
$F^{p, q}_{\alpha}$ (or Besov spaces $B^{p, q}_{\alpha}$ if $0<\beta<1$) 
with negative index $\alpha$ (see, e.g., \cite{Tri}).

More precisely, it is easy to see that, for $\beta=1$,
  conditions {\rm (\ref{E:3.24})} and {\rm (\ref{E:3.25})} 
 in the homogeneous case are equivalent to 
$Q \in \dot F^{-1, 2}_\infty (\R^n)= \text{BMO}_{-1} (\R^n)$. 
For $0<\beta<1$, an equivalent condition on $Q$ is given by:
$Q \in \dot F^{- \frac{2\beta} {1+\beta}, \infty}_\infty (\R^n)= 
\dot B^{- \frac{2\beta} {1+\beta}, \infty}_\infty(\R^n)$. 
Similarly, in the inhomogeneous case we have $Q \in F^{-1, 2}_\infty(\R^n)$
for $\beta=1$ and $Q \in  
B^{- \frac{2\beta} {1+\beta}, \infty}_\infty(\R^n)$ 
respectively. 

Finally, we are in a position to characterize 
multiplicative quadratic form inequalities of 
Nash's type: 
\begin{equation}\label{E:corr}
\left \vert \langle Q u, \, u\rangle \right \vert  \leq C \, 
||\nabla u||^{2p}_{L^2(\R^n)}  \,||u||^{2(1-p)}_{L^1(\R^n)}, 
\quad \forall u 
\in C^\infty_0(\R^n), 
\end{equation}
where $p  \in (0, \, 1)$.

\begin{corollary} \label{Corollary 3.9} Suppose $Q\in D'(\R^n)$, $n \ge 2$, and 
 $0<p<1$.  Let $p_* = \frac n {n+2}$ and $p^* = \frac {n+1} {n+2}$. 
Then {\rm (\ref{E:corr})} holds if  anyone of the following 
conditions is valid:

\noindent {\rm (a)}  $0< p < p_*$ and $Q=0$.

\noindent {\rm (b)}  $p =  p_*$ and $Q \in L^\infty(\R^n)$.

\noindent {\rm (c)}  $p_*<p<p^*$ and $Q = {\rm div} \, \vec \Gamma$, where 
$\vec \Gamma \in {\rm Lip}_{n+1-p(n+2)}(\R^n)$.

\noindent {\rm (d)} $p =  p^*$ and  $Q = {\rm div} \, \vec \Gamma$, where 
$\vec \Gamma \in {\rm BMO}(\R^n)$. 

\noindent {\rm (e)} $p^*<p<1$ and  $Q = {\rm div} \, \vec \Gamma$, where 
$\vec \Gamma  \in \mathcal{L}^{2, \, \lambda}(\R^n)$, and 
$\lambda= 3n + 2 - 2p(n+2)$.

Conversely, if {\rm (\ref{E:corr})} holds then  conditions {\rm (a)}--{\rm (e)} 
are valid where one may set $\vec \Gamma = \nabla \Delta^{-1} \, Q$.

\end{corollary}

Corollary~\ref{Corollary 3.9} follows by combining  Theorem~\ref{Theorem 3.2} and   
Theorem~\ref{Theorem 3.6} with 
Corollary~\ref{Corollary 2.3} and Remark 3.1.


\begin{thebibliography}{ChWW}


\bibitem[AH]{AH} D.~R. Adams and L.~I. Hedberg,
{\em Function Spaces 
 and Potential Theory}, Springer-Verlag, 
Berlin--Heidelberg--New York, 
 1996. 

\bibitem[Agm]{Agm} S.~Agmon, 
{\em On positive solutions of elliptic equations with periodic coefficients 
in $\R^n$, spectral results and extensions to elliptic operators on Riemannian manifolds}, 
Differential Equations (Birmingham, Ala., 1983), 7--17. 
 North-Holland Math. Stud. {\bf 92},
North-Holland, Amsterdam, 1984. 


\bibitem[Agr]{Agr} M.~S. Agranovich, 
{\em On series in root vectors of operators 
defined by forms with a selfadjoint principal part},
Funct. Anal Appl.,  {\bf 28}  (1994), 151--167.


\bibitem[AiSi]{AiSi} M.~Aizenman and B.~Simon, 
{\em Brownian motion and Harnack inequality for 
Schr\"odinger operators},
Comm. Pure Appl. Math.,  {\bf 35}  (1982), 209--273.

\bibitem[BiSo]{BiSo} 
M.~S. Birman and M.~Z. Solomjak, 
{\em Schr\"odinger operator.  Estimates for  number of
 bound states as function-theoretical problem}, 
Amer. Math. Soc. Transl., Ser. 2, {\bf 150}, 1--54. 
Amer. Math. Soc., Providence, RI, 1992.




\bibitem[BSh]{BSh} F.~A. Berezin and M.~A. Shubin, 
{\em  The Schr\"odinger Equation}, 
Kluwer Academic Publishers, Dordrecht, 1991.

\bibitem[BoBr]{BoBr} J. Bourgain and H. Brezis,
{\em  On the equation ${\rm div}\,Y=f$ 
and application to control of phases},   J. Amer. Math. Soc., {\bf 16}  (2003),   393--426.



\bibitem[BrK]{BrK} 
H. Brezis and T. Kato, 
{\em Remarks on the Schr\"odinger operator with singular complex potentials}, 
 J. Math. Pures Appl., {\bf 58}  (1979), 137--151.


\bibitem[COV]{COV} C. Cascante, J.~M. Ortega, and I.~E. Verbitsky,  
{\em Nonlinear potentials and two weight trace inequalities for
general dyadic  and radial kernels,\/} preprint math.FA/0309286, 
 to appear in Indiana Univ. Math. J.



\bibitem[ChWW]{ChWW} S.-Y.~A. Chang, J.~M. Wilson, and T.~H. Wolff, 
{\em Some weighted norm inequalities concerning the Schr\"odinger operators}, 
Comment. Math. Helv. {\bf 60} (1985), 217--246.

\bibitem[CZh]{CZh} K.~L. Chung and Z. Zhao, 
{\em From Brownian motion to
Schr\"odinger's equation}, Springer--Verlag, 
Berlin--Heidelberg--New York, 1995.  



\bibitem[CLMS]{CLMS} R. Coifman, P.~L. Lions, Y. Meyer
and S. Semmes, 
{\em Compensated compactness and Hardy spaces}, 
J. Math. Pures Appl., {\bf 72} (1993), 247--286.





\bibitem[Co]{Co} P. Constantin, 
{\em Remarks on the Navier-Stokes equations}, 
New Perspectives in Turbulence,  229--261.  
Springer-Verlag, New York, 1991.





\bibitem[Dav1] {Dav1} E.~B. Davies,  
{\em $L^p$ spectral theory of higher order elliptic 
differential operators}, 
Bull. London Math. Soc. {\bf 29} (1997), 513--546.

\bibitem[Dav2] {Dav2} E.~B. Davies,  
{\em Non-self-adjoint differential operators}, 
Bull. London Math. Soc. {\bf 34} (2002), 513--532.


\bibitem[EE] {EE}  D.~E. Edmunds and W.~D. Evans,  
{\em Spectral Theory and Differential Operators}, 
Clarendon Press, Oxford, 1987.


\bibitem[Fef]{Fef} C. Fefferman, 
{\em The uncertainty principle},
Bull. Amer. Math. Soc. {\bf 9}
 (1983), 129--206.


\bibitem[GT] {GT} D. Gilbarg and N.~S. Trudinger, 
{\em Elliptic Partial Differential Equations of Second Order},  
Reprint of the 1998 edition. Classics in Mathematics, 
Springer-Verlag, Berlin--Heidelberg-New York, 2001.


\bibitem[Gr1]{Gr1} E. Grinshpun, 
{\em On spectral properties of Schr\"odinger-type operator with complex
potential},  Oper. Theory: Adv. Appl. {\bf 87}
 (1996), 164--176.

\bibitem[Gr2]{Gr2} E. Grinshpun, 
{\em Asymptotics of spectrum under infinitesimally form-bounded 
perturbation},
Integral Eqs. Oper. Theory {\bf 19}
 (1994), 240--250.

\bibitem[HMV]{HMV} K. Hansson, V.~G. Maz'ya,  and I.~E. Verbitsky, 
{\em Criteria of solvability for
multidimensional Riccati's equations}, 
Arkiv f\"or Matem. {\bf 37} (1999), 87--120.

\bibitem[HW]{HW}
L.~I. Hedberg and Th.~H. Wolff, 
{\em Thin sets in nonlinear
potential theory},
Ann. Inst. Fourier (Grenoble),
{\bf 33} (1983), 161--187.

\bibitem[Ka1] {Ka1}  T. Kato, 
{\em Schr\"odinger operators with singular potentials}, Israel J. Math. {\bf 13} 
(1972), 135--148. 


\bibitem[Ka2] {Ka2}  T. Kato, 
{\em Perturbation Theory for Linear Operators}, 
Springer-Verlag, 
Berlin--Heidelberg--New York, 1995.

\bibitem[KS] {KS}  R. Kerman and E. Sawyer, 
{\em The trace inequality and eigenvalue estimates for Schr\"odinger
operators}, Ann. Inst. Fourier, Grenoble  {\bf 36} (1987), 207--228.


\bibitem[Ko] {Ko}  P. Koosis, 
{\em The Logarithmic Integral I},  
Cambridge Studies in Advanced Mathematics {\bf 12}, 
Cambridge, Cambridge University Press, 1998.


\bibitem[LL] {LL} E.~H. Lieb and M. Loss, 
{\em Analysis}, Second Edition,
Amer. Math. Soc., 
Providence, RI, 2001.

\bibitem[LPS] {LPS}  V.A. Liskevich, M.A. Perelmuter, and Yu. A. 
Semenov, 
{\em Form-bounded perturbations 
of generators of sub-Markovian semigroups}, 
Acta Appl. Math.  {\bf 44} (1996), 353--377.

\bibitem[MM1]  {MM1} A.~S. Marcus and V.~I. Matsaev, 
{\em Operators associated with sesquilinear forms and spectral asymptotics},
 Mat. Issled. {\bf 61} (1981), 86--103.

\bibitem[MM2]  {MM2} A.~S. Marcus and V.~I. Matsaev, 
{\em Comparison theorems for spectra of linear operators and spectral 
asymptotics},
 Trans. Moscow Math. Soc. {\bf 45} (1984), 139--187.


\bibitem[M1] {M1}  V.~G. Maz'ya,
{\em On the theory of the
$n$-dimensional Schr\"odinger operator},
 Izv. Akad. Nauk SSSR,
ser. Matem. {\bf 28} (1964), 1145-1172.




\bibitem[M2] {M2} V.~G. Maz'ya, 
{\em Sobolev Spaces}, 
Springer-Verlag, 
Berlin--Heidelberg--New York, 1985.




\bibitem[MV1]  {MV1} V.~G. Maz'ya and I.~E. Verbitsky, 
{\em Capacitary
estimates for
fractional integrals, with applications to partial differential
equations and Sobolev multipliers},
 Arkiv f\"or Matem. {\bf 33} (1995), 81--115.


\bibitem[MV2]  {MV2} V.~G. Maz'ya and I.~E. Verbitsky, 
{\em The Schr\"odinger operator on the energy space: boundedness
and compactness criteria},  Acta Math. {\bf 188} (2002), 
263--302.

\bibitem[MV3]  {MV3} V.~G. Maz'ya and I.~E. Verbitsky, 
{\em Boundedness and compactness 
criteria for the one-dimensional Schr\"odinger operator}.
 In: Function Spaces, Interpolation Theory and Related Topics.  
Proc. Jaak Peetre Conf., 
Lund, Sweden, August 17-22, 2000. Eds. M. Cwikel, A. Kufner, 
G. Sparr. De Gruyter, Berlin, 2002, 369--382.

\bibitem[MV4]  {MV4} V.~G. Maz'ya and I.~E. Verbitsky, 
{\em The form boundedness
criterion  for the relativistic Schr\"odinger operator}, preprint 
math-ph/0309031, to appear in Ann. Inst. Fourier. 



\bibitem[RSi] {RSi} M. Reed and B. Simon, 
{\em Methods of Modern 
Mathematical 
Physics. II: Fourier Analysis, Self-Adjointness}, 
Academic Press, New York--London, 1975.


\bibitem[RSS] {RSS} G. V. Rozenblum,  M.A. Shubin, and M.Z.  Solomyak, 
{\em Spectral Theory of Differential Operators}, Encyclopaedia 
of Math. Sci., {\bf 64}. Partial Differential Equations VII. 
Ed. M.A. Shubin. Springer-Verlag, Berlin--Heidelberg, 1994.  



\bibitem[Sch] {Sch}  M. Schechter, 
{\em Operator Methods in Quantum Mechanics},  
Dover Publications,
 Mineola, NY, 2002.

\bibitem[Sim] {Sim} B. Simon, 
{\em Schr\"odinger semigroups}, 
Bull. Amer. Math. Soc. {\bf 7} (1982), 447--526. 

\bibitem[St1] {St1} E.~M. Stein,
{\em  Singular Integrals and Differentiability Properties of Functions}, 
Princeton University Press, 
Princeton, New Jersey, 1970. 
93

\bibitem[St2] {St2} E.~M. Stein,
{\em  Harmonic Analysis: Real-Variable 
Methods, Orthogonality, and Oscillatory Integrals}, 
Princeton University Press, 
Princeton, New Jersey, 1993. 


\bibitem[Tri] {Tri} H. Triebel, 
{\em The Structure of Functions.} Monographs in Mathematics, {\bf 97}. 
Birkh\"auser Verlag, Basel, 2001. 


\bibitem[Tru] {Tru} N.~S. Trudinger, 
{\em Linear elliptic operators with measurable coefficients},  
Ann. Scuola Norm. Sup. Pisa {\bf  27} (1973), 265--308.


\bibitem[V1] {V1} I.~E. Verbitsky, 
{\em Imbedding and multiplier theorems for discrete
Littlewood-Paley spaces}, Pacific J. Math., 
 {\bf 176} (1996),  
 529--556.

\bibitem[V2] {V2} I.~E. Verbitsky, 
{\em Nonlinear potentials and trace 
inequalities}, 
The Maz'ya Anniversary
Collection. Eds. J. Rossmann,  P. Tak\'ac, and G. Wildenhain.  
Operator Theory: Advances and Applications {\bf 110} (1999),  
 323--343.

\end{thebibliography}
\end{document}